\documentclass[11pt]{article}
\usepackage[centertags,intlimits]{amsmath}                     
\usepackage{amsfonts}                     
\usepackage{amsthm}
\usepackage{graphics}
\usepackage[dcucite]{harvard}
\allowdisplaybreaks[2]
\numberwithin{equation}{section}
\newtheorem{theorem}{Theorem}[section]
\newtheorem{lemma}{Lemma}[section]
\newtheorem{corollary}{Corollary}[section]

\theoremstyle{remark}
\newtheorem{remark}{Remark}[section]

\providecommand{\abs}[1]{\lvert #1\rvert}
\providecommand{\norm}[1]{\lVert #1\rVert}

\newcommand{\nc}{\newcommand}
\nc{\vb}{\mathbf{v}}
\nc{\bx}{\mathbf{x}}
\nc{\by}{\mathbf{y}}
\nc{\bz}{\mathbf{z}}
\nc{\bu}{\mathbf{u}}
\nc{\bfr}{\mathbf{r}}
\nc{\bA}{\mathbf{A}}
\nc{\R}{\mathbb R}
\nc{\N}{\mathbb N}
\nc{\C}{\mathbb C}
\nc{\D}{\mathbb D}
\nc{\F}{\mathbf F}
\nc{\B}{\cal B}
\nc{\br}{\bigr}
\nc{\bl}{\bigl}
\nc{\Bl}{\Bigl}
\nc{\Br}{\Bigr}
\nc{\ind}{\mathbf{1}}
\DeclareMathOperator*{\esssup}{ess\,sup\;}  

\title{Stochastic processes in random graphs}
\author{Anatolii A. Puhalskii\\
University of Colorado, Denver, U.S.A. \\
 Institute for Problems in Information
Transmission, Moscow, Russia}
\date{January 26, 2004}
\begin{document}
\sloppy
\maketitle
\begin{abstract}  
We study the asymptotics of large, moderate and normal deviations for
the connected components of the sparse random graph by {\em the method
  of stochastic processes}. We obtain the
logarithmic asymptotics of large deviations 
of the joint distribution of the number of 
connected components, of  the sizes of the giant components, and of the
numbers of the excess edges of the giant components. For the
supercritical case, we obtain
 the asymptotics of normal deviations and 
the logarithmic asymptotics of  large and moderate deviations
of the joint distribution of the
number of components, of the size of the largest component, and of the
number of the excess edges of the largest component. 
For the critical case, we obtain the
logarithmic asymptotics of moderate deviations 
of the joint distribution of the sizes of  connected components and of  the
numbers of the excess edges. 
Some  related asymptotics are also
 established. 
The proofs of the large and moderate
deviation asymptotics employ methods of idempotent probability theory.
As a byproduct of the results, we provide some additional insight into the
nature of phase  transitions in sparse random graphs.
\end{abstract}  
\footnotetext[1]{
\indent
{\em Short title:} Stochastic processes in random graphs
}

\thispagestyle{empty}

\section{Introduction}
\label{sec:introduction}
The random graph $\mathcal{G}(n,p)$ 
is defined as a non-directed graph on $n$ vertices where
every two vertices are independently connected by an edge with
probability $p$. The graph is said to be sparse if $p=c/n$ for $c>0$
and $n$ large. 
Properties of  sparse random graphs
have been studied at length and major developments have been
  summarised in the recent
monographs by Bollobas \citeyear{Bol01}, 
 Janson, {\L}uczak and
Ruci{\'n}ski \citeyear{JanLucRuc00}, and Kolchin \citeyear{Kol99}. 
The focus of this paper is on the
asymptotics as $n\to\infty$
of the sizes of the giant connected components, i.e., 
components  of  order  $n$ in size,  
of $\mathcal{G}(n,c_n/n)$, where $c_n\to c>0$. 
It is known that for $c>1$ with probability tending to $1$ as
$n\to\infty$ there exists a unique giant
component of $\mathcal{G}(n,c/n)$, which is asymptotically 
 $\beta n$ in size, where $\beta\in(0,1)$ is the positive
root to the equation $1-\beta=\exp(-\beta c)$, 
the rest of the
components being of sizes not greater than of  order  $\log n$.
For $c<1$ with probability tending to $1$ there are no connected components 
of sizes greater than of  order  $\log n$, while for $c=1$ the size of
the largest component is of order $n^{2/3}$. 
Our primary objective is to evaluate the probabilities that there exist 
several giant connected components. As to be expected,
these probabilities are 
 exponentially small in $n$, so we study the decay rates and state
our results in the form of the large deviation principle (LDP). 
In addition, influenced by the papers of  Stepanov
\citeyear{Ste70} and Aldous \citeyear{Ald97}, we concern ourselves with the 
large deviation asymptotics of the  number of the 
 connected components  and of the numbers
of the  excess edges of the connected components.
Thus, the main result is an LDP for the joint
distribution of the normalised number of the connected components of
$\mathcal{G}(n,c_n/n)$,  
of the normalised sizes   of  the connected components, and 
 of the normalised numbers of the excess edges.
Projecting  yields LDPs for  the sizes and
for the number      of the connected components. 
Stepanov \citeyear{Ste70} and later Bollob{\'a}s, Grimmett, and Janson 
\citeyear{BolGriJan96}, analysing a more general setting, have obtained 
 the logarithmic asymptotics
 of the moment generating function of the number of the connected
 components of $\mathcal{G}(n,c/n)$. 
If $c\le2$, the latter asymptotics also yield the LDP for the number of
components, as  Bollob{\'a}s, Grimmett, and Janson 
\citeyear{BolGriJan96} demonstrate,
 but not  for arbitrary $c>0$.
This anomaly  is caused by a phase transition occurring at $c=2$
 discovered by Stepanov \citeyear{Ste70}, which results, as we show, in
the action functional becoming non-convex as $c$ passes
through the value of $2$. Moreover, the
 phase transition turns out to consist
in  a giant component breaking up.

Another group  of results presented in the paper has to do with
the properties of the  largest connected component. 
We establish  normal
deviation, moderate deviation, and large deviation
 asymptotics for the joint distribution of the size of the  largest
 connected component, of
the number of its excess edges, and of 
the number of the connected components.
In  related work,  O'Connell \citeyear{Oco98} proves
 an LDP for the size of the largest connected component 
of $\mathcal{G}(n,c/n)$ and  Stepanov
\citeyear{Ste70a,Ste72}
obtains  central limit theorems for the size of the largest
component  and the
number of components; different proofs of the 
central limit theorem for the size of the largest
component  are  given in Pittel \citeyear{Pit90} and 
Barraez,  Boucheron, and Fernandez de la Vega \citeyear{BarBouFer00}, the
latter authors also provide estimates of the rate of convergence.
Our third group of results concerns the critical random graph when $c=1$.
We complement the result of Aldous \citeyear{Ald97} on the convergence in
distribution of the suitably normalised sizes and  numbers of the excess
edges of the connected components with
moderate deviation asymptotics for these random variables.

Our analysis employs   a surprising (to us) connection to queueing
 theory.   The results outlined above
are derived as consequences of the asymptotic properties of
 a ``master'' stochastic process, which captures the partitioning of the random
 graph into   connected components and builds on 
an earlier construction of a similar sort, cf., Janson,
{\L}uczak and Ruci\'nski \citeyear{JanLucRuc00}. This
 stochastic process is intimately related to the waiting-time
 process (or the queue-length process)
in a certain 
time- and state-dependent queueing system
 and  the connected components correspond to the busy cycles of the
 system. 
We capitalise on this connection by invoking
 our  intuition for the behaviour of queues as well as some standard
queueing theory tools such as properties
 of the Skorohod reflection mapping.
Thus,  at first we apply  the methods of the asymptotic theory of
stochastic processes, namely, the methods of weak and 
large deviation convergence, in order 
to establish asymptotics of the master
process and  then translate them into the properties of
the connected components of the random graph.
 In the context of the random graph theory,
the present paper can thus  be
 considered as developing the approach  pioneered by Aldous
 \citeyear{Ald97} of deriving asymptotic properties of 
random graphs  as  consequences of asymptotics of  associated stochastic
processes. On the technical side, we extensively use the observation
also made by Aldous  \citeyear{Ald97}
 that the connected components can be identified with the excursions
 of a certain  stochastic process.
Yet, the specific  construction in this paper is
different from the one of Aldous \citeyear{Ald97}.
It is actually   much the same as
 the one of Barraez,  Boucheron, and Fernandez de la Vega
 \citeyear{BarBouFer00}, as we have learnt after the paper had
 been submitted,  except for an important distinction, which we
 discuss below. 

There are also other interesting technical aspects of the proofs,
which concern all three types of asymptotics: large deviations,
moderate deviations and normal deviations.
The proof of the LDP for the master
 process  relies on the results of
the large deviation theory of semimartingales, %
Puhalskii \citeyear{Puh01}, which seem to be called for since the
action functional  is  ``non-Markovian'' and ``non-time homogeneous''. 
 The cumulant that characterises the action functional
is not non-degenerate, which   is known to present
 certain difficulties for establishing the LDP.
In the standard approach the problem reveals itself when
 the large deviation  lower bound  is proved and  is
 usually tackled via a perturbation argument: an extra  term is added
 to the process under study    so that the perturbed process has a
 non-degenerate cumulant and then a limit is taken in the lower bound for
 the perturbed process as the perturbation term tends to zero, cf.,
 Liptser \citeyear{Lip96},  de Acosta \citeyear{Aco00}, 
Liptser, Spokoiny and Veretennikov
 \citeyear{LipSpoVer02}.  
 Our approach to proving the LDP
 replaces establishing  the upper and
 lower bounds with the requirement that the limiting maxingale problem have a
 unique solution. The degeneracy of the cumulant presents
 a problem here too. We  cope with it via a perturbation argument as well
 the important difference being that the perturbation is applied to the
 limit idempotent process that specifies the maxingale problem
 rather than to the pre-limiting stochastic processes.
 This  change of the object has important methodological advantages.
 Firstly, the  proof of the LDP is simplified as compared
with the case where the perturbation is introduced at the
 pre-limiting stage. Secondly, 
once the perturbation argument has been carried out for a given cumulant,
one can use it to prove LDPs for a range of stochastic processes that
produce the same cumulant in the limit. We expand on these ideas
in Puhalskii \citeyear{Puh03un}.
The actual implementation of the perturbation approach for the setting
in the paper
relies on the techniques  of idempotent probability theory, Puhalskii
\citeyear{Puh01},  and also
draws on time-change arguments 
in Ethier and Kurtz \citeyear[Chapter 6]{EthKur86},
thus applying probabilistic ideas to an idempotent
 probability setting.
Idempotent probability theory techniques are also instrumental in the proofs
of the moderate-deviation asymptotics. These proofs are modelled on the
preceding proofs of the normal-deviation asymptotics and
to a large degree replicate them by replacing limit stochastic
processes with their idempotent counterparts.

An interesting feature of the proof of the normal deviation
asymptotics for  the largest component 
is that it provides an instance of
 convergence in distribution of stochastic processes 
``with unmatched jumps in the limit
 process'',  Whitt \citeyear{Whi02}, i.e., though the jumps of the
 pre-limiting  processes vanish, the limit process is discontinuous, moreover,
it is not right-continuous with left-hand limits. We thus do not have
convergence in distribution in the Skorohod topology and have to use
some ad-hoc techniques to obtain the needed conclusions.
 As it is explained in Whitt \citeyear{Whi02}, convergence with unmatched
 jumps often occurs in the study of diffusion approximation of time
dependent queues, so it is not surprising (but is amusing) to see it here.
Incidentally,  we
are faced with  a similar situation 
in the proof of the moderate-deviation asymptotics
when no LDP for the Skorohod topology
is available and the corresponding limit theorem can be viewed as an
example of large deviation convergence in distribution of stochastic
processes with unmatched jumps in the limit idempotent process.

We now outline the structure of the
paper. In Section~\ref{sec:main-results}
 we define the underlying
stochastic  processes,
derive queue-like  equations for them, 
state the  results on the properties of the connected components, 
and comment on
them. Section~\ref{sec:techn-prel} contains  technical preliminaries.
 Section~\ref{sec:large-devi-asympt} is concerned
with proving the LDP for the basic processes. In Section~\ref{sec:large-con} 
 the LDPs for the connected
components are proved. Section~\ref{sec:norm-moder-devi} contains  proofs of
the normal and moderate deviation
 asymptotics for the largest component. 
Section~\ref{sec:crigra} considers critical random
 graphs. 
  The appendix provides
an overview of the  notions and facts of idempotent probability theory
invoked in the proofs.

\section{The model equations and main results}
\label{sec:main-results}

  We model the formation 
of the sparse random graph on $n$
  vertices with edge probability $p_n=c_n/n$ via 
  stochastic processes $V^n=(V^n_i,\,i=0,1,\ldots,n)$
  and $E^n=(E^n_i,\,i=0,1,\ldots,n)$. At time $0$ the processes are
  at $0$. At time 1 an arbitrary
  vertex of the graph is picked and is connected by  edges
   to the other vertices
   independently with probability $p_n$. We say that this
   vertex has been first {\em generated} and then 
 {\em saturated}. The vertices, to which it has been
  connected, are called 
  {\em generated}. The value of $V^n_1$ is defined as
the number of vertices in the
  resulting connected component, i.e., the number of the
  generated vertices at time 1; $E^n_1=0$. 
At time 2 we pick one of the generated
  non-saturated vertices if any and saturate it  by connecting 
  it independently with probability $p_n$ to the vertices that either have not
  been generated yet or have been generated but not saturated.
If there are no generated non-saturated
  vertices, we pick an arbitrary non-generated vertex, declare it
  generated and saturate it by attempting to connect it to the
  non-generated vertices, thus generating those of these vertices 
  connection to which is established. We denote as
 $V^n_2$   the total number of vertices
  generated at times 1 and 2 and we denote as $E^n_2$ the number of
  edges connecting the vertex that was saturated at time 2 with
  the other vertices generated at time 1 if any. 
  We proceed in this fashion by saturating
  one vertex per unit of time until
  time $n$. Thus, at time $i$ a generated non-saturated vertex is
  picked and is connected by edges with probability $p_n$ to the
  non-saturated vertices, both  generated 
  and not yet generated; if there are no generated non-saturated
  vertices available, then an arbitrary non-generated vertex is
  chosen, is declared generated and is then saturated. The 
  increment $V^n_i-V^n_{i-1}$ is defined as
  the number of vertices generated at $i$, the 
  increment $E^n_i-E^n_{i-1}$ is defined as
  the number of edges drawn at $i$ between the vertex being saturated
  and the vertices generated by $i$. 
Thus, $V^n_i-V^n_{i-1}$ equals either the number of new vertices
  joined to a connected component at time $i$ if $V^n_{i-1}>i-1$ 
 or it is the number
  of vertices that start a new component  at $i$ if
  $V^n_{i-1}=i-1$. Accordingly, the increment $E^n_i-E^n_{i-1}$ 
either equals the number of   excess edges in a connected component
appeared at time $i$, i.e., the edges in excess of those 
that are necessary to maintain connectedness, or $E^n_i-E^n_{i-1}=0$.  
Since during this process every two vertices independently 
attempt connection with probability $p_n$ exactly once, 
the resulting configuration of edges at time $n$ 
has the same distribution as the one in the random graph $\mathcal{G}(n,p_n)$. 
 In fact, the sizes of the  connected components of
 $\mathcal{G}(n,p_n)$ can be  recovered from the process $V^n$ as 
 time-spans  between successive moments when $V^n_i$ is equal to $i$.
 The numbers of the excess edges in the
connected components are equal to the increments of the process 
$E^n$ 
over such time periods.  In addition, the number of times when
$V^n_i$ is equal to $i\in\{1,2,\ldots,n\}$ equals the number of the connected
  components of $\mathcal{G}(n,p_n)$. We now turn this description
  into equations.

 Since at time $i$ there are  $V^n_i$
    generated  vertices,
the evolution of $V^n$ is given by the following recursion
  \begin{multline}
    \label{eq:1}
    V^n_{i}=
\bl(V^n_{i-1}+\sum_{j=1}^{n-V^n_{i-1}}\xi^n_{ij}\br)
\ind(V^n_{i-1}>i-1)+\bl(i+\sum_{j=1}^{n-i}\xi^n_{ij}\br)
\ind(V^n_{i-1}=i-1),\\i=1,2,\ldots,n,\,
 V^n_0=0,
  \end{multline}
where 
the $\xi^n_{ij},\,i\in\N,\,j\in\N,\,n\in\N$,
 are mutually independent Bernoulli random variables with
$\mathbf{P}(\xi^n_{ij}=1)=p_n$ and 
$\ind(\Gamma)$ is the indicator function of an event
$\Gamma$ that equals $1$ on $\Gamma$ and $0$ outside of $\Gamma$. 
Let $Q^n_i$ denote the number of non-saturated 
generated vertices at time $i$. Since $Q^n_i=V^n_i-i$, 
\eqref{eq:1} implies that
\begin{multline}
  \label{eq:10}
  Q^n_{i}=\bl(Q^n_{i-1}+\sum_{j=1}^{n-Q^n_{i-1}-(i-1)}\xi^n_{ij}-1\br)
\ind(Q^n_{i-1}>0)+\sum_{j=1}^{n-i}\xi^n_{ij}
\,\ind(Q^n_{i-1}=0),\\i=1,2,\ldots,n,\,
 Q^n_0=0.
\end{multline}
The evolution  of the process $E^n$ is governed by the recursion
\begin{equation}
  \label{eq:25}
  E^n_i=E^n_{i-1}+\sum_{j=1}^{Q^n_{i-1}-1}\zeta^n_{ij},\,i=1,2,\ldots,n,\,
E^n_0=0,
\end{equation}
where the $\zeta^n_{ij},\,i\in\N,\,j\in\N,\,n\in\N$,
are mutually independent Bernoulli random variables with
$\mathbf{P}(\zeta^n_{ij}=1)=p_n$, 
which are independent of the $\xi^n_{ij}$, and
 sums  are assumed to be equal to $0$ if the upper summation
index is less than the lower one.

We use for the analysis of  \eqref{eq:10} the following
insight. 
Let us introduce a related process 
${Q'}^n=({Q'}^n_i,\,i=0,1,\ldots, n)$ by 
\begin{equation}
  \label{eq:11}
  {Q'}^n_{i}=\bl({Q'}^n_{i-1}
+\sum_{j=1}^{n-{Q'}^n_{i-1}-i}\xi^n_{ij}-1\br)^+,
\,i=1,2,\ldots,n,\,
 {Q'}^n_0=0,
\end{equation}
where $a^+=\max(a,0)$. We note that ${Q'}^n_i$
is the waiting time of the $i$-th request, where $i=0,1,\ldots,n-1$,
in the queueing system that
starts empty, has  $\sum_{j=1}^{n-{Q'}^n_{i}-(i+1)}\xi^n_{i+1,j}$ 
as the $i$-th request's service time 
and $1$ as the interarrival times. (Alternatively, ${Q'}^n_i$
can be considered as the queue length at time $i$ 
 for the discrete-time queueing system  that
serves one request per unit time,  the number of arrivals in
$[i,i+1]$ being equal to $\sum_{j=1}^{n-{Q'}^n_{i}-(i+1)}\xi^n_{i+1,j}$.)
It is seen  that ${Q'}^n_i=(Q^n_i-1)^+$,
so the asymptotic properties of the process 
$Q^n=(Q^n_i,\,i=0,1,\ldots,n)$  multiplied by a
vanishing constant are the same as those of the process
${Q'}^n=({Q'}^n_i,\,i=0,1,\ldots,n)$. 
In addition, 
 connected components of the random graph correspond to  busy
cycles of this queueing system, i.e., the excursions of ${Q'}^n$.
Thus, a possible way to study the random graph is through
the process ${Q'}^n$. 
This  approach is, in effect, pursued by
Barraez,  Boucheron, and Fernandez de la Vega \citeyear{BarBouFer00}
who study what  in our notation is the process
$({Q'}^n_i+1,\,i=0,1,\ldots,n)$. 
It is, however, inconvenient for our purposes 
because ${Q'}^n_i=0$ not only when $Q^n_i=0$ but also when
$Q^n_i=1$, so the queueing system may have more busy cycles than there
are connected components.  For this reason, we
choose to work with $Q^n$ directly.
Yet, the queueing theory connection serves us as a guide.
Let us recall that the solution of \eqref{eq:11} is  given by 
$  {Q'}^n=\mathcal{R}(\tilde{S}^n),$
where the process $\tilde{S}^n=(\tilde{S}^n_i,
\,i=0,1,\ldots,n)$ with $\tilde{S}^n_0=0$ 
is defined  by
$  \tilde{S}^n_i
=\sum_{k=1}^{i}\sum_{j=1}^{n-{Q'}^n_{k-1}-k}\xi^n_{kj}-i$, 
and $\mathcal{R}$ is the Skorohod reflection operator:
$  \mathcal{R}(\mathbf{x})_t=\mathbf{x}_t-\inf_{s\in[0,t]}\mathbf{x}_s\wedge 0
$ for $\bx=(\bx_t,\,t\in\R_+)$, where $\wedge$ denotes the minimum.
We find it productive to express $Q^n$ as a reflection too.
The idea is to sacrifice the Markovian character
of recursion  \eqref{eq:10} for the nice  properties
of the  reflection mapping. 

A manipulation of \eqref{eq:10} yields the following equality:
\begin{equation}
  \label{eq:2}
  Q^n_{i}=S^n_i+\epsilon^n_{i}+\Phi^n_i,\,i=0,1,\ldots,n,\,
\end{equation}
where
\begin{align}
  \label{eq:3}
  S^n_i&=\sum_{k=1}^i\Bl(\sum_{j=1}^{n-Q^n_{k-1}-(k-1)}\xi^n_{kj}-1\Br),\\
  \label{eq:35}
  \epsilon^n_{i}&=\ind(Q^n_i>0)
-\sum_{k=1}^i\xi^n_{k,n-k+1}\,\ind(Q^n_{k-1}=0),\\  
  \label{eq:36}
  \Phi^n_i&=\sum_{k=1}^i\ind(Q^n_{k}=0).
\end{align}
For the sequel it is useful to note that $\Phi^n_n$ equals the number of
the connected components of $\mathcal{G}(n,p_n)$. 

Denoting 
as $\lfloor x\rfloor$ the integer part of $x\in\R_+$,  we
introduce  continuous-time processes
$\overline{Q}^n=(\overline{Q}^n_t,
\, t\in[0,1])$, $\overline{S}^n=(\overline{S}^n_t,\, t\in[0,1])$, 
 $\overline{\Phi}^n=(\overline{\Phi}^n_t,\, t\in[0,1])$,
and  $\overline{E}^n=(\overline{E}^n_t,\, t\in[0,1])$  
by the respective equalities 
$\overline{Q}^n_t=Q^n_{\lfloor  nt\rfloor}/n$,
$\overline{S}^n_t=S^n_{\lfloor  nt\rfloor}/n$, 
$\overline{\Phi}^n_t=\Phi^n_{\lfloor  nt\rfloor}/n$, and
$\overline{E}^n_t=E^n_{\lfloor  nt\rfloor}/n$. By \eqref{eq:36}
$\overline{\Phi}^n_t=\int_0^t\ind(\overline{Q}^n_s=0)\,d\overline{\Phi}^n_s$,
so, by \eqref{eq:2}
 the pair $(\overline{Q}^n,\overline{\Phi}^n)$  solves the Skorohod
 problem in $\R$ for $\overline{S}^n+\overline{\epsilon}^n$, consequently,
\begin{align}
  \label{eq:33}
\overline{Q}^n&=\mathcal{R}(\overline{S}^n+\overline{\epsilon}^n),\\
  \label{eq:33a}
\overline{\Phi}^n&=\mathcal{T}(\overline{S}^n+\overline{\epsilon}^n),
  \end{align}
where $\mathcal{T}(\bx)_t=-\inf_{s\in[0,t]}\bx_s\wedge0$ for 
$\bx=(\bx_t,t\in\R_+)$
and
  $\overline{\epsilon}^n=(\overline{\epsilon}^n_t,\,t\in[0,1])$ is
  defined by
  \begin{equation}
    \label{eq:1a}
\overline{\epsilon}^n_t=\frac{\epsilon^n_{\lfloor nt\rfloor}}{n} .
  \end{equation}
Equation \eqref{eq:25} yields the representation
\begin{equation}
  \label{eq:27}
  \overline{E}^n_t=\frac{1}{n}\sum_{i=1}^{\lfloor nt\rfloor}
\sum_{j=1}^{Q^n_{i-1}-1}\zeta^n_{ij},\,t\in[0,1].
\end{equation}
Equations \eqref{eq:3}, \eqref{eq:33}, \eqref{eq:33a}, and \eqref{eq:27} play a
central part in establishing the main  results of the paper. In some more
detail, the processes $\overline{\epsilon}^n$ prove to be 
inconsequential and may
be disregarded (see Lemma~\ref{le:eps}), so \eqref{eq:3}, (\ref{eq:27}), and
\eqref{eq:33} enable us to obtain functional limit theorems for the
processes $(\overline{S}^n,\overline{E}^n)$, which on making another
use of (\ref{eq:33}) and (\ref{eq:33a}) yield the
asymptotics of the connected components (we note that the latter step
does not reduce to a mere application of the continuous mapping principle).
Before  embarking on this programme, we  state and discuss the results.

We will say that  a sequence
$ \mathbf{P}_n,\, n\in\N,$ of probability measures on the Borel 
$\sigma$-algebra of a metric space $\Upsilon$ (or a sequence of
random elements $ X_n,\, n\in\N,$ with values in $\Upsilon$ 
and distributions $\mathbf{P}_n$) obeys the large deviation principle (LDP)
for scale $k_n$, where $k_n\to\infty$ as $n\to\infty$,
with action functional $I:\,\Upsilon\to[0,\infty]$  if
 the sets $ \{ \upsilon \in \Upsilon : ~ I(\upsilon)  \leq a \}$ are compact 
for all $a\in\R_+$, 
\begin{align*}
\limsup_{n\to\infty}\frac{1}{k_n}\log \mathbf{P}_n (F)  
\leq - \inf_{\upsilon \in F}
I(\upsilon)\quad&\text{ for all closed sets $F  \subset \Upsilon$}\\
\intertext{ and}
\liminf_{n\to\infty}\frac{1}{k_n}\log \mathbf{P}_n (G)  \geq  -  \inf_{\upsilon
  \in G} 
I(\upsilon)\quad&\text{ for all open sets $G  \subset \Upsilon$}.  
\end{align*}
Let for $u\in[0,1]$, $\rho\in\R_+$, and $c>0$
\begin{align}
\label{eq:k}
  K_\rho(u)&=
u\log\frac{\rho u}{1-e^{-\rho u}}-\frac{\rho u^2}{2},
\\
\label{eq:l}
L_c(u)&=(1-u)\log(1-u)+(c-\log c)u-\frac{c u^2}{2},
\end{align}
where we adopt the conventions $0/0=1$ and  $0\cdot\infty=0$.
We also denote $a\vee b=\max(a,b)$,
$\pi(x)=x\log x-x+1,\,x\in\R_+$, and assume that $\pi(\infty)\cdot 0=\infty$.

Let $\mathbb{S}$ denote the subset of $\R_+^\N$ of sequences
$\bu=(u_1,u_2,\ldots)$ such that $\sum_{i=1}^\infty u_i<\infty$.
Given a convex function $\chi:\,\R_+\to\R_+$
such that $\chi(0)=0$, $\chi(x)>0$ for $x>0$, 
and $\chi(x)/x\to0$ as $x\to0$, we endow $\mathbb{S}$ with an  Orlicz
space topology that is  generated by 
a Luxembourg metric $d_\chi(\bu,\bu')=\inf\{b\in\R_+:\,
\sum_{i=1}^\infty \chi(\abs{u_i-u'_i}/b)\le 1\}$, where
$\bu=(u_1,u_2,\ldots)$ and $\bu'=(u'_1,u'_2,\ldots)$
 (cf., e.g., Krasnosel'skii and Rutickii \citeyear{KraRut}, 
Bennett and Sharpley \citeyear{BenSha88}). 
 Let also $\mathbb{S}_1$ denote the subspace of $\mathbb{S}$
of non-increasing sequences  $\bu=(u_1,u_2,\ldots)$ 
with $\sum_{i=1}^\infty u_i\le 1$.  
It is endowed with induced topology which 
is  equivalent to the  product topology.

Let
$(U^n_1,U^n_2,\ldots,)$ be the  sequence of
 the sizes of the  connected components of  the random
graph $\mathcal{G}(n,c_n/n)$ 
arranged in descending order appended with zeros to
make it infinite, $(R^n_1,R^n_2,\ldots)$ be the sequence of the corresponding
 numbers of the excess edges appended with zeros,  and
 $\alpha^n$ be the number 
of the connected components of $\mathcal{G}(n,c_n/n)$.
  We define $\overline{U}^n=(U^n_1/n,U^n_2/n,\ldots)$ and
  $\overline{R}^n=(R^n_1/n,R^n_2/n,\ldots)$,  and consider 
 $(\alpha^n/n,\,\overline{U}^n,\,\overline{R}^n)$  as a random
element of $[0,1]\times \mathbb{S}_1\times\mathbb{S}$, which is assumed to
be equipped with product topology.

\begin{theorem}
  \label{the:jointld}
Let $c_n\to c>0$ as $n\to\infty$. 
Then the sequence $(\alpha^n/n,\,\overline{U}^n,\,\overline{R}^n),
\, n\in\N,$ obeys the LDP
  in $[0,1]\times \mathbb{S}_1\times\mathbb{S}$ 
for scale $n$ with action functional
  $I_c^{\alpha,U,R}$ defined for $a\in[0,1]$, $\bu=(u_1,u_2,\ldots)\in
  \mathbb{S}_1$, and $\bfr=(r_1,r_2,\ldots)\in
  \mathbb{S}$ by
  \begin{multline*}
              I_c^{\alpha,U,R}(a,\bu,\bfr)=
\sum_{i=1}^\infty\sup_{\rho\in\R_+}\Bl( K_{\rho}(u_i)+
r_i\log\frac{\rho}{c}\Br)
+L_c\bl((1-2a)\vee\sum_{i=1}^\infty  u_i\br)\\
+\frac{c}{2}\bl(1-(1-2a)\vee\sum_{i=1}^\infty u_i\br)^2
\;\pi\biggl(\frac{2\bl(1-a-(1-2a)\vee\sum_{i=1}^\infty u_i\br)}{%
 c\bl(1-(1-2a)\vee\sum_{i=1}^\infty u_i\br)^2}\biggr)
\end{multline*}
if $\sum_{i=1}^\infty u_i\le 1-a$, 
and $I_c^{\alpha,U,R}(a,\bu,\bfr)=\infty$ otherwise.
\end{theorem}
As a consequence, we obtain some marginal LDPs.
\begin{corollary}
  \label{the:4} 
Let $c_n\to c>0$ as $n\to\infty$.
Then the sequences $(\overline{U}^n,\overline{R}^n),\, n\in\N$, 
and $(\alpha^n/n,\overline{U}^n),\, n\in\N,$  obey the LDPs
   for scale $n$ in the respective spaces 
 $\mathbb{S}_1\times\mathbb{S}$ and  $[0,1]\times\mathbb{S}_1$
with respective action functionals
  $I_c^{U,R}$ and $I_c^{\alpha,U}$
 defined for $\bu=(u_1,u_2,\ldots)\in \mathbb{S}_1$,
$\bfr=(r_1,r_2,\ldots)\in \mathbb{S}$, and $a\in[0,1]$ by
\begin{equation*}
  I_c^{U,R}(\bu,\bfr)=
\sum_{i=1}^\infty\sup_{\rho\in\R_+}\Bl( K_{\rho}(u_i)+
r_i\log\frac{\rho}{c}\Br)
+ L_c\Bl(\Bl(1-\frac{1}{c}\Br)\vee\sum_{i=1}^\infty
u_i\Br),
\end{equation*}
and
  \begin{multline*}
              I_c^{\alpha,U}(a,\bu)=
\sum_{i=1}^\infty K_{c}(u_i)
+L_c\bl((1-2a)\vee\sum_{i=1}^\infty  u_i\br)\\
+\frac{c}{2}\bl(1-(1-2a)\vee\sum_{i=1}^\infty u_i\br)^2
\;\pi\biggl(\frac{2\bl(1-a-(1-2a)\vee\sum_{i=1}^\infty u_i\br)}{%
 c\bl(1-(1-2a)\vee\sum_{i=1}^\infty u_i\br)^2}\biggr)
\end{multline*}
if $\sum_{i=1}^\infty u_i\le 1-a$,
and $I_c^{\alpha,U}(a,\bu)=\infty$ otherwise.
\end{corollary}
  \begin{corollary}
    \label{the:5}
Let $c_n\to c>0$ as $n\to\infty$.
Then the sequences $\alpha^n/n,\, n\in\N,$ and 
$\overline{U}^n/n,\, n\in\N,$ obey the LDPs in the respective spaces
$[0,1]$ and $\mathbb{S}_1$ for scale $n$ with the respective action
functionals 
\begin{equation*}
  I_c^\alpha(a)=\inf_{\tau\in[(1-2a)^+,1-a]}
      \biggl( K_c(\tau)+L_c(\tau)+
\frac{c(1-\tau)^2}{2}
\;\pi\Bl(\frac{2(1-a-\tau)}{c(1-\tau)^2}\Br)\biggr)
\end{equation*}
and 
\begin{align*}
  I_c^U(\bu)=&
      \sum_{i=1}^\infty K_c(u_i)+L_c\Bl(\Bl(1-\frac{1}{c}\Br)
\vee\sum_{i=1}^\infty u_i\Br).
  \end{align*}
\end{corollary}
The next corollary clarifies the structure of the most probable
  configurations of the giant components.
Let, given  $\delta>0$,  $m\in\mathbb{N}$, and
$u_i\in(0,1],\,i=1,2,\ldots,m$ with $\sum_{i=1}^mu_i\le 1$,
$A_{\delta}^n(u_1,\ldots,u_m)$
denote the event that there exist  
$m$ connected components of $\mathcal{G}(n,c_n/n)$,
whose respective sizes are
between $n(u_i-\delta)$ and $n(u_i+\delta)$ for
$i=1,2,\ldots,m$. For $\epsilon>0$, we define event
$\tilde{A}_{\delta,\epsilon}^n(u_1,\ldots,u_m)$ as follows. Let 
 $r^\ast_i=cu_i^2/(1-\exp(-cu_i))-cu_i^2/2-u_i$.
If $\sum_{i=1}^mu_i\ge 1-1/c$, then
 $\tilde{A}_{\delta,\epsilon}^n(u_1,\ldots,u_m)$ equals
 the intersection of $A_{\delta}^n(u_1,\ldots,u_m)$, the
 event that the numbers of the excess edges of the $m$ components are 
within the respective intervals
 $(n(r^\ast_i-\epsilon),n(r^\ast_i+\epsilon))$, 
 and the event
that any other connected component is of size less than  $n\epsilon$.
 If $\sum_{i=1}^mu_i<
 1-1/c$, then $\tilde{A}_{\delta,\epsilon}^n(u_1,\ldots,u_m)$
equals the intersection of $A_{\delta}^n(u_1,\ldots,u_m)$, the event that
there exists another connected component whose  size is in the interval
$(n(u^\ast-\epsilon),n(u^\ast+\epsilon))$, where 
$u^\ast/\br(1-\exp(-cu^\ast)\br)=1-\sum_{i=1}^mu_i$, the
 event that the numbers of the excess edges of these $m+1$ components are 
within the respective intervals
 $(n(r^\ast_i-\epsilon),n(r^\ast_i+\epsilon))$
for $i=1,2,\ldots,m$
and $(n(r^\ast-\epsilon),n(r^\ast+\epsilon))$, where 
 $r^\ast=c{u^\ast}^2/(1-\exp(-cu^\ast))-c{u^\ast}^2/2-u^\ast$,
 and  the event that any other connected component is of size less 
than  $n\epsilon$. 
\begin{corollary}
    \label{the:2}
Let $c_n\to c>0$ as $n\to\infty$. Then
\begin{align*}
        \lim_{\delta\to0}\limsup_{n\to\infty}\frac{1}{n}\log
      \mathbf{P}\bl(A_\delta^n(u_1,\ldots,u_m)\br)&=
\lim_{\delta\to0}\liminf_{n\to\infty}\frac{1}{n}\log
      \mathbf{P}\bl(A_\delta^n(u_1,\ldots,u_m)\br)\\
=\lim_{\substack{\delta\to0\\\epsilon\to0}}\limsup_{n\to\infty}\frac{1}{n}\log
      \mathbf{P}\bl(\tilde{A}_{\delta,\epsilon}^n(u_1,\ldots,u_m)\br)&= 
\lim_{\substack{\delta\to0\\\epsilon\to0}}\liminf_{n\to\infty}\frac{1}{n}\log
      \mathbf{P}\bl(\tilde{A}_{\delta,\epsilon}^n(u_1,\ldots,u_m)\br)\\
&=-\Bl(\sum_{i=1}^mK_c(u_i)+
L_c\bl(\sum_{i=1}^mu_i\br)\Br)
\end{align*}
and
\begin{align*}
      \lim_{\delta\to0}\liminf_{n\to\infty}
      \mathbf{P}\bl(\tilde{A}_{\delta,\epsilon}^n(u_1,\ldots,u_m)|
A_\delta^n(u_1,\ldots,u_m)\br)&= 1\,.
    \end{align*}
\end{corollary}
Let  $\beta^n$ denote the size of the largest connected component of 
$\mathcal{G}(n,c_n/n)$ and $\gamma^n$ denote the number of its excess edges. 
We  state results on the asymptotics of $(\alpha^n/n,\beta^n/n,\gamma^n/n)$.
\begin{corollary}
  \label{co:oc}Let $c_n\to c>0$ as $n\to\infty$. Then the following holds.
  \begin{enumerate}
  \item  The sequence $(\alpha^n/n,\beta^n/n,\gamma^n/n),\, n\in\N,$ 
obeys the LDP in
  $[0,1]^2\times\R_+$ for scale $n$ with action functional defined by
\begin{equation*}
    I_c^{\alpha,\beta,\gamma}(a,0,0)=L_c\bl((1-2a)^+\br)
+\frac{c}{2}\,(1-(1-2a)^+)^2%
\;\pi\Bl(\frac{2(1-a-(1-2a)^+)}{c(1-(1-2a)^+)^2}\Br),
\end{equation*}
\begin{multline*}
    I_c^{\alpha,\beta,\gamma}(a,u,r)=\sup_{\rho\in\R_+}\Bl( K_{\rho}(u)+
r\log\frac{\rho}{c}\Br)-K_c(u)\\+\inf_{\tau\in[(1-2a)\vee u,1-a]}\biggl(
  \Bigl\lfloor\frac{\tau}{u}\Bigr\rfloor K_c(u)
+K_c\Bl(\tau-u\Bigl\lfloor\frac{\tau}{u}\Bigr\rfloor\Br)
+L_c( \tau)+\frac{c}{2}(1-\tau)^2
\;\pi\Bl(\frac{2(1-a-\tau)}{%
 c(1-\tau)^2}\Br)\biggr)
\end{multline*}
if $u\in(0,1-a]$,
and $I_c^{\alpha,\beta,\gamma}(a,u,r)=\infty$ otherwise.
  \item  The sequence $(\beta^n/n,\gamma^n/n),\, n\in\N,$ 
obeys the LDP in
  $[0,1]\times\R_+$ for scale $n$ with action functional
  $I_c^{\beta,\gamma}$
defined by  $I_c^{\beta,\gamma}(0,0)=L_c((1-1/c)^+)$,
$I_c^{\beta,\gamma}(0,r)=\infty$ if $r>0$,
  \begin{multline*}
    I_c^{\beta,\gamma}(u,r)=\sup_{\rho\in\R_+}\Bl( K_{\rho}(u)+
r\log\frac{\rho}{c}\Br)+
\Bl(\Bl\lfloor \frac{1}{u}
\Bl(1-\frac{1}{c}\Br)\Br\rfloor-1\Br) K_c(u)+
K_c(\hat{u}\wedge u)\\+L_c\biggl(\Bl\lfloor
\frac{1}{u}\Bl(1-\frac{1}{c}\Br)\Br\rfloor u+\hat{u}\wedge u\biggr)
\end{multline*}
 if  $u\in(0,(1-1/c)^+)$, where $\hat{u}\in[0,1]$ satisfies the equality 
$\hat{u}/\bl(1-\exp(-c\hat{u})\br)=1-\lfloor\bl(1-1/c\br)/u\rfloor u$, and 
 by $I_c^{\beta,\gamma}(u,r)=\sup_{\rho\in\R_+}\bl( K_{\rho}(u)+
r\log(\rho/c)\br)+L_c(u)$ if $u\ge(1-1/c)^+$.
\item
 The sequence $\beta^n/n,\, n\in\N,$ obeys the LDP in $[0,1]$
for scale $n$  with action functional $I_c^\beta$ defined as follows:
 $I_c^\beta(0)=L_c\bl((1-1/c)^+\br)$,
  \begin{equation*}
    I_c^\beta(u)=
\Bl\lfloor \frac{1}{u}\Bl(1-\frac{1}{c}\Br)\Br\rfloor K_c(u)+
K_c(\hat{u}\wedge u)+L_c\biggl(\Bl\lfloor
\frac{1}{u}\Bl(1-\frac{1}{c}\Br)\Br\rfloor u+\hat{u}\wedge u\biggr)
\end{equation*}
 if  $u\in\bl(0,(1-1/c)^+\br)$,  and 
 $I_c^\beta(u)=K_c(u)+L_c(u)$ if $u\ge(1-1/c)^+$.
  \end{enumerate}
\end{corollary}
The next theorem considers normal  deviations of 
$(\alpha^n,\beta^n,\gamma^n)$.    We  recall that $\beta\in(0,1)$ is
 defined as the positive solution of  the equation
 $1-\beta=\exp(-\beta c)$ if $c>1$. For $c\le 1$ we define
  $\beta=0$. Let also $\alpha=1-\beta-c(1-\beta)^2/2$ and 
$\gamma=(c-1)\beta-c\beta^2/2$. 
\begin{theorem}
  \label{the:7}
Let $\sqrt{n}(c_n-c)\to\theta\in\R$ as $n\to\infty$, where $c>0$. Then the
following holds.
\begin{enumerate}
\item 
The sequence
$\sqrt{n}(\alpha^n/n-\alpha),\, n\in\N,$ 
converges in distribution in $\R$ as $n\to\infty$ to a 
 Gaussian random variable $\tilde{\alpha}$
with $\mathbf{E}\tilde{\alpha}=-\theta(1-\beta^2)/2$
and $\text{\bf{Var}}\,\tilde{\alpha}=\beta(1-\beta)+c(1-\beta)^2/2$.
\item
If, in addition, $c>1$, then the sequence
$(\sqrt{n}(\alpha^n/n-\alpha),\,\sqrt{n}(\beta^n/n-\beta),
\,\sqrt{n}(\gamma^n/n-\gamma)),\, n\in\N,$ 
converges in distribution in $\R^3$ as $n\to\infty$ to a 
 Gaussian random variable $(\tilde{\alpha},\,
 \tilde{\beta},\,\tilde{\gamma})$ with
 $\mathbf{E}\tilde{\beta}=\theta\beta(1-\beta)/\bl(1-c(1-\beta)\br)$,
$\mathbf{E}\tilde{\gamma}=\theta\beta^2/2$,
$\text{\bf{Var}}\,\tilde{\beta}=\beta(1-\beta)/\bl(1-c(1-\beta)\br)^2$, 
$\text{\bf{Var}}\,\tilde{\gamma}=\beta(1-\beta)+c\beta(3\beta/2-1)$,
 $\text{\bf{Cov}}\,(\tilde{\alpha},\tilde{\beta})=-
\beta(1-\beta)/\bl(1-c(1-\beta)\br)$,
$\text{\bf{Cov}}\,(\tilde{\alpha},\tilde{\gamma})=-\beta(1-\beta)(c-1) $, and 
 $\text{\bf{Cov}}\,(\tilde{\beta},\tilde{\gamma})=\beta(1-\beta)(c-1)/%
\bl(1-c(1-\beta)\br)$.
\end{enumerate}
\end{theorem}
We now state a moderate deviation asymptotics result for
$(\alpha^n,\beta^n,\gamma^n)$. We assume as given a  real-valued sequence
 $b_n,\,n\in\mathbb{N},$ such that $b_n\to\infty$ and
 $b_n/\sqrt{n}\to0$ as $n\to\infty$. Let $y^T$ denote the transpose of
 $y\in\R^3$. 
 \begin{theorem}
   \label{the:6}Let $(\sqrt{n}/b_n)(c_n- c)\to\hat{\theta}\in\R$ 
as $n\to\infty$,
   where $c>0$. Then the following holds.
   \begin{enumerate}
   \item 
The sequence $(\sqrt{n}/b_n)(\alpha^n/n-\alpha),\, n\in\N,$
obeys the LDP in $\R$ for scale $b_n^2$ with action functional
$(x-\mu_\alpha)^2/(2\sigma_{\alpha}^2),\,x\in\R$, where
$\mu_\alpha=-\hat{\theta}(1-\beta^2)/2$ and $\sigma_{\alpha}^2=
\beta(1-\beta)+c(1-\beta)^2/2$.
\item
If, in addition,
 $c>1$, then the sequence $\bl((\sqrt{n}/b_n)(\alpha^n/n-\alpha),\,
(\sqrt{n}/b_n)(\beta^n/n-\beta),(\sqrt{n}/b_n)(\gamma^n/n-\gamma)\br),\, 
n\in\N,$
obeys the LDP in $\R^3$ for scale $b_n^2$ with action functional
$(y-\mu)^T\Sigma^{-1}(y-\mu)/2,\,y\in\R^3$, where
$\mu=(\mu_\alpha,\mu_\beta,\mu_\gamma)^T$ and
$\Sigma=\biggl(
\begin{smallmatrix}
  \sigma_{\alpha}^2&\sigma_{\alpha\beta}&\sigma_{\alpha\gamma}\\
\sigma_{\alpha\beta}&\sigma_{\beta}^2&\sigma_{\beta\gamma}\\
\sigma_{\alpha\gamma}&\sigma_{\beta\gamma}&\sigma_{\gamma}^2
\end{smallmatrix}\biggr)$ are given by 
$\mu_\beta=\hat{\theta}\beta(1-\beta)/\bl(1-c(1-\beta)\br)$,
$\mu_\gamma=\hat{\theta}\beta^2/2$,
$\sigma_{\beta}^2=\beta(1-\beta)/\bl(1-c(1-\beta)\br)^2$, 
$\sigma_\gamma^2=\beta(1-\beta)+c\beta(3\beta/2-1)$,
 $\sigma_{\alpha\beta}=-\beta(1-\beta)/\bl(1-c(1-\beta)\br)$,
$\sigma_{\alpha\gamma}=-\beta(1-\beta)(c-1) $, and 
 $\sigma_{\beta\gamma}=\beta(1-\beta)(c-1)/%
\bl(1-c(1-\beta)\br)$.
   \end{enumerate}
 \end{theorem}
The list of results is concluded with the critical-graph case. 
We recall that 
excursions of a non-negative function $\bx=(\bx_t,\,t\in\R_+)$ are
defined as intervals $[s_i,t_i]$, where $s_i<t_i$,
 such that $\bx_{s_i}=\bx_{t_i}=0$ and
$\bx_p>0$ for $p\in(s,t)$, $t_i-s_i$ is called the excursion's length;
 continuous functions have at most countably many excursions.
Let, given $\tilde{\theta}\in\R$,
 process $\tilde{X}=(\tilde{X}_t,\,t\in\R_+)$ 
be defined as the Skorohod reflection of the process
$(W_t+\tilde{\theta} t-t^2/2,\,t\in\R_+)$,
where $W=(W_t,\,t\in\R_+)$ is a  Wiener process. 
In the next theorem,
$\tilde{U}=
(\tilde{U}_1,\tilde{U}_2,\ldots)$ is the sequence of the excursion
lengths of  $\tilde{X}$ arranged in descending order and
 $\tilde{R}=(\tilde{R}_1,\tilde{R}_2,\ldots)$ is  the sequence of the
 increments of the process 
$\bl(N_{\int_0^t\tilde{X}_s\,ds},\,t\in\R_+\br)$
over these excursions,
 where $(N_t,\,t\in\R_+)$ is a  Poisson process (independent of $W$).
Let $\breve{\mathbb{S}}$ denote the subspace of $\R_+^\N$ of
non-increasing  sequences
$\bu=(u_1,u_2,\ldots)$ equipped with induced topology.
The sequence $b_n$  is defined as in Theorem~\ref{the:6}.
\begin{theorem}
  \label{the:crit}
  \begin{enumerate}
  \item 
Let $n^{1/3}(c_n-1)\to\tilde{\theta}\in\R$ as $n\to\infty$. 
Then the sequences $\tilde{U}^n=
(U^n_1/n^{2/3},U^n_2/n^{2/3},\ldots)$ and
$\tilde{R}^n=(R^n_1/n^{2/3},R^n_2/n^{2/3},\ldots)$ jointly converge in
distribution in
$\breve{\mathbb{S}}\times\R_+^\N$
 to the respective sequences 
$\tilde{U}=(\tilde{U}_1,\tilde{U}_2,\ldots)$ and
$\tilde{R}=(\tilde{R}_1,\tilde{R}_2,\ldots)$. If, moreover,
$\sqrt{n}(c_n-1)\to\theta\in\R$ as $n\to\infty$ (so
$\tilde{\theta}=0)$, then the 
$(\tilde{U}^n,\tilde{R}^n)$
 are asymptotically independent
of $\sqrt{n}(\alpha^n/n-1/2)$ so that
$\bl(\sqrt{n}(\alpha^n/n-1/2),\tilde{U}^n,\tilde{R}^n\br)$
  jointly converge in distribution in 
$\R\times\breve{\mathbb{S}}\times\R_+^\N$ to 
$(\tilde{\alpha},\tilde{U},\tilde{R})$, where $(\tilde{U},\tilde{R})$
correspond to $\tilde{\theta}=0$,
$\tilde{\alpha}$ is independent of
$(\tilde{U},\tilde{R})$ and is Gaussian 
with $\mathbf{E}\tilde{\alpha}=-\theta/2$
and $\text{\bf{Var}}\,\tilde{\alpha}=1/2$.
  \item 
Let  $(n^{1/3}/b_n^{2/3})(c_n-1)\to\breve{\theta}\in\R $ as $n\to\infty$. 
Then the sequences $\breve{U}^n=
(U^n_1/(nb_n)^{2/3},U^n_2/(nb_n)^{2/3},\ldots)$ and
$\breve{R}^n=(R^n_1/(nb_n)^{2/3},R^n_2/(nb_n)^{2/3},\ldots)$ jointly
obey the LDP
 in $\breve{\mathbb{S}}\times\R_+^\N$ for scale $b_n^2$ 
with action functional
\begin{equation*}
  \breve{I}^{U,R}_{\breve{\theta}}(\bu,\bfr)=
-\frac{1}{24}\sum_{i=1}^\infty u_i^3
+\frac{1}{6}\,
\Bl(\sum_{i=1}^\infty u_i-\breve{\theta}\Br)^3\vee0
+\frac{\breve{\theta}^3}{6}+
\frac{1}{24}\sum_{i=1}^\infty
u_i^3\,\pi\Bl(\frac{12\,r_i}{u_i^3}\Br)
\end{equation*}
if $\sum_{i=1}^\infty u_i<\infty$ and $r_i=0$ when $u_i=0$, and 
$  \breve{I}^{U,R}_{\breve{\theta}}(\bu,\bfr)=\infty$ otherwise, where
$\bu=(u_1,u_2,\ldots)$ and $\bfr=(r_1,r_2,\ldots)$. 
If, moreover,
$(\sqrt{n}/b_n)(c_n-1)\to\hat{\theta}$ as $n\to\infty$ (so
$\breve{\theta}=0$),  then  the 
$\bl((\sqrt{n}/b_n)(\alpha^n/n-1/2),\breve{U}^n,\breve{R}^n\br)$ obey the
LDP in $\R\times\breve{\mathbb{S}}\times\R_+^\N$ with action functional
\begin{equation*}
  \breve{I}^{\alpha,U,R}_{\theta}(a,\bu,\bfr)=
\Bl(a+\frac{\hat{\theta}}{2}\Br)^2+  \breve{I}^{U,R}_{0}(\bu,\bfr).
\end{equation*}
  \end{enumerate}
\end{theorem}
\begin{corollary}
  \label{co:crit}
Let  $(n^{1/3}/b_n^{2/3})(c_n-1)\to\breve{\theta}\in\R $ as
$n\to\infty$. Then the following holds.
\begin{enumerate}
\item 
The sequence $\breve{U}^n,\,n\in\N,$ obeys
the LDP in $\breve{\mathbb{S}}$ 
for scale $b_n^2$ with action functional
\begin{equation*}
  \breve{I}^U_{\breve{\theta}}(\bu)=-\frac{1}{24}\sum_{i=1}^\infty u_i^3
+\frac{1}{6}\,
\Bl(\sum_{i=1}^\infty u_i-\breve{\theta}\Br)^3\vee0
+\frac{\breve{\theta}^3}{6}
\end{equation*}
if $\sum_{i=1}^\infty u_i<\infty$ and
$\breve{I}^U_{\breve{\theta}}(\bu)=\infty$ otherwise.
\item
The sequence $\beta^n/(nb_n)^{2/3},n\in\N,$ obeys the LDP in $\R_+$ for scale
$b_n^2$ with action functional $\breve{I}^\beta_{\breve{\theta}}$ given by
$\breve{I}^\beta_{\breve{\theta}}(0)=\breve{\theta}^3\vee0/6$,
\begin{equation*}
  \breve{I}^\beta_{\breve{\theta}}(u)=
-\Bigl\lfloor\frac{\breve{\theta}}{u}\Bigr\rfloor
  \frac{u^3}{24}-\frac{1}{24}\biggl(\Bl(2\Bl(\breve{\theta}-
\Bigl\lfloor\frac{\breve{\theta}}{u}\Bigr\rfloor u\Bl)\Br)\wedge u\biggr)^3
  +\frac{1}{6}\,\Bl(\Bigl\lfloor\frac{\breve{\theta}}{u}\Bigr\rfloor u
+\Bl(2\Bl(\breve{\theta}-
\Bigl\lfloor\frac{\breve{\theta}}{u}\Bigr\rfloor u\Bl)\Br)\wedge u
-\breve{\theta}\Br)^3+\frac{\breve{\theta}^3}{6}
\end{equation*}
if $u\in(0,\breve{\theta}^+)$
and 
\begin{equation*}
  \breve{I}^\beta_{\breve{\theta}}(u)=
-\frac{u^3}{24}+\frac{(u-\breve{\theta})^3}{6}+
\frac{\breve{\theta}^3}{6}
\end{equation*}
if $u\ge\breve{\theta}^+$.
\end{enumerate}
\end{corollary}
We now comment on the results and relate them to earlier ones. 
Equation \eqref{eq:11} in a slightly different form appears in 
Barraez,  Boucheron, and Fernandez de la Vega \citeyear{BarBouFer00}.
Considering the sequences $\overline{U}^n$ and 
$\overline{R}^n$ has been prompted  by the form of 
the results of Aldous \citeyear{Ald97}.
Corollary~\ref{the:2} implies, in particular, that provided there
 exist $m$ components  asymptotically $nu_i$ in size, where
 $\sum_{i=1}^mu_i<1-1/c$, then with probability close to 1 there exists another
 giant component. This can be explained by noting that the number of
 vertices outside of the $m$ components is asymptotically equal to 
$n(1-\sum_{i=1}^mu_i)$, so the  ``effective'' expected  degree of an 
outside node is
$c(1-\sum_{i=1}^mu_i)>1$, which means there is enough potential for
 another giant component.

 Part 3 of Corollary~\ref{co:oc}  is due to
 O'Connell \citeyear{Oco98}, who provides an alternative, elegant form of
 the action functional for $c>1$ and $u>0$:
$       I_c^\beta(u)=
   kK_c(u)+L_c(ku)$ for $ u\in[x_{k},x_{k-1}],\,k\in\N$,
where   $x_0=1$ and the
$x_k,\,k\in\N,$  are the solutions of the equations
$    x_k/\bl(1-\exp(-cx_k)\br)=1-kx_k$.
(Note that the expression for the action functional
 in Theorem~3.1 of  O'Connell \citeyear{Oco98}
  has a misprint.)
O'Connell \citeyear{Oco98} also noted that the action functional
 $I_c^\beta$ is not convex.
The advantage of the form of $I_c^\beta$  used in
 Corollary~\ref{co:oc} is that it is
 suggestive of  the structure of
the most probable configuration with the
  largest component  asymptotically  $nu$ in size: 
if   $u\ge 1-1/c$, then the component of  size
 $nu$ is the only giant component,
     while if $u< 1-1/c$, then there are $\lfloor (1-1/c)/u\rfloor$
  components, whose sizes are  asymptotically
 $nu$, and one component asymptotically  $n(\hat{u}\wedge u)$ in size.
(A similar remark has been made by O'Connell \citeyear{Oco98}.) 
This conjecture  is confirmed by the
 proof of  Corollary~\ref{co:oc}. In
    addition, the number of components in an optimal configuration is
    asymptotically equal to $n\bl(1-u-c(1-u)^2/2\br)$ if $u\ge1-1/c$ and
$n\bl(1-\hat{\tau}-c(1-\hat{\tau})^2/2\br)$  if $u<1-1/c$, where
$\hat{\tau}=\lfloor (1-1/c)/u\rfloor u+\hat{u}\wedge u$.

\begin{figure}
\centering
\scalebox{0.7}{\includegraphics{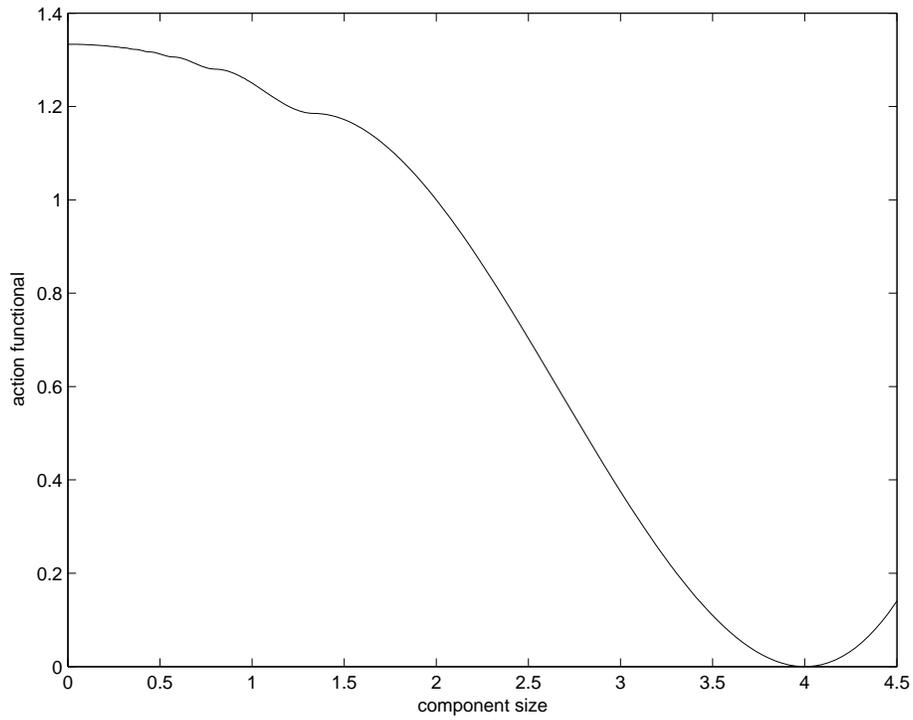}}
\caption{ Moderate deviations of
 the size of the largest component of the critical
  graph ($\breve{\theta}=2$)}\label{fi:3}
\end{figure}

\begin{figure}
\centering
\scalebox{0.7}{\includegraphics{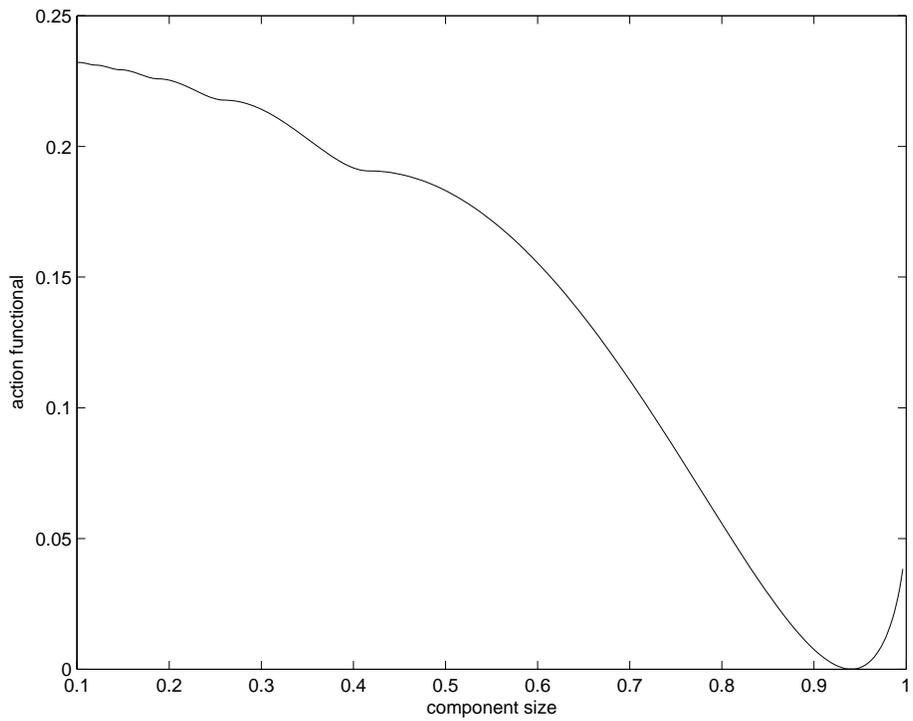}}
\caption{ Large deviations of the size of the largest component of
  $\mathcal{G}(n,3/n)$
\text{(O'Connell (1998))}}\label{fi:4}
\end{figure}

Corollary~\ref{co:crit} leads to similar conclusions.
The action functional $\breve{I}^\beta_{\breve{\theta}}(u)$ 
 can be written for $\breve{\theta}>0$ and 
$u\in(0,2\breve{\theta}]$ as 
$\breve{I}^\beta_{\breve{\theta}}(u)=-ku^3/24+(ku-\breve{\theta})^3/6+
\breve{\theta}^3/6$, where $k\in\N$ is such that 
$u\in[\breve{\theta}/(k+1/2),\breve{\theta}/(k-1/2)]$. It is not
convex for $\breve{\theta}>0$ either. 
Fig.~\ref{fi:3}  shows the action functional for $\breve{\theta}=2$.
(Note that the form of the curve is the
same for all $\breve{\theta}>0$ since 
$\breve{I}^\beta_{x\breve{\theta}}(xu)=x^3\breve{I}^\beta_{\breve{\theta}}(u)$
for $x>0$.)
Interestingly,  the graph of $\breve{I}^\beta_{\breve{\theta}}$
is reminiscent of the one of $I^\beta_{c}$  given
in O'Connell's paper \citeyear{Oco98}, 
which we reproduce in Fig.~\ref{fi:4} for 
comparison's sake.
For $\breve{\theta}>0$, 
the most probable configuration with the largest component
asymptotically  $(nb_n)^{2/3}u$ in   size
consists of only one such component if 
$u\ge\breve{\theta}$ and 
has $\lfloor\breve{\theta}/u\rfloor$ components asymptotically 
$(nb_n)^{2/3}u$ in  size  along with one component asymptotically 
$(nb_n)^{2/3}\bl((2(\breve{\theta}-
\lfloor\breve{\theta}/u\rfloor u))\wedge u\br)$ in size if $u< \breve{\theta}$.
Since the action functional $\breve{I}^\beta_{\breve{\theta}}(u)$ 
equals zero at the only
point $u=2\breve{\theta}^+$,  the $\beta^n/(nb_n)^{2/3}$ converge
 in probability to $2\breve{\theta}^+$ as $n\to\infty$, which is
 consistent with the asymptotics $\beta/(c-1)\to2$ as $c\downarrow1$.
There is also an analogue  for the critical
graph of Corollary~\ref{the:2} on the 
most probable ``conditional'' configurations. In particular,
given there exists a component asymptotically  $(nb_n)^{2/3}u$ in size,
 with probability tending to 1 it has asymptotically 
$(nb_n)^{2/3}u^3/12$ excess edges, and if $u<\breve{\theta}$ also with
probability tending to $1$ there exists another component
asymptotically  $(nb_n)^{2/3}2(\breve{\theta}-u)$ in size.

The first assertion of part 1 of Theorem~\ref{the:crit} is due to
Aldous \citeyear{Ald97}, who establishes the convergence of the sizes of the
connected components   for the stronger 
$\ell_2$ topology. 
 Our proof  uses similar ideas.
Part 2 of Theorem~\ref{the:crit} can also be expressed as a statement
on a certain type of convergence of excursions.
Let idempotent
process $\breve{X}=(\breve{X}_t,\,t\in\R_+)$ 
be defined as the reflection of the idempotent process
$(\breve{W}_t+\breve{\theta} t-t^2/2,\,t\in\R_+)$, 
where $\breve{W}=(\breve{W}_t,\,t\in\R_+)$ is an
idempotent Wiener process, 
and let  $(\breve{N}_t,\,t\in\R_+)$ be an
 idempotent Poisson process independent of $\breve{W}$
(for the notions of idempotent probability the reader is referred
either to  Puhalskii \citeyear{Puh01} or  the appendix). Let 
$(\breve{U}_1,\breve{U}_2,\ldots)$ be the sequence of the excursion
lengths of $\breve{X}$ arranged in descending order and
$(\breve{R}_1,\breve{R}_2,\ldots)$ be the sequence of 
the increments of 
$\bl(\breve{N}_{\int_0^t\breve{X}_p\,dp},\,t\in\R_+\br)$ over these excursions.
Then the sequences $(U^n_1/(nb_n)^{2/3},U^n_2/(nb_n)^{2/3},\ldots)$ and
$(R^n_1/(nb_n)^{2/3},R^n_2/(nb_n)^{2/3},\ldots)$ jointly
large deviation  converge in
distribution in $\breve{\mathbb{S}}\times\R_+^\N$ at rate $b_n^2$ 
to the respective sequences 
$(\breve{U}_1,\breve{U}_2,\ldots)$ and
$(\breve{R}_1,\breve{R}_2,\ldots)$ as $n\to\infty$. (The definition of
large deviation convergence is recalled in Section~\ref{sec:techn-prel}.)
Thus, the actual assertion combines  statements on
large deviation convergence and on the idempotent distribution of the limit.
The LDP for $(\overline{U}^n,\overline{R}^n)$ of Corollary~\ref{the:4}
admits a similar reformulation.

Part 1 of  Theorem~\ref{the:7} for the case 
where $c_n=c$ and accordingly $\theta=0$
  is due to
Stepanov \citeyear{Ste70a,Ste72}. 
 Part 2 of  Theorem~\ref{the:7} complements the results of
Stepanov \citeyear{Ste70a,Ste72} 
(see also Pittel \citeyear{Pit90}, Barraez,  Boucheron, and Fernandez de
  la Vega  \citeyear{BarBouFer00})  by allowing for
  $\theta\not=0$, 
incorporating $\gamma^n$ and indicating
the  covariance of $\tilde{\alpha}$ and $\tilde{\beta}$.  
As to be expected, the latter two random variables are negatively correlated.
Parts 1 and 2 of    Theorem~\ref{the:6}, equivalently, state that the 
$(\sqrt{n}/b_n)(\alpha^n/n-\alpha)$ and 
$\br((\sqrt{n}/b_n)(\alpha^n/n-\alpha),\,
(\sqrt{n}/b_n)(\beta^n/n-\beta),(\sqrt{n}/b_n)(\gamma^n/n-\gamma)\br)$
 large deviation  converge
at rate $b_n^2$ as $n\to\infty$ to  Gaussian idempotent 
variables  with respective parameters 
$(\mu_\alpha,\sigma_\alpha^2)$ and
$(\mu,\Sigma)$.   
This formulation not only emphasises  analogy with
Theorem~\ref{the:7} but is instrumental in the proof below.

\begin{figure}
\centering
\scalebox{0.69}{\includegraphics{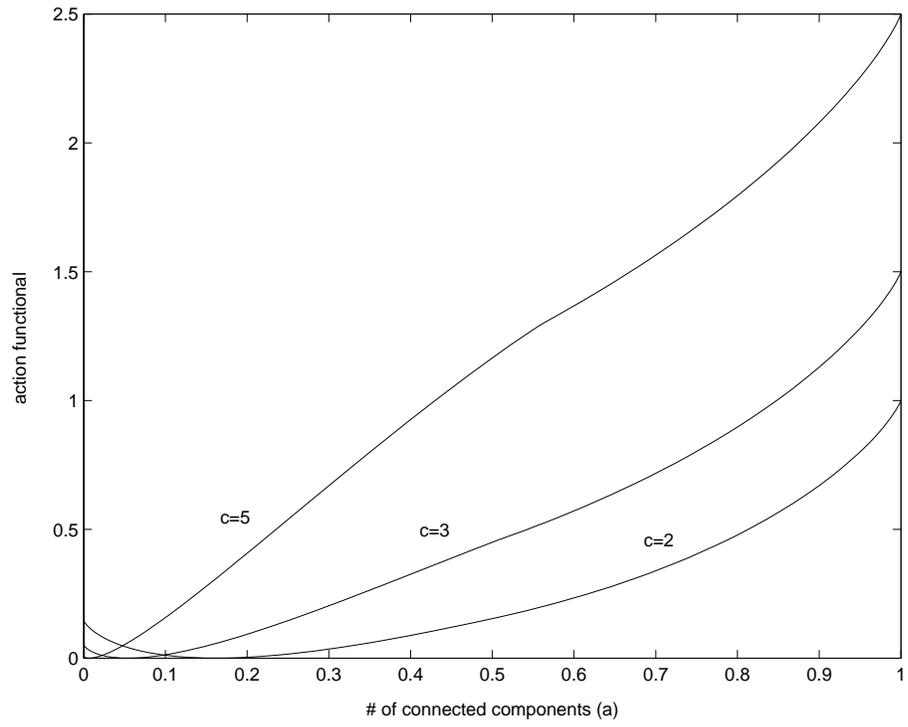}}
\caption{ LDP for the number of connected components}\label{fi:1}
\end{figure}

\begin{figure}
\centering\scalebox{0.69}{\includegraphics{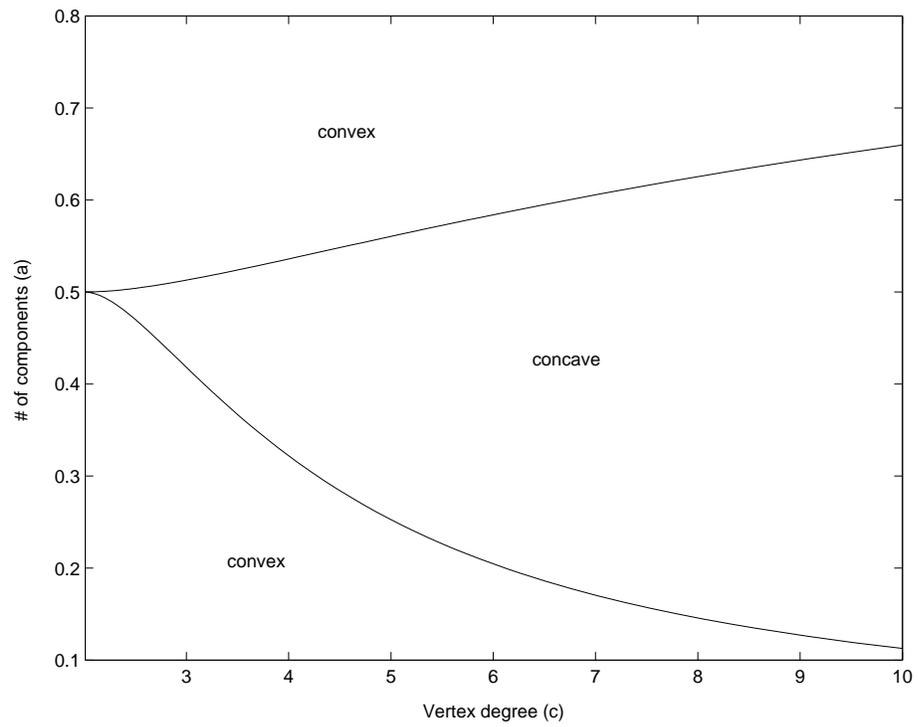}}
\caption{ Convexity-concavity regions for $I^\alpha_c$}\label{fi:2}
\end{figure}

We now consider implications of the LDP for $\alpha^n/n$ 
of Corollary~\ref{the:5}, which provide some revealing insights.
The derivative with respect to $\tau$ of the function in the infimum on
the right of the expression for $I_c^\alpha$  equals 
\begin{equation*}
2\Bl(1-\frac{a}{1-\tau}\Br)-\frac{c\tau}{e^{c\tau}-1}
-\log\Bl(2\Bl(1-\frac{a}{1-\tau}\Br)\Br)
+\log\frac{c\tau}{e^{c\tau}-1}.
\end{equation*}
Since $\tau\ge 1-2a$, the derivative is non-negative if and only if 
$a\ge (1-\tau)\bl(1-c\tau/(2(e^{c\tau}-1)\br)\br)$.
The  function on the right of the latter inequality, 
as a function of $\tau\in[0,1]$,
is concave, is decreasing if $c\le 2$, 
and is first increasing and then decreasing
if $c>2$. Let  $a^\ast\in[1/2,1]$ denote
  the maximum of this function  
 on $[0,1]$. For $a\in[0,a^\ast]$  the equation
  \begin{equation}
    \label{eq:63}
a=(1-\tau)\Bl(1-\frac{1}{2}\frac{c\tau}{e^{c\tau}-1}\Br)
\end{equation}
has one root if either $a<1/2$ or $a=a^\ast$
  and has two roots otherwise. 
Let $\tau^\ast(a)\in[0,1]$, where $a\in[0,a^\ast]$, denote 
the greatest  root of \eqref{eq:63}.
 Then the infimum on
the right of the expression for $I_c^\alpha$ 
 is attained  at $\tau=\tau^\ast(a)$
  if $a\in[0,1/2]$, at $\tau=0$ if $a\in[a^\ast,1]$, and
either at $\tau=\tau^\ast(a)$ or $\tau=0$ if $a\in(1/2,a^\ast)$.
Accordingly,  the optimal configuration  has either
one giant component asymptotically  $n\tau^\ast(a)$ in size or no giant
  components.
 We can therefore write
\begin{equation}
  \label{eq:65}
    I_c^\alpha(a)=K_c(\tau^\ast(a))+L_c(\tau^\ast(a))+
\frac{c(1-\tau^\ast(a))^2}{2}
\;\pi\Bl(\frac{2(1-a-\tau^\ast(a))}{c(1-\tau^\ast(a))^2}\Br)
\end{equation}
if $a\in[0,1/2]$,
\begin{equation}
  \label{eq:67}
    I_c^\alpha(a)=\biggl(\frac{c}{2}\;\pi\Bl(\frac{2(1-a)}{c}\Br)\biggr)
\wedge 
\biggl( K_c(\tau^\ast(a))+L_c(\tau^\ast(a))+
\frac{c(1-\tau^\ast(a))^2}{2}
\;\pi\Bl(\frac{2(1-a-\tau^\ast(a))}{c(1-\tau^\ast(a))^2}\Br)\biggr)
\end{equation}
if $a\in(1/2,a^\ast)$, and 
\begin{equation}
  \label{eq:68}
    I_c^\alpha(a)= \frac{c}{2}
\;\pi\Bl(\frac{2(1-a)}{c}\Br)
\end{equation}
if $a\in[a^\ast,1]$. 
If $c\le 2$,  the action functional is in fact given by \eqref{eq:65} and
\eqref{eq:68} since $a^\ast=1/2$. 
It is seen to be convex and differentiable in $a$. 
If $c>2$, then
$a^\ast>1/2$,  the difference between the  first and the second 
functions in the minimum on the right of \eqref{eq:67} is positive for
$a=1/2$, is decreasing
in $a$ for $a>1/2$, 
 and there exists a unique $\hat{a}\in(1/2,a^\ast)$ where these two functions
are equal. Thus,  for $c>2$ 
\begin{equation*}
  I_c^\alpha(a)=
  \begin{cases}\displaystyle
    K_c(\tau^\ast(a))+L_c(\tau^\ast(a))+
\frac{c(1-\tau^\ast(a))^2}{2}
\;\pi\Bl(\frac{2(1-a-\tau^\ast(a))}{c(1-\tau^\ast(a))^2}\Br)&\text{ if }
a\in[0,\hat{a}],\\\displaystyle
\frac{c}{2}\;\pi\Bl(\frac{2(1-a)}{c}\Br)&\text{ if }a\in[\hat{a},1].
  \end{cases}
\end{equation*}
For $c>2$ the function $I_c^\alpha(a)$ is strictly convex to
 the right of $\hat{a}$ and is strictly  concave in a neighbourhood 
to its left.
As a matter of fact, there exists $\tilde{a}\in(0,1/2)$
 such that  $I_c^\alpha(a)$ is strictly convex for $a<\tilde{a}$ (and
 $a>\hat{a}$), and is strictly concave for $\tilde{a}<a<\hat{a}$. 
The value of $\tilde{a}$ is
 given by \eqref{eq:63} for $\tau=\tilde{\tau}$, where
 $\tilde{\tau}$  solves the equation $\exp(-c\tau)-1+c\tau=c\tau^2$. In
addition, $\tilde{a}\downarrow 0$ and $\hat{a}\uparrow 1$ as $c\to\infty$ 
(in fact, $\tilde{a}< 2/c$ for $c>2$), so the concavity region grows as
 $c$ does. Fig.~\ref{fi:1}
 shows the action functionals for various values of $c$
 and Fig.~\ref{fi:2} shows the regions of
 convexity and  concavity for $I^\alpha_c$.

Another distinguishing feature of
point $\hat{a}$ is that at it the left derivative of
$I_c^\alpha(a)$ is greater than the right one,
$I^\alpha_c(a)$ being differentiable in $a$ elsewhere. It is, moreover, a
 point of phase transition:  for
 $a<\hat{a}$ the most probable configuration has one giant component
 asymptotically  $n\tau^\ast(a)$ in  size 
while for $a>\hat{a}$ it is optimal to
 have no giant components. Hence, for $c>2$ 
we have the following structure of the random graph with a given number  of
 components of order $na$: for small values of $a$, with probability close to 1
 there is one giant component asymptotically   $n\tau^\ast(a)$ in size
 and many (actually asymptotically $na$) 
small components of sizes not greater than of order $o(n)$ (it can be
conjectured  their sizes are of order $\log n$ or less);  
as $a$ increases, more small components split 
off from the giant
 component and the size of the giant component decreases gradually;
 however, at $a=\hat{a}$ the giant component breaks up in that its
 size drastically reduces from being of order 
$n\tau^\ast(\hat{a})$ to being of order $o(n)$, and for
 $a>\hat{a}$ only small components remain which disintegrate further
 as $a$ increases. If $c\le 2$, then as $a$
 increases the giant component, which is
 asymptotically     $n\tau^\ast(a)$ in size, gradually
 decreases in size and disappears at $a=1/2$, so no drastic changes
 occur. We thus shed new light on
 the observation by Stepanov \citeyear{Ste70} (see
 also Bollob{\'a}s, Grimmett, and Janson 
\citeyear{BolGriJan96}) of
 $c=2$ being a critical point.


There is another connection between our results and those of Stepanov
\citeyear{Ste70} as well as of Bollob{\'a}s, Grimmett, and Janson 
\citeyear{BolGriJan96}, to which we alluded in the introduction. 
The above observation has been made by Stepanov
\citeyear{Ste70} on the basis of the asymptotics
 \begin{equation}
   \label{eq:64}
    \lim_{n\to\infty}
\frac{1}{n}\log \mathbf{E} e^{\lambda \alpha^n}=S_c(\lambda),\,\lambda\in\R,
    \end{equation}
where
    \begin{multline}
              \label{eq:6}
            S_c(\lambda)=
\sup_{\tau\in[(1-e^\lambda/c)^+,1]}\Bl(\lambda(1-\tau)+
\frac{c}{2}\,(1-\tau)^2e^{-\lambda}
-(1-\tau)\log(1-\tau)\\-\frac{c}{2}\,(1-\tau^2)
-\tau\log\tau+\tau\log(1-e^{-c\tau})\Br),
    \end{multline}
and a subsequent analysis of the function $S_c(\lambda)$.
We are able to reproduce \eqref{eq:64}
by using the LDP for $\alpha^n/n$ 
   and Varadhan's lemma, see, e.g., Dembo and Zeitouni \citeyear{DemZei}.
Moreover, 
since $I^\alpha_c$ is strictly convex for $c\le 2$,
it is possible to derive the LDP for $\alpha^n/n$ of
Corollary~\ref{the:5} from limit \eqref{eq:64} via 
G{\"a}rtner's theorem, see G{\"a}rtner \citeyear{gar} or Freidlin and Wentzell
\citeyear{wf2}, so that $I_c(\alpha)$ is given by the Legendre-Fenchel
transform of $S_c(\lambda)$.
 This has been  done actually by Bollob{\'a}s, Grimmett, and Janson 
\citeyear{BolGriJan96}, who  obtain asymptotics \eqref{eq:64} independently of 
Stepanov \citeyear{Ste70} and, in effect, provide a
 solution to the optimisation
problem \eqref{eq:6}, though they do not find the form of
$I^\alpha_c$ in Corollary~\ref{the:5}.
However, for $c>2$ G{\"a}rtner's  theorem is not applicable because of
``the onset of concavity'' described above. The
 Legendre-Fenchel transform of $S_c(\lambda)$, being the convex hull of 
$I_c(\alpha)$, no longer coincides with $I_c(\alpha)$, which
 causes Bollob{\'a}s, Grimmett, and Janson's
\citeyear{BolGriJan96}  stopping short of obtaining the above LDP.

\section{Technical preliminaries}
\label{sec:techn-prel}
In this section we collect pieces of terminology and notation used
throughout the paper, recall some results on weak convergence and
large deviation asymptotics pertinent to the developments below, and 
provide a number of auxiliary lemmas.

We denote by  $\mathbb{D}_C([a,b],\R^d)$, where  $d\in\N$,   the
space of right-continuous with left-hand limits
 $\R^d$-valued functions on an interval $[a,b]$
equipped with uniform metric and 
Borel $\sigma$-algebra. 
Space $\D(\R_+,\R^d)$  is defined as
the  space of $\R^d$-valued right-continuous with
  left-hand limits functions on $\R_+$ equipped with the
  Skorohod topology and Borel $\sigma$-algebra.  
Spaces $\C([0,1],\R^d)$ and $\C(\R_+,\R^d)$ are the
subspaces of the respective spaces $\D_C([0,1],\R^d)$ and $\D(\R_+,\R^d)$
consisting of continuous    functions with induced topologies.
Elements of these spaces are mostly denoted by boldface lower-case
Roman letters, e.g., $\bx=(\bx_t,\,t\in[a,b])$; $\bx_{t-}$ denotes the
left-hand limit at $t$; $\dot{\bx}_t$ 
denotes the  Radon-Nykodim derivative at $t$ 
  with respect to Lebesgue measure of an absolutely continuous
$\bx$. We denote by $p_1$ the
projection $(\bx_t,\,t\in\R_+)\to(\bx_t,\,t\in[0,1])$ from $\D(\R_+,\R^d)$
to $\D_C([0,1],\R^d)$ and note that it is continuous at 
$\bx\in\C(\R_+,\R^d)$. 
Maps $\mathcal{R}$ and $\mathcal{T}$ from $\D(\R_+,\R)$ to
$\D(\R_+,\R)$ are defined by
$\mathcal{R}(\bx)_t=\bx_t-\inf_{s\in[0,t]}\bx_s\wedge0$ and 
$\mathcal{T}(\bx)_t=-\inf_{s\in[0,t]}\bx_s\wedge0$.  If $\bx_0\ge0$,
then the functions $\by=\mathcal{R}(\bx)$ and 
$\phi =\mathcal{T}(\bx)$ can be equivalently defined as a solution to
a Skorohod problem in that $\by=\bx+\phi$, $\by_t\ge0$, $\phi$ is
 non-decreasing with $\phi_0=0$ and
 $\phi_t=\int_0^t\ind(\by_s=0)\,d\phi_s,\,t\in\R_+$.
Unless specified otherwise,
 ``almost everywhere (a.e.)''  refers to Lebesgue measure and
product topological spaces are  
equipped with product topologies; besides,
$\inf_\emptyset$ is understood as $\infty$ and
 $\mathcal{B}(\R)$ denotes the
 Borel $\sigma$-algebra on $\R$.

We assume that all the random objects we consider
are defined on a complete probability space
$(\Omega,\mathcal{F},\mathbf{P})$, 
the expectation of a random variable $\xi$ is
   denoted as $\mathbf{E}\xi$. 
For a sequence of $\R^d$-valued random variables $\xi_n,\, n\in\N,$ and a
  sequence of real numbers $k_n\to\infty$
  we write $\xi_n\overset{\mathbf{P}^{1/k_n}}{\to}0$ and say that the $\xi_n$
  tend to zero super-exponentially in probability at rate $k_n$ 
if $\lim_{n\to\infty}
\mathbf{P}(\abs{\xi_n}>\epsilon)^{1/k_n}=0$ for arbitrary $\epsilon>0$.
We also let $\overset{\mathbf{P}}{\to}$ denote convergence in probability,
 $\overset{d}{\to}$  denote convergence in distribution in the
associated metric space, and $\xrightarrow[k_n]{ld}$  denote
large deviation (LD) 
convergence in distribution at rate $k_n$. To recall the definition
of the latter (see, e.g., Puhalskii \citeyear{Puh01}), we say that a 
$[0,1]$-valued function 
$\mathbf{\Pi}$ defined on the power set of a metric space $\Upsilon$ is a
deviability on $\Upsilon$ if
 $\mathbf{\Pi}(\Gamma)=\sup_{\upsilon\in\Gamma}
\exp(-I(\upsilon)),\,\Gamma\subset \Upsilon$, where 
$I$  is an action functional on $\Upsilon$, i.e.,  a
$[0,\infty]$-valued function on $\Upsilon$ such that the sets
$\{\upsilon\in\Upsilon:\, I(\upsilon)\le a\}$ are compact for $a\in\R_+$.
We say that a sequence $\mathbf{P}_n,\,n\in\N,$ of probability measures on
the Borel $\sigma$-algebra of $\Upsilon$ 
 LD converges at rate $k_n$ to a deviability $\mathbf{\Pi}$ on
$\Upsilon$ if 
$\lim_{n\to\infty}\bl(\int_\Upsilon f(\upsilon)^{k_n}\,d\mathbf{P}_n(\upsilon)
\br)^{1/k_n}=\sup_{\upsilon\in\Upsilon}f(\upsilon)\mathbf{\Pi}(\{\upsilon\})$
for every continuous bounded $\R_+$-valued function $f$ on
$\Upsilon$.
 Equivalently, the sequence $\mathbf{P}_n,\,n\in\N,$ 
LD converges at rate $k_n$ to $\mathbf{\Pi}$ if
it obeys the LDP with action functional $I$ for scale $k_n$.
We recall that if the sequence $\mathbf{P}_n$ is exponentially tight of
order $k_n$, i.e., for every $\epsilon>0$ there exists a compact
$K\subset\Upsilon$ such that
$\limsup_{n\to\infty}\mathbf{P}_n(\Upsilon\setminus K)^{1/k_n}<\epsilon$, then
it is LD relatively sequentially compact, i.e., there exists a
subsequence $\mathbf{P}_{n'}$ that LD converges at rate $k_{n'}$ to a
deviability $\mathbf{\Pi}'$; every such deviability is called an LD
accumulation point of the $\mathbf{P}_n$.
We also say that  a sequence of
random variables $ X_n,\, n\in\N,$ with values in $\Upsilon$ 
 LD converges in distribution at rate $k_n$ to a Luzin idempotent
variable $X$ with values in $\Upsilon$  
if the sequence of laws of the $X_n$ LD converges
at rate $k_n$ to the idempotent distribution of $X$.

Let $H^n=(H^n_t,\,t\in\R_+),\, n\in\N,$ 
be a sequence of $\R^d$-valued stochastic processes having
right-continuous with left-hand limits paths. 
The sequence $H^n$ 
is said to be $\C$-tight if the sequence of the distributions of the $H^n$ on
$\D(\R_+,\R^d)$  is tight for weak convergence of probability measures
on $\D(\R_+,\R^d)$ with its every  accumulation point being the law of
a continuous process. The following limits provide
 necessary and sufficient conditions for  $\C$-tightness:
\begin{align*}
  \lim_{B\to\infty}\limsup_{n\to\infty}
\mathbf{P}\bl(\abs{H^n_0}>B\br)=0,&&
  \lim_{\delta\to0}\limsup_{n\to\infty}
\mathbf{P}\bl(\sup_{\substack{s,t\in[0,T]:\\\abs{s-t}\le\delta}}
\abs{H^n_t-H^n_s}>\epsilon\br)=0,\;T\in\R_+,\,\epsilon>0.
\end{align*}
The sequence  $H^n$ is said to be $\C$-exponentially tight of
order $k_n$ if the sequence of the distributions of the $H^n$
  is exponentially tight of order $k_n$ 
as a sequence of probability measures on $\D(\R_+,\R^d)$
 and its every LD accumulation point  $\mathbf{\Pi}$ is such that $\mathbf{\Pi}(\bx)=0$
for every $\bx\in\D(\R_+,\R^d)\setminus\C(\R_+,\R^d)$.
The sequence of laws of the 
$H^n$ is $\C$-exponentially tight of order $k_n$ if and
only if
\begin{align*}
  \lim_{B\to\infty}\limsup_{n\to\infty}
\mathbf{P}\bl(\abs{H^n_0}>B\br)^{1/k_n}=0,&&
  \lim_{\delta\to0}\limsup_{n\to\infty}
\mathbf{P}\bl(\sup_{\substack{s,t\in[0,T]:\\\abs{s-t}\le\delta}}
\abs{H^n_t-H^n_s}>\epsilon\br)^{1/k_n}=0,\;T\in\R_+,\,\epsilon>0.
\end{align*}

We 
denote by $\xi^n_{ij}$ and $\zeta^n_{ij}$, where 
$i\in\N,\,j\in\N,\,n\in\N$, i.i.d. Bernoulli random
variables on $(\Omega,\mathcal{F},\mathbf{P})$ 
with $\mathbf{P}(\xi^n_{ij}=1)=c_n/n$
and 
define  $\mathcal{F}^n_t,\,t\in\R_+,$ as
the $\sigma$-algebras generated by the 
$\xi^n_{ij}$ and $\zeta^n_{ij}$ for $
i=1,2,\ldots, \lfloor n(t\wedge1)\rfloor,\,j\in\mathbb{N}$, 
completed with sets of $\mathbf{P}$-measure zero, and introduce  filtrations
$\F^n=(\mathcal{F}^n_t,\,t\in\R_+)$. 
\begin{lemma}
  \label{le:eps}
Let $c_n\to c>0$ as $n\to\infty$.
Let  $b_n\to\infty$ and 
$b_n/\sqrt{n}\to0$ as $n\to\infty$. The following convergences
 hold as $n\to\infty$:
\begin{equation*}
  \sup_{t\in[0,1]}\abs{\overline{\epsilon}^n_t}
\stackrel{\mathbf{P}^{1/n}}{\to}0,\qquad
  \sup_{t\in[0,1]} \frac{\sqrt{n}}{b_n}
\abs{\overline{\epsilon}^n_t}\stackrel{\mathbf{P}^{1/b_n^2}}{\to}0,\qquad
  \sup_{t\in[0,1]}\sqrt{n}
\abs{\overline{\epsilon}^n_t}\stackrel{\mathbf{P}}{\to}0,
\end{equation*}
and 
\begin{equation*}
\sup_{t\in\R_+}
\frac{\abs{\epsilon^n_{\lfloor n^{2/3}t\rfloor\wedge n}}}{n^{1/3}}\,
\stackrel{\mathbf{P}}{\to}0,\qquad
 \sup_{t\in\R_+}
\frac{\abs{\epsilon^n_{\lfloor(nb_n)^{2/3}%
t\rfloor\wedge n}}}{n^{1/3}b_n^{4/3}}\,
\stackrel{\mathbf{P}^{1/b_n^2}}{\to}0\,.
\end{equation*}
\end{lemma}
\begin{proof}
We prove the convergences on the first line. By \eqref{eq:35} and \eqref{eq:1a}
\begin{equation}
  \label{eq:28}
    \sup_{t\in[0,1]}\abs{\overline{\epsilon}^n_t}\le 
\frac{1}{n}+\frac{1}{n}\sum_{k=1}^n\xi^n_{k,n-k+1}.
\end{equation}
The right-most convergence follows since $\mathbf{E}\xi^n_{k,n-k+1}=c_n/n$.
Next, by \eqref{eq:28} and the exponential Markov inequality for 
$\delta>0$ and $\lambda>0$
\begin{equation*}
\mathbf{P}\bl(  \sup_{t\in[0,1]}\abs{\overline{\epsilon}^n_t}>\delta\br)^{1/n}
\le e^{\lambda/n}\mathbf{E}e^{\lambda\xi^n_{1,1}}e^{-\lambda\delta}\to
e^{-\lambda\delta}\quad\text{ as }n\to\infty.
\end{equation*}
The left-most convergence in the statement of the theorem
follows since $\lambda$ is arbitrary.
Finally,
  \begin{equation*}
\mathbf{P}\Bl(  \sup_{t\in[0,1]}\frac{\sqrt{n}}{b_n}
\abs{\overline{\epsilon}^n_t}>\delta\Br)^{1/b_n^2}
\le e^{1/b_n^2}\bl(\mathbf{E}e^{\xi^n_{1,1}}\br)^{n/b_n^2}
e^{-\delta\sqrt{n}/b_n}\to 0
\quad\text{ as }n\to\infty,
\end{equation*}
proving the convergence in the middle.

The convergences on the second line are proved similarly.
\end{proof}
In the next three lemmas  we assume that $c>0$.
\begin{lemma}
  \label{le:min}
  \begin{trivlist}
  \item[    1.] The function $K_\rho(u),\,u\in[0,1],\,\rho\in\R_+,$ 
equals $0$ when either $u=0$ or $\rho=0$, 
is strictly decreasing, strictly concave and strictly 
subadditive in each of the
variables  $u$ and $\rho$ when the other variable is positive.
The function $L_c(u),\,u\in[0,1],$ equals $0$ at $u=0$ and
is strictly increasing  in $u$.
\item[ 2.]
If $u\in [(1-1/c)^+,1]$, then the function $K_c(x)+L_c(u+x)$ as a
function of $x$ is
strictly increasing for $x\in [0,1-u]$. If $c>1$ and $u\in[0,1-1/c)$,
then $K_c(x)+L_c(u+x)$  is  strictly increasing for 
$x\in[0,\tilde{u}]$, is strictly
decreasing  for $x\in[\tilde{u},u^\ast]$, and is strictly increasing for $x\in
[u^\ast,1-u]$, where $\tilde{u}\in[0,1-1/c-u]$ is the solution of the
equation 
\begin{equation*}
  \frac{x}{1-e^{-cx}}+x=1-u
\end{equation*}
and $u^\ast\in(1-1/c-u,\beta-u]$ is the solution of the equation
\begin{equation*}
   \frac{x}{1-e^{-cx}}=1-u.
\end{equation*}
 The values of the function  at $x=u^\ast$ and $x=0$ 
  coincide:  $K_c(u^\ast)+L_c(u+u^\ast)=L_c(u)$.
  \end{trivlist}
\end{lemma}
\begin{proof}
  Part 1 follows from the definitions. Part 2 follows by the equality
\begin{equation*}
  \frac{\partial}{\partial x}\bl(K_c(x)+L_c(u+x)\br)=
\bl(c(1-u-x)-\log(c(1-u-x))\br)
-\Bl(\frac{cx}{1-e^{-cx}}-\log\frac{cx}{1-e^{-cx}}\Br)
\end{equation*}
and the fact that the function $x-\log x$ is decreasing for
$x\in(0,1)$ and is increasing for $x>1$.
\end{proof}
Let $0\le s< t\le1$ and    $\Lambda_{s,t}$ denote the set of  absolutely   
continuous real-valued functions $\bx=(\bx_p,\,p\in[s,t])$ 
with  $\dot{\bx}_p\ge-1$  a.e. and 
$1-p-\bx_p\ge0$ on $[s,t]$. We denote for
$\bx\in\Lambda_{s,t}$
\begin{equation*}
  I^S_{s,t}(\bx)=\int_s^t\pi\Bl(\frac{\dot{\bx}_p+1}%
{c(1-p-\bx_p)}\Br)\,c(1-p-\bx_p)\,dp\,.
\end{equation*}
Let also for $0<\breve{s}<\breve{t}$, absolutely continuous 
real-valued $\bx=(\bx_p,\,p\in[\breve{s},\breve{t}])$, and
$\breve{\theta}\in\R$
\begin{equation*}
  \breve{I}^S_{\breve{s},\breve{t}}(\bx)=
\frac{1}{2}\int_{\breve{s}}^{\breve{t}}(\dot{\bx}_p+p-
\breve{\theta})^2\,dp.
\end{equation*}
\begin{lemma}
  \label{le:var} 
  \begin{enumerate}
  \item 
Given  $w\in(0,(t-s)^2/2)$,
the infimum of $I^S_{s,t}(\bx)$ over $\bx\in\Lambda_{s,t}$ such that 
$\bx_s=\bx_t=0$ and $\int_s^t\bx_p\,dp=w$ is attained at 
\begin{equation*}
\tilde{\bx}_p(s,t)=
s-p+\frac{t-s}{1-e^{-\tilde{\rho}(t-s)}}
\bl(1-e^{-\tilde{\rho}(p-s)}\br),\,
  p\in[s,t]\,,
\end{equation*}
where $\tilde{\rho}\in\R_+$ satisfies 
the equality $\partial K_{\rho}(t-s)/\partial \rho=-w$, i.e., 
\begin{equation*}
  \frac{\tilde{\rho}(t-s)}{1-e^{-\tilde{\rho}(t-s)}}=1+
\frac{w\tilde{\rho}}{t-s}+\frac{1}{2}\,\tilde{\rho}(t-s).
\end{equation*}
The value of the infimum equals 
$K_{\tilde{\rho}}(t-s)+(\tilde{\rho}-c)w+L_c(t)-L_c(s)=
\sup_{\rho\in\R_+}\bl(K_\rho(t-s)+(\rho-c) w\br)+L_c(t)-L_c(s)$.

If $w=0$, then the infimum is attained at
$\tilde{\bx}_p(s,t)=0,\,p\in[s,t],$ and is equal to $L_c(t)-L_c(s)$.
\item
Given $\breve{w}\in\R_+$,
the infimum of $\breve{I}^S_{\breve{s},\breve{t}}(\bx)$ 
over absolutely continuous real-valued functions
$\bx=(\bx_p,\,p\in[\breve{s},\breve{t}])$ such that 
$\bx_s=\bx_t=0$ and 
$\int_{\breve{s}}^{\breve{t}} \bx_p\,dp=\breve{w}$ is
attained at 
\begin{equation*}
\breve{\bx}_p(\breve{s},\breve{t})=
6\breve{w}\,\frac{(p-\breve{s})(\breve{t}-p)}{(\breve{t}-\breve{s})^3},\,
  p\in[\breve{s},\breve{t}]\,,
\end{equation*}
and equals
\begin{equation*}
  \frac{6\breve{w}^2}{(\breve{t}-\breve{s})^3}
-\breve{w}+\frac{(\breve{t}-\breve{\theta})^3-
(\breve{s}-\breve{\theta})^3}{6}.
\end{equation*}
  \end{enumerate}
\end{lemma}
\begin{proof}
Let $C$ denote the closed convex subset of the Banach space of real-valued
Lebesgue measurable
functions $h=(h_p,\,p\in[s,t])$ with norm $\norm{h}=\int_s^t\abs{
  h_p}\,dp$ specified by the conditions $h_p\ge0$ a.e.,
$\int_s^t h_p\,dp=t-s$, and $\int_s^t\int_s^ph_q\,dq\,dp=w+(t-s)^2/2$. 
We define a $[0,\infty]$-valued 
functional $F$ on $C$  by
\begin{equation*}
F(h)=\int_s^t\pi\Bl(\frac{h_p}{c\bl(1-s-\int_s^p h_q\,dq\br)}\Br)
\,c\bl(1-s-\int_s^p h_q\,dq\br)\,dp.
\end{equation*}
On noting that for $h\in C$
\begin{multline*}
  F(h)=
\int_s^t\bl(h_p\log\frac{
    h_p}{c}+c(1-s-\int_s^ph_q\,dq)\br)\,dp
+(1-t)\log(1-t)-(1-s)\log(1-s)\\
=\int_s^th_p\log    h_p\,dp+(t-s)\bl(c(1-s)-\log c\br)
-c\Bl(w+\frac{(t-s)^2}{2}\Br)+(1-t)\log(1-t)-(1-s)\log(1-s),
\end{multline*}
we see that  $F$ is strictly  convex on $C$. Therefore, 
the infimum of $F$ on $C$ is attained at a stationary point if the latter
exists. The method of Lagrange multipliers shows that
$  \tilde{h}_p=\bl(\tilde{\rho}(t-s)/(1-e^{-\tilde{\rho}(t-s)})\br)
e^{-\tilde{\rho}(p-s)}$ is such a point.
The assertion of part 1 of the
 lemma  follows.

For part 2 we apply the classical method of solving the isoperimetric
problem, see, e.g., Alekseev, Tikhomirov and Fomin 
\citeyear{AleTikFom87}.
\end{proof}
\begin{lemma}
  \label{le:optim}
  \begin{enumerate}
  \item 
Let $a\in[0,1]$ and $\tau\in[0,1]$. Then the infimum of
\begin{equation*}
\int_0^1\pi\Bl(\frac{1-\dot{\phi}_t}{c(1-t)}\Br)\,c(1-t)\,dt
\end{equation*}
over absolutely continuous non-decreasing functions
$\phi=(\phi_t,\,t\in[0,1])$ such that $\phi_0=0$, $\phi_1=a$, 
$\dot{\phi}_t\le1$ a.e., and the Lebesgue measure of the set
where $\dot{\phi}_t=0$  is at least $\tau$, equals
\begin{equation*}
  L_c((1-2a)\vee \tau)
+\frac{c}{2}(1-(1-2a)\vee\tau)^2
\;\pi\Bl(\frac{2(1-a-(1-2a)\vee\tau)}{%
 c(1-(1-2a)\vee\tau)^2}\Br).
\end{equation*}
\item
Let $\breve{\tau}\in\R_+$ and $\breve{\theta}\in\R$. Then the infimum of
$\int_0^\infty(-\dot{\phi}_t-\breve{\theta}+t)^2\,dt/2$
over absolutely continuous non-decreasing functions
$\phi=(\phi_t,\,t\in\R_+)$ such that $\phi_0=0$
 and the Lebesgue measure of the set
where $\dot{\phi}_t=0$  is at least $\breve{\tau}$, equals
$\bl((\breve{\tau}-\breve{\theta})^3\vee0+\breve{\theta}^3\br)/6$.
\end{enumerate}
\end{lemma}
\begin{proof}
We prove part 1.  The optimising integral 
can be written for a suitable function $g$ as
\begin{equation*}
  \int_0^1\bl( c(1-t)-\log(c(1-t))\br)\,dt+
\int_0^1 g(\dot{\phi}_t)\,dt
+\int_0^1\dot{\phi}_t\,\log(1-t)\,dt.
\end{equation*}
Let $\dot{\phi}^\ast$ denote the increasing rearrangement of $\dot{\phi}$
defined by $\dot{\phi}^\ast_t=\sup\{\lambda\in\R_+:\,\mu_\phi(\lambda)\le
  t\}$, where $\mu_\phi(\lambda)$ is the 
Lebesgue measure of those $t\in[0,1]$ for
  which $\dot{\phi}_t\le\lambda$. Since the function $\log(1-t)$ is
  decreasing, by  a  Hardy--Littlewood inequality, see Bennett and
  Sharpley \citeyear{BenSha88} or DeVore and Lorentz \citeyear{DeVLor93}, 
$\int_0^1\dot{\phi}_t\,\log(1-t)\,dt\ge 
\int_0^1\dot{\phi}^\ast_t\,\log(1-t)\,dt$. Also $\int_0^1
g(\dot{\phi}_t)\,dt=\int_0^1 g(\dot{\phi}^\ast_t)\,dt$.
  Therefore, the
function $\dot{\phi}$ can be assumed  non-decreasing, so
$\dot{\phi}_t=0$ for $t\in[0,\tau]$ and by the definition of $L_c$
\begin{equation}
     \label{eq:118p}
 \int_0^1\pi\Bl(\frac{1-\dot{\phi}_t}{c(1-t)}\Br)\,c(1-t)\,dt=
L_c(\tau)+ I(\phi,\tau),
\end{equation}
where 
\begin{equation}
  \label{eq:2420p}
  I(\phi,\tau)=\int_{\tau}^1
\pi\Bl(\frac{1-\dot{\phi}_t}{c(1-t)}\Br)\,c(1-t)\,dt\,.
\end{equation}
We now minimise $I(\phi,\tau)$
on the set $\Xi(\tau)$  of absolutely continuous
 functions $\phi$ with $\phi_\tau=0$, $\phi_1=a$,
 $\dot{\phi}_t\in [0,1]$ a.e., and
 $\dot{\phi_t}$ non-decreasing. 
Convexity considerations provide us with the lower bound
\begin{equation}
  \label{eq:66p}
I(\phi,\tau)\ge \frac{c(1-\tau)^2}{2}\:
\pi\biggl(\frac{2(1-\tau-a)}{c(1-\tau)^2}\biggr),
\end{equation}
which is attained at 
\begin{equation*}
\dot{\tilde{\phi}}_t=
1-\frac{2(1-\tau-a)}{(1-\tau)^2}\,(1-t),\;t\in[\tau,1].
\end{equation*}
If $\tau\ge 1-2a$, this function belongs to
$\Xi(\tau)$ and delivers the infimum to $I(\phi,\tau)$ on $\Xi(\tau)$
implying the required.

However, if $\tau<1-2a$ 
(hence, $2a<1$), then $\tilde{\phi}_t$ is negative for 
$t\in(\tau,2-(1-\tau)^2/(1-\tau-a)-\tau)$.
We prove that for those $\tau$ 
the infimum of $I(\phi,\tau)$ 
over $\phi\in \Xi(\tau)$   is attained at $\hat{\phi}$ defined by
 $\dot{\hat{\phi}}_t=0$ when $t\in[\tau,1-2a]$
and  $\dot{\hat{\phi}}_t=1-(1-t)/(2a)$ when 
$t\in[1-2a,1]$.
Let us consider $\dot{\phi}=(\dot{\phi}_t,\,t\in[\tau,1])$ for
 $\phi\in \Xi(\tau)$
as an element of the Banach space of Lebesgue measurable functions
$h=(h_t,\,t\in[\tau,1])$ with norm $\esssup_{t\in[\tau,1]}
\abs{h_t}$.
Let  functional $F$ on the subset of functions with $0\le h_t\le1$
a.e.  be defined by
$F(h)=\int_{\tau}^1\pi\bl((1-h_t)/(c(1-t))\br)\,c(1-t)\,dt.$ 
It is convex and has 
 a G\^{a}teau derivative at $\dot{\hat{\phi}}$  given by 
$\langle F'(\dot{\hat{\phi}}),h\rangle=
-\int_\tau^1\log\bl((1-\dot{\hat{\phi}}_t)/(c(1-t))\br)h_t\,dt$.
Therefore, for $\phi\in\Xi(\tau)$
\begin{multline*}
\langle  F'(\dot{\hat{\phi}}),\dot{\phi}-\dot{\hat{\phi}}\rangle=
\int_{\tau}^{1-2a}\log(c(1-t))\dot{\phi}_t\,dt
+\log(2ac)\int_{1-2a}^1(\dot{\phi}_t-\dot{\hat{\phi}}_t)\,dt\\
\ge\log(2ac)\int_\tau^{1-2a}\dot{\phi}_t\,dt+
\log(2ac)\int_{1-2a}^1(\dot{\phi}_t-\dot{\hat{\phi}}_t)\,dt =0,
\end{multline*}
implying (see, e.g., Ekeland and Temam \citeyear{EkeTem76}) that
 $I(\hat{\phi},\tau)\le I(\phi,\tau)$ for $\phi\in \Xi(\tau)$ as claimed. 
The definition of $\hat{\phi}$  and  \eqref{eq:2420p} 
yield $   I(\hat{\phi},\tau)=L_c(1-2a)-L_c(\tau)+2a^2c\,\pi\bl(1/(2ac)\br)$,
which in view of  \eqref{eq:118p} 
implies the assertion of the lemma for the case
$\tau< 1-2a$.

The proof of part 2 is similar, the infimum being attained at
$\breve{\phi}$ with $\dot{\breve{\phi}}_t=0$ for
$t\in[0,\tau\vee\breve{\theta}]$ and
$\dot{\breve{\phi}}_t=t-\breve{\theta}$ for $t>\tau\vee\breve{\theta}$.
\end{proof}

\begin{lemma}
  \label{le:comp}
Subsets $K$  of $\R_+^\N$ of sequences $\bu=(u_1,u_2,\ldots)$  such that 
$\sup_{\bu\in K}\sum_{i=1}^\infty u_i<\infty$
and $\lim_{i\to\infty}\sup_{\bu\in K}u_i=0$ are compact subsets 
of $\mathbb{S}$.
\end{lemma}
\begin{proof}
  It suffices to check sequential compactness. Let
  $\mathbf{u}^n,\,n\in\N,$ be a sequence of elements of
  $K$. The sequence
$\mathbf{u}^n,\,n\in\N,$ being compact for the product topology, let 
$\tilde{\mathbf{u}}=(\tilde{u}_1,\tilde{u}_2,\ldots)$ 
denote an accumulation point. Passing if necessary to a subsequence,
we may assume that $u^n_i\to \tilde{u}_i$ as $n\to\infty$ for
$i\in\N$. We have that $\tilde{\mathbf{u}}\in K$. Let $B=\sup_{\bu\in
  K}\sum_{i=1}^\infty u_i$. 
Given $\epsilon>0$, let $\delta>0$ be such that
$\chi(x)\le \epsilon x/(2B)$ 
for $x\in[0,\delta]$ (we use that $\chi(x)/x\to0$ as $x\to0$), 
let $k$ be such that $u_i\le
\delta\epsilon$ for $i\ge k$ and $\mathbf{u}\in K$, 
 and let $n_0$ be such that
$\abs{u^n_i-\tilde{u}_i}\le\delta\epsilon$ for $i=1,2,\ldots,k$ and $n\ge
n_0$. We then have that for $n\ge n_0$
\begin{equation*}
  \sum_{i=1}^\infty
  \chi\Bl(\frac{\abs{u^n_i-\tilde{u}_i}}{\epsilon}\Br)
\le\frac{1}{2B}\sum_{i=1}^\infty\abs{u^n_i-\tilde{u}_i}\le 1,
\end{equation*}
proving by $\epsilon$ being arbitrary that
$d_\chi(\mathbf{u}^n,\tilde{\mathbf{u}})\to 0$ as $n\to\infty$.
\end{proof}
\section{Large deviation asymptotics for the basic processes}
\label{sec:large-devi-asympt}
The main results of this  section are LDPs for
the stochastic processes  $\overline{S}^n$ and $\overline{E}^n$.
 We also give without proofs  LDPs for the
$\overline{\Phi}^n$ and $\overline{Q}^n$, which are not used further.
All these  processes are well-defined  random elements of
  $\mathbb{D}_C([0,1],\R)$.
For the notions and facts of 
idempotent probability theory  used extensively 
in the below argument, the reader
is referred to the appendix (or Puhalskii \citeyear{Puh01}).
\begin{theorem}
  \label{the:1}
Let $c_n\to c>0$ as $n\to\infty$.
Then the processes  $\overline{S}^n$ 
obey the LDP
for scale $n$ in $\mathbb{D}_C([0,1],\R)$ with 
action functional  $I^S$  given  by
\begin{equation*}
I^S(\bx)=\int_0^1\pi\Bl(\frac{\dot{\bx}_t+1}{c(1-t-\mathcal{R}(\bx)_t)}\Br)\, 
c\bl(1-t-\mathcal{R}(\bx)_t\br)\,dt\,
\end{equation*}
for absolutely continuous $\bx=(\bx_t,\,t\in[0,1])$ with
$\bx_0=0$,  $\dot{\bx}_t\ge
-1\,\text{a.e.}$, and $\mathcal{R}(\bx)_t\le1-t$ for $t\in[0,1],$  
and $I^S(\bx)=\infty$ for other $\bx$.
\end{theorem}
\begin{proof}
Let $A^n=(A^n_t,\,t\in[0,1])$  be defined by 
 \begin{equation}
   \label{eq:9}
     A^n_t=\frac{1}{n}\sum_{k=1}^{\lfloor nt\rfloor}
\sum_{j=1}^{n-Q^n_{k-1}-(k-1)}\xi^n_{kj}.
\end{equation}
 We note that by \eqref{eq:3} and
the definition of $\overline{S}^n_t$
\begin{equation}
  \label{eq:47}
\overline{S}^n_t=  A^n_t-\frac{\lfloor nt\rfloor}{n},\,t\in[0,1],
\end{equation}
so an LDP for the $\overline{S}^n$ would follow from an LDP for the $A^n$.
Let $\mathbf{e}=(t,\,t\in\R_+)$.
We prove that the $A^n$ as elements of $\D_C([0,1],\R)$
obey the LDP for scale
$n$ with action functional
\begin{equation*}
  I^{A}(\bx)=\int_0^1
\pi\Bl(\frac{\dot{\bx}_t}{c(1-t-\mathcal{R}(\bx-\mathbf{e})_t)}\Br)\,
c\bl(1-t-\mathcal{R}(\bx-\mathbf{e})_t\br)\,dt
\end{equation*}
if $\bx$ is absolutely continuous, $\bx_0=0$, $\dot{\bx}_t\ge0$ a.e.,
and $\mathcal{R}(\bx-\mathbf{e})_t\le 1-t$ for $t\in[0,1]$,
and $I^{A}(\bx)=\infty$ otherwise.

Let us extend  the time-domain of the processes $A^n$ to 
$\R_+$ by letting $A^n_t=A^n_{1}$ for $t\ge 1$. We
show that the extended $A^n$ 
  satisfy the hypotheses of
Theorem~5.1.5 in Puhalskii~\citeyear{Puh01}.
By \eqref{eq:9}  $A^n$ is a totally
discontinuous $\F^n$-adapted semimartingale with predictable measure of jumps
$\bl(\nu^n([0,t],\Gamma),\,t\in\R_+,\,\Gamma\in\mathcal{B}(\R)\br)$ given by
\begin{equation*}
  \nu^n([0,t],\Gamma)=\sum_{k=0}^{\lfloor
    n(t\wedge 1)\rfloor-1}F^{n}\Bl(1-\overline{Q}^n_{k/n}-
\frac{k}{n},\Gamma\setminus\{0\}\Br),\,\Gamma\in\mathcal{B}(\R),
\end{equation*}
where
\begin{equation}
  \label{eq:21}
      F^n(s,\Gamma')=\mathbf{P}\Bl(\frac{1}{n}\sum_{j=1}^{\lfloor
  ns\rfloor}\xi^n_{1j}\in \Gamma'\Br),\,s\in\R_+,\,\Gamma'\in\mathcal{B}(\R).
\end{equation}
Since the jumps of $A^n$ are bounded from above by $1$,
$A^n$ satisfies 
the Cram{\'e}r condition, so its stochastic
(or Dol{\'e}ans-Dade) exponential is well defined and has the form
\begin{equation}
  \label{eq:37}
    \mathcal{E}^n_t(\lambda)=\prod_{k=1}^{\lfloor n(t\wedge 1)\rfloor }
\Bl(1+\int_{\R}\bl(e^{\lambda
  x}-1\br)\,\nu^n\bl(\bl\{\frac{k}{n}\br\},dx\br)\Br)
=\prod_{k=0}^{\lfloor n(t\wedge 1)\rfloor-1 }
\int_{\R}e^{\lambda  x}\,F^n\bl(1-\overline{Q}^n_{k/n}-\frac{k}{n},\,dx\br),
\end{equation}
where $\lambda\in\R$.
 By  \eqref{eq:33} and \eqref{eq:47}
\begin{equation}
  \label{eq:4}
\overline{Q}^n=\mathcal{R}(A^n-\mathbf{e}^n+\overline{\epsilon}^n),
  \end{equation} 
where $\mathbf{e}^n =(\lfloor nt\rfloor/n,\,t\in\R_+)$.
Hence, recalling that the $\xi^n_i$ are Bernoulli and equal 1 with
probability $c_n/n$, we have by \eqref{eq:21} and \eqref{eq:37}
\begin{equation}
  \label{eq:56}
    \frac{1}{n}\log  \mathcal{E}^n_t(n\lambda)=
n\log\bl(1+(e^\lambda-1)\frac{c_n}{n}\br)
\int_0^{\lfloor n(t\wedge 1)\rfloor/n}
 \bl(1-\mathcal{R}(A^n-\mathbf{e}^n+\overline{\epsilon}^n)_s
-\frac{\lfloor ns
  \rfloor}{n}\br)\,ds.
\end{equation}
Let us  note that by the fact  that $Q^n_k+k\le
n$ and  \eqref{eq:33}
\begin{equation}
  \label{eq:48}
1-\mathcal{R}(A^n-\mathbf{e}^n-\overline{\epsilon}^n)_s-
\frac{\lfloor ns  \rfloor}{n}\ge 0 \text{ for }s\in[0,1].  
\end{equation}
 Thus, denoting for $\bx\in\D(\R_+,\R)$
  \begin{equation}
    \label{eq:12}
G_t(\lambda,\bx)=
c(e^\lambda-1)\int_0^{t\wedge 1}(1-
\mathcal{R}(\bx-\mathbf{e})_s-s)\,ds,
\end{equation}
 we conclude by \eqref{eq:56}, \eqref{eq:48},
the convergence $c_n\to c$, Lipshitz continuity of the
 reflection mapping  on $\D_C([0,1],\R_+)$, and
 Lemma~\ref{le:eps} that for arbitrary
$T>0$
\begin{equation*}
  \sup_{t\in[0,T]}\abs{\frac{1}{n}\log \mathcal{E}^n_t(n\lambda)-
G_t(\lambda,A^n)}\overset{\mathbf{P}^{1/n}}{\to}0
\text{ as $n\to\infty$.}
\end{equation*}
Since $G_t(\lambda,\bx)$ satisfies the uniform continuity and
majoration conditions of Theorem 5.1.5 of Puhalskii \citeyear{Puh01}, 
by the theorem
the sequence of laws of the $A^n$ on  $\D(\R_+,\R)$ 
is $\C$-exponentially
tight (of order $n$),
 and its every large deviation accumulation point solves the
maxingale problem $(0,G)$ with
$G=\bl(G_t(\lambda,\bx),\,t\in\R_+,\lambda\in\R,\bx\in\D(\R_+,\R)\br)$. 
Let   deviability $\mathbf{\Pi}^A$ on $\D(\R_+,\R)$ be a solution of $(0,G)$. 
We note  that 
$\mathbf{\Pi}^A\bl(\D(\R_+,\R)\setminus\C(\R_+,\R)\br)=0$ 
by the $\C$-exponential
 tightness  of the laws of the $A^n$. 
Let deviability $\hat{\mathbf{\Pi}}^A$ 
be the restriction of $\mathbf{\Pi}^A$ on $\C(\R_+,\R)$.
The  claimed LDP will follow if  for $\bx\in\C(\R_+,\R)$
\begin{equation}
  \label{eq:96}
  \hat{\mathbf{\Pi}}^{A}(\bx)=
  \begin{cases}
\exp(-I^{A}(p_1\bx))&\text{ if } \bx_{t}=\bx_{t\wedge 1},\,t\in\R_+,\\
\infty&\text{ otherwise.}    
  \end{cases}
\end{equation}
 The idea of the proof of \eqref{eq:96}
 is to translate the problem into a problem on
 uniqueness of  idempotent processes. 
Let $\Upsilon=\C(\R_+,\R)\times\C(\R_+,\R)$ and
 component idempotent processes 
$A=(A_t(\bx,\bx'),\,t\in\R_+,(\bx,\bx')\in\Upsilon)$
and $N=(N_t(\bx,\bx'),\,t\in\R_+,(\bx,\bx')\in\Upsilon)$ be
defined by the respective equalities $A_t(\bx,\bx')=\bx_t$
and $N_t(\bx,\bx')=\bx'_t$.  We will
 prove that  there exists
deviability $\mathbf{\Pi}$ on $\Upsilon$ such that $A$ and $N$  satisfy
 \begin{equation}
  \label{eq:18}
  A_t=N_{B_t(A)},\,t\in\R_+,\;\;\mathbf{\Pi}\text{-a.e.},
\end{equation}
where  
\begin{equation}
  \label{eq:71}
      B_t(\bx)=
  c\int_0^{t} (1-\mathcal{R}(\bx-\mathbf{e})_s-s)^+\,ds,
\end{equation}
$A$ has idempotent distribition $\hat{\mathbf{\Pi}}^A$ 
and $N$ is idempotent Poisson, i.e.,
$\sup_{\bx'\in \C(\R_+,\R)}\mathbf{\Pi}(\bx,\bx')
=\hat{\mathbf{\Pi}}^A(\bx)$ and 
$\sup_{\bx\in \C(\R_+,\R)}\mathbf{\Pi}(\bx,\bx')
=\mathbf{\Pi}^N(\bx')$, where $\mathbf{\Pi}^N$ is
the  Poisson deviability.
 After that we will draw on Ethier and Kurtz
\citeyear[Theorem 1.1, Chapter 6]{EthKur86} 
to conclude that \eqref{eq:18} has a
unique strong solution. That will imply that \eqref{eq:18} has a unique weak
   solution in the sense that the idempotent distribution  of $A$ is
specified uniquely and is given by \eqref{eq:96}.
 The reasoning used to establish \eqref{eq:18} is also
 along the lines of  the approaches developed in 
Ethier and Kurtz \citeyear{EthKur86}.

By \eqref{eq:48},
Lemma~\ref{le:eps}, and $\mathbf{\Pi}^A$ being an LD accumulation point of the
laws of the $A^n$, we have that
\begin{equation}
  \label{eq:97}
 1-\mathcal{R}(\bx-\mathbf{e})_s-s\ge0,\,s\in[0,1], 
\mathbf{\Pi}^A\text{-a.e.}, 
\end{equation}
so
\begin{equation}
  \label{eq:19}
  G_t(\lambda;\bx)=\tilde{G}_t(\lambda;\bx),
\,t\in\R_+,\;\;\hat{\mathbf{\Pi}}^A\text{-a.e.},
\end{equation}
where for $\lambda\in\R$
\begin{equation*}
  \tilde{G}_t(\lambda;\bx)=(e^\lambda-1)B_t(\bx),\,t\in\R_+,\,
\bx\in\C(\R_+,\R).
\end{equation*}
 Given $\epsilon>0$,
we define for $\bx,\bx'\in\C(\R_+,\R)$
\begin{equation}
\label{eq:gtil}
  G_t^\epsilon(\lambda;(\bx,\bx'))=\tilde{G}_t(\lambda;\bx)
+(e^{\lambda}-1)\epsilon t
\end{equation}
and introduce an idempotent process
$\hat{A}=(\hat{A}_t(\bx,\bx'),\,t\in\R_+,(\bx,\bx')\in\Upsilon)$ 
 by
\begin{equation}
  \label{eq:7}
  \hat{A}_t(\bx,\bx')=\bx_t+\bx'_{t}.
\end{equation}
As the deviability $\mathbf{\Pi}^A$ is a solution of the maxingale problem 
$(0,G)$, $\mathbf{\Pi}^A$ is concentrated on $\C(\R_+,\R)$, 
 $\mathbf{\Pi}^A$ and $\hat{\mathbf{\Pi}}^A$ coincide on 
$\C(\R_+,\R)$, and Lemma~\ref{le:expmax}  and
\eqref{eq:19} hold, it follows that the idempotent process $\bl(\exp(\lambda
\bx_t-\tilde{G}_t(\lambda;\bx)),\,t\in\R_+,\bx\in \C(\R_+,\R)\br)$ 
 is a  $\mathbf{C}$-uniformly maximable
 exponential maxingale on $(\C(\R_+,\R),\hat{\mathbf{\Pi}}^A)$, where
    $\mathbf{C}=(\mathcal{C}_t,\,t\in\R_+)$ is the canonical $\tau$-flow.
Next, the fact that $\bl(\exp(\lambda \bx_{ t}-
(e^{\lambda}-1) t),\,t\in\R_+,\bx\in \C(\R_+,\R)\br)$ 
 is a $\mathbf{C}$-exponential maxingale on
 $(\C(\R_+,\R),\mathbf{\Pi}^N)$ implies that
$\bl(\exp(\lambda \bx_{ t}-
(e^{\lambda}-1)\epsilon t),\,t\in\R_+,\bx\in \C(\R_+,\R)\br)$ 
 is a $\mathbf{C}$-exponential maxingale on
 $(\C(\R_+,\R),\mathbf{\Pi}^{N,\epsilon})$, where
$\mathbf{\Pi}^{N,\epsilon}((\bx_t,\,t\in\R_+))=
\mathbf{\Pi}^{N}((\bx_{t/\epsilon},\,t\in\R_+))$.
By Lemma~\ref{le:prodmax}, \eqref{eq:gtil} and \eqref{eq:7}
under product deviability $\hat{\mathbf{\Pi}}^A\times \mathbf{\Pi}^{N,\epsilon}$ the
idempotent process
$\bl(\exp(\lambda \hat{A}_t(\bx,\bx')-
G_t^\epsilon(\lambda;(\bx,\bx'))),
\,t\in\R_+,(\bx,\bx')\in\Upsilon\br)$ 
is an exponential maxingale relative to the $\tau$-flow
$\mathbf{A}=(\mathcal{A}_t,\,t\in\R_+)$, where
 $\mathcal{A}_t=\mathcal{C}_{ t}\otimes
\mathcal{C}_{ t}$. Let
\begin{equation}
  \label{eq:8}
    \sigma^\epsilon_t(\bx,\bx')=\inf\{s\in\R_+:\;B_s(\bx)+\epsilon s
\ge t\}.
\end{equation}
The idempotent variables $\sigma^\epsilon_t,\,t\in\R_+,$ are
 bounded idempotent 
$\mathbf{A}$-stopping times and
$G_{\sigma^\epsilon_t(\bx,\bx')}^\epsilon(\lambda;(\bx,\bx'))=
(e^\lambda-1)t$, so by
 Lemma~\ref{le:timechange} the idempotent process
$\bl(\exp(\lambda N^\epsilon_t(\bx,\bx')-
(e^\lambda-1)t),\,t\in\R_+,(\bx,\bx')\in\Upsilon\br)$, where
$N^\epsilon_t(\bx,\bx')=
 \hat{A}_{\sigma^\epsilon_t(\bx,\bx')}(\bx,\bx')$, is an
exponential maxingale on 
$(\Upsilon,\hat{\mathbf{\Pi}}^A\times \mathbf{\Pi}^{N,\epsilon})$
 relative to the $\tau$-flow $\mathbf{A}^\epsilon=
(\mathcal{A}_{\sigma^\epsilon_t},\,t\in\R_+)$. 
Hence, 
$N^\epsilon=(N^\epsilon_t(\bx,\bx'),\,t\in\R_+,(\bx,\bx')\in\Upsilon)$
is an $\mathbf{A}^\epsilon$-Poisson 
idempotent process, so it is a Poisson
idempotent process on 
$(\Upsilon,\hat{\mathbf{\Pi}}^A\times\mathbf{\Pi}^{N,\epsilon})$. In view of 
\eqref{eq:7}, \eqref{eq:8} and the
definition of $A_t$ we can  write that on $\Upsilon$
\begin{equation}
  \label{eq:13}
  A_t+\bx'_{ t}=N^\epsilon_{B_t(A)+\epsilon t},\,t\in\R_+.
\end{equation}
We now show that  \eqref{eq:18} is obtained as a limit of \eqref{eq:13}.
The pair $(A,N^\epsilon)$ specifies a mapping of $\Upsilon$ into
itself. Let $\mathbf{\Pi}^\epsilon$ denote the image of
$\hat{\mathbf{\Pi}}^A\times\mathbf{\Pi}^{N,\epsilon}$ 
under this mapping, i.e., 
$\mathbf{\Pi}^\epsilon(\bx,\bx')=\sup_{\substack{(\by,\by')\in\Upsilon:\,
A(\by,\by')=\bx,\\N^\epsilon(\by,\by')=\bx'}}\hat{\mathbf{\Pi}}^A(\by)
\mathbf{\Pi}^{N,\epsilon}(\by')$;
briefly, $\mathbf{\Pi}^\epsilon$ is the joint idempotent distribution of
$(A,N^\epsilon)$ on 
$(\Upsilon,\hat{\mathbf{\Pi}}^A\times\mathbf{\Pi}^{N,\epsilon})$. Since 
the idempotent distributions of $A$ and $N^\epsilon$ are deviabilities
and do not depend on $\epsilon$, 
the net
$\mathbf{\Pi}^\epsilon,\,\epsilon\to0,$ of deviabilities on $\Upsilon$
  is tight.
It is thus relatively compact for
weak convergence of idempotent probabilities.
Let $\mathbf{\Pi}$ denote an accumulation point 
of the $\mathbf{\Pi}^\epsilon$. 
By the continuous mapping theorem 
the marginal idempotent distributions 
of $\mathbf{\Pi}$ are equal to $\hat{\mathbf{\Pi}}^A$
and $\mathbf{\Pi}^N$:
$\sup_{\bx'\in\C(\R_+,\R)}\mathbf{\Pi}(\bx,\bx')
=\hat{\mathbf{\Pi}}^A(\bx)$ and 
$\sup_{\bx\in\C(\R_+,\R)}\mathbf{\Pi}(\bx,\bx')=\mathbf{\Pi}^N(\bx')$.
Next, by the definition of $\mathbf{\Pi}^\epsilon$, 
\eqref{eq:13}, and \eqref{eq:71}
 for $T>0$ and $\eta>0$, 
\begin{multline}\label{eq:75}
\mathbf{\Pi}^\epsilon\bl((\bx,\bx'):\,\sup_{t\in[0,T]}
\abs{\bx_t-\bx'_{B_t(\bx)}}\ge\eta\br)=
  (\hat{\mathbf{\Pi}}^A\times\mathbf{\Pi}^{N,\epsilon})\bl(\sup_{t\in[0,T]}
\abs{A_t-N^\epsilon_{B_t(A)}}\ge\eta\br)\\\le
(\hat{\mathbf{\Pi}}^A\times\mathbf{\Pi}^{N,\epsilon})
\Bl(\sup_{t\in[0,T]}\abs{\bx'_{ t}}\ge\frac{\eta}{2}\Br)\vee
(\hat{\mathbf{\Pi}}^A\times\mathbf{\Pi}^{N,\epsilon})
\Bl(\sup_{\substack{s,t\in[0,(c+\epsilon)T]:\\\abs{s-t}\le \epsilon T}}
\abs{N^\epsilon_{s}-N^\epsilon_{t}}\ge\frac{\eta}{2}\Bl)\\
=\mathbf{\Pi}^{N}\Bl(\sup_{t\in[0,\epsilon T]}\abs{\bx_{ t}}\ge\frac{\eta}{2}\Br)\vee
\sup_{\substack{s,t\in[0,(c+\epsilon)T]:\\\abs{s-t}\le \epsilon T}}
\mathbf{\Pi}^{N}\Bl(\abs{\bx_{s}-\bx_{t}}\ge\frac{\eta}{2}\Bl)=
\mathbf{\Pi}^{N}\Bl(\bx_{\epsilon T}\ge\frac{\eta}{2}\Br),
\end{multline}
where the latter two equalities use the definition of $\mathbf{\Pi}^{N,\epsilon}$,
the facts that $N^\epsilon$ is idempotent Poisson 
under $\hat{\mathbf{\Pi}}^A\times\mathbf{\Pi}^{N,\epsilon}$ and that
idempotent Poisson processes have stationary increments.
Given $L>0$, we have 
by an exponential Markov inequality and the fact that 
$\bl(\exp\bl(L \bx_{ t}-(e^L-1) t\br),\,t\in\R_+
\br)$ is an exponential maxingale under $\mathbf{\Pi}^{N}$
\begin{equation*}
  \mathbf{\Pi}^{N}\Bl(\bx'_{\epsilon T}\ge\frac{\eta}{2}\Br)\le
\mathbf{S}_{\mathbf{\Pi}^{N}}\bl(\exp(L\bx'_{\epsilon T})\br)
\exp\Br(-\frac{L\eta}{2}\Br)=\exp\Bl((e^L-1)\epsilon
T-\frac{L\eta}{2}\Br)\,,
\end{equation*}
where $\mathbf{S}_{\mathbf{\Pi}^{N}}$ denotes idempotent expectation
 with respect to $\mathbf{\Pi}^{N}$.
 Letting $\epsilon\to0$ and $L\to\infty$,
we conclude that
$\lim_{\epsilon\to0}\mathbf{\Pi}^{N}(\bx'_{\epsilon T}\ge\eta/2)=0$, 
so by \eqref{eq:75} 
$\lim_{\epsilon\to0}  \mathbf{\Pi}^\epsilon\bl((\bx,\bx'):\,\sup_{t\in[0,T]}
\abs{\bx_t-\bx'_{B_t(\bx)}}\ge\eta\br)=0.
$
Since the $\mathbf{\Pi}^\epsilon$ weakly converge along a subnet to $\mathbf{\Pi}$
and $\sup_{t\in[0,T]}\abs{\bx_t-\bx'_{B_t(\bx)}}$ is a continuous
function of $(\bx,\bx')\in\Upsilon$ so that the set
$\{(\bx,\bx')\in\Upsilon:\,\sup_{t\in[0,T]}
\abs{\bx_t-\bx'_{B_t(\bx)}}>\eta\}$ is open, we  conclude 
 that $  \mathbf{\Pi}\bl((\bx,\bx'):\,\sup_{t\in[0,T]}
\abs{\bx_t-\bx'_{B_t(\bx)}}>\eta\br)=0$. Consequently,
$\mathbf{\Pi}\bl((\bx,\bx'):\,\sup_{t\in[0,T]}
\abs{\bx_t-\bx'_{B_t(\bx)}}>0\br)=\sup_{\eta>0}  
\mathbf{\Pi}\bl((\bx,\bx'):\,\sup_{t\in[0,T]}
\abs{\bx_t-\bx'_{B_t(\bx)}}>\eta\br)=0,$
which is equivalent to \eqref{eq:18} by $A$ and $N$ being the
 first
and second component processes on $\Upsilon$, respectively.

Equation \eqref{eq:18} is of the form considered 
in Ethier and Kurtz \citeyear[Theorem 1.1, Chapter 6]{EthKur86}. 
The hypotheses of the theorem are seen to be met, which implies that
 \eqref{eq:18} has a unique (strong) solution for $A$ given by
$  A_t=N_{\sigma_t(N)},$ where
$\sigma_t(\bx')=
\inf\bl\{s\in[0,1]:\,
\int_0^s \bl(c (1-\mathcal{R}(\bx'-\mathbf{e})_p-p)^+\br)^{-1}\,dp\ge
 t\br\},\,\bx'\in\C(\R_+,\R)$.
Therefore, $\mathbf{\Pi}(\bx,\bx')=0$ if $(\bx_t,\,t\in\R_+)
\not=(\bx'_{\sigma_t(\bx')},\,t\in\R_+)$, so  the fact that
$\sup_{\bx\in\C(\R_+,\R)}\mathbf{\Pi}(\bx,\bx')=\mathbf{\Pi}^{N}(\bx')$ yields
$\mathbf{\Pi}(\bx,\bx')=\mathbf{\Pi}^{N}(\bx')$ if  $(\bx_t,\,t\in\R_+)
=(\bx'_{\sigma_t(\bx')},\,t\in\R_+)$.
Consequently, for $\bx\in\C(\R_+,\R)$
\begin{equation}
  \label{eq:16}
  \hat{\mathbf{\Pi}}^A(\bx)=\sup_{\bx'\in\C(\R_+,\R)}\mathbf{\Pi}(\bx,\bx')=
\sup_{\substack{\bx'\in\C(\R_+,\R):\\\,\bx_t
=\bx'_{\sigma_t(\bx')}}}
\mathbf{\Pi}^{N}(\bx')=
\sup_{\substack{\bx'\in\C(\R_+,\R):\\\,\bx_t=\bx'_{B_t(\bx)}}}
\mathbf{\Pi}^{N}(\bx').
\end{equation}
We have thus proved that $\hat{\mathbf{\Pi}}^A$ is uniquely specified by the
right-most side of \eqref{eq:16}. In particular, if
$\bx_t\not=\bx_{t\wedge 1}$ for some $t\in\R_+$, the set over which
the latter supremum is evaluated is empty, so $\hat{\mathbf{\Pi}}^A(\bx)=0$.
Let $\bx_t=\bx_{t\wedge 1},\,t\in\R_+$. 
Recalling that $\mathbf{\Pi}^{N}(\bx')=\exp(-I^{N}(\bx'))$, where
$  I^{N}(\bx')=\int_0^\infty\pi(\dot{\bx}'_t)\,dt
$ if $\bx'$ is absolutely continuous, $\bx'_0=0$, and
$\dot{\bx}'_t\ge0$ a.e., and
$I^{N}(\bx')=\infty$ otherwise, we derive by
a change of variables and \eqref{eq:71}
 that the right-most side of \eqref{eq:16} equals
$\exp(-I^A(p_1\bx))$ provided 
$1-\mathcal{R}(\bx-\mathbf{e})_s-s\ge0,\,s\in[0,1]$. If $\bx$ does not
meet the latter condition, then $\hat{\mathbf{\Pi}}^A(\bx)=0$ according to
\eqref{eq:97}. 
Equality \eqref{eq:96} has been proved, so the LDP for 
 the (extended) processes
$A^n$ has been proved.
By the contraction principle the (non-extended)
$A^n$ obey the LDP
in $\D_C([0,1],\R)$ with $I^A$. (Note that the $A^n$ are random elements of
$\D_C([0,1],\R)$.)  The LDP for the $\overline{S}^n$ follows by 
\eqref{eq:47} and the
contraction principle.
\end{proof}
\begin{corollary}
  \label{co:excess}
Let $c_n\to c>0$ as $n\to\infty$.
Then the processes  $(\overline{S}^n,\overline{E}^n)$ 
obey the LDP
for scale $n$ in $\mathbb{D}_C([0,1],\R^2)$ with 
action functional  $I^{S,E}$  given  by
\begin{equation*}
  I^{S,E}(\bx,\by)=I^S(\bx)+I^E_\bx(\by),
\end{equation*}
where $I^E_\bx(\by)=\int_0^1
\pi\bl(\dot{\by}_t/(c\mathcal{R}(\bx)_t)\br)\,c\,\mathcal{R}(\bx)_t\,dt$
if $\by=(\by_t,\,t\in[0,1])$ is non-decreasing and absolutely continuous
with $\by_0=0$ and $I^{E}_\bx(\by)=\infty$ otherwise.
\end{corollary}
\begin{proof}
Given a sequence $\bx^n,\,n\in\N$,   of elements of $\D_C([0,1],\R)$, let 
\begin{equation*}
\overline{E}'^n_t=\frac{1}{n}\sum_{i=1}^{\lfloor nt\rfloor}
\sum_{j=1}^{\lfloor 
n\mathcal{R}(\bx^n)_{(i-1)/n}\rfloor-1}\zeta^n_{ij},\,t\in[0,1].  
\end{equation*}
A standard argument (e.g.,  Theorem 2.3 in Puhalskii \citeyear{Puh94b})
 shows that if $\bx^n\to \bx$ as $n\to\infty$, then the
sequence $\overline{E}'^n,\,n\in\N$, obeys the LDP in $\D_C([0,1],\R)$ 
for scale $n$ with
action functional $I^E_\bx(\by),\,\by\in\D_C([0,1],\R)$. The claim now
follows by an argument as in
Puhalskii \citeyear[Theorem 2.2]{Puh95} (see also 
Chaganty \citeyear{Cha97}), \eqref{eq:33}, \eqref{eq:27}, and
Lemma~\ref{le:eps}.
\end{proof}
\begin{remark}\label{re:1}
  An application of the contraction principle 
    yields  LDPs for the 
$\overline{Q}^n$ and $\overline{\Phi}^n$.
\begin{enumerate}
\item
The processes  $\overline{Q}^n$ 
obey  the LDP for scale $n$ in  $\mathbb{D}_C([0,1],\R)$ with
action functional   $I^Q$  given  by
\begin{equation*}
  I^Q(\bx)=\int_0^1
\pi\Bl(\frac{\dot{\bx}_t+1}{c(1-t-\bx_t)}\Br)\,
c(1-t-\bx_t)\,\ind(\bx_t>0)\,dt 
+\int_0^{(1-1/c)^+}\pi\Bl(\frac{1}{c(1-t)}\Br)\,c(1-t)
\,\ind(\bx_t=0)\,dt
\end{equation*}
for absolutely continuous $\bx=(\bx_t,\,t\in[0,1])$  with
$\bx_0=0$, $\dot{\bx}_t\ge -1\,\text{a.e.}$ and $\bx_t\in[0,1-t]$
$t\in[0,1],$  
and $I^Q(\bx)=\infty$ for other $\bx$.
\item 
The processes $\overline{\Phi}^n$ 
obey  the LDP for scale $n$ in 
$\mathbb{D}_C([0,1],\R)$  with 
action functional   $I^\Phi$ given  by
 \begin{equation*}
   I^\Phi(\phi)=\int_0^1
\pi\Bl(\frac{1-\dot{\phi}_t}{c(1-t)}\Br)\,c(1-t)\,dt
+\sum K_c(l_i)
 \end{equation*}
if $\phi=(\phi_t,\,t\in[0,1])$ 
is absolutely continuous and non-decreasing, $\phi_0=0$ and
$\dot{\phi}_t\le 1$ a.e., where  the $l_i$  
are the lengths of the maximal  intervals where $\phi$ is constant and
 summation is performed over all such intervals, and
$I^\Phi(\phi)=\infty$ otherwise.
\end{enumerate}
\end{remark}
\section{Large deviations for connected components}
\label{sec:large-con}
In this section we prove Theorem~\ref{the:jointld}
and Corollaries~\ref{the:4} -- \ref{co:oc}.
We need the following lemma.
Let  $a\in[0,1]$, $m\in\N$, $u_1,\ldots,u_m$ be such that $u_i\in(0,1]$ 
 and $\sum_{i=1}^mu_i\le 1$, $r_1,\ldots,r_m$ belong to $\R_+$, and 
$\delta>0$.
 We denote by $B^n_{\delta}(a;\{u_i,r_i\}_{i=1}^m)$
 the event that
there exist $m$ connected components 
of $\mathcal{G}(n,c_n/n)$ of sizes in the intervals
$(n(u_i-\delta),n(u_i+\delta))$ for $i=1,2,\ldots,m$,  
the numbers of the excess edges of these components belong to the
 respective intervals $(n(r_i-\delta),n(r_i+\delta))$,
the other connected components are of sizes less  than $n\delta$, and the
total number of components of the random graph belongs to the interval
$(n(a-\delta),n(a+\delta))$. 
Let also $\tilde{B}^n_{\delta}(a)$ denote the event that 
all the  connected components are of sizes less  than $n\delta$ and the
total number of components  belongs to the interval
$(n(a-\delta),n(a+\delta))$.
\begin{lemma}
  \label{le:konechngia}
Let $c_n\to c>0$ as $n\to\infty$. If $\sum_{i=1}^mu_i\le 1-a$, then
\begin{equation*}
  \begin{split}
& \lim_{\delta\to0}\limsup_{n\to\infty}\frac{1}{n}\log
      \mathbf{P}\bl(B_{\delta}^n(a;\{u_i,r_i\}_{i=1}^m)\br)=
\lim_{\delta\to0}\liminf_{n\to\infty}\frac{1}{n}\log
      \mathbf{P}\bl(B_{\delta}^n(a;\{u_i,r_i\}_{i=1}^m)\br)
\\=&-\Biggl[\sum_{i=1}^m\sup_{\rho\in\R_+}
\Bl( K_{\rho}(u_i)+r_i\log\frac{\rho}{c}\Br)
+L_c\br((1-2a)\vee \sum_{i=1}^m u_i\bl)
\\&+
\frac{c}{2}\bl(1-(1-2a)\vee\sum_{i=1}^m u_i\br)^2
\;\pi\biggl(\frac{2\bl(1-a-(1-2a)\vee\sum_{i=1}^m u_i\br)}{%
 c\bl(1-(1-2a)\vee\sum_{i=1}^m u_i\br)^2}\biggr)\Biggr].
  \end{split}
\end{equation*}
If $\sum_{i=1}^mu_i> 1-a$, then
\begin{equation*}
  \lim_{\delta\to0}\limsup_{n\to\infty}\frac{1}{n}\log
      \mathbf{P}\bl(B_{\delta}^n(a;\{u_i,r_i\}_{i=1}^m)\br)=-\infty.
\end{equation*}
Also
\begin{multline*}
 \lim_{\delta\to0}\limsup_{n\to\infty}\frac{1}{n}\log
      \mathbf{P}\bl(\tilde{B}_{\delta}^n(a)\br)=
\lim_{\delta\to0}\liminf_{n\to\infty}\frac{1}{n}\log
      \mathbf{P}\bl(\tilde{B}_{\delta}^n(a)\br)\\
=-\Biggl[L_c\bl((1-2a)^+\br)+
\frac{c}{2}\bl(1-(1-2a)^+\br)^2\;
\pi\biggl(\frac{2\bl(1-a-(1-2a)^+\br)}{c\bl(1-(1-2a)^+\br)^2}\biggr)\Biggr].
 \end{multline*}
\end{lemma}
\begin{proof}We carry out the proof for the sets
      $B_{\delta}^n(a;\{u_i,r_i\}_{i=1}^m)$. A similar (and actually
      simpler) reasoning applies to the $\tilde{B}_{\delta}^n(a)$.
We denote throughout 
      $B_{\delta}^n(a;\{u_i,r_i\}_{i=1}^m)$ as $B^n_\delta$.
Upper bounds are addressed first. 
 Let $\delta\in(0,\min_{i=1,2,\ldots,m}u_i)$,
 $\sigma=(\sigma(1),\sigma(2),\ldots,\sigma(m))$ 
denote a permutation of the set
  $\{1,2,\ldots,m\}$ and $B'_{\delta,\sigma}$ denote the set
of functions $\bx\in\D_C([0,1],\R)$ with $\bx_0=0$
  such that  $\abs{\mathcal{T}(\bx)_1-a}\le\delta$ and
  there exist points $0=t'_0\le t'_1\le t'_2
\le\ldots\le t_{2m}'\le1=t'_{2m+1}$ with
$\abs{t_{2i}'-t_{2i-1}'-u_{\sigma(i)}}
\le\delta$ for $i=1,2\ldots,m$ for which
$\mathcal{R}(\bx)_{t_{2i-1}'}=
\mathcal{R}(\bx)_{t_{2i}'}=0$, $\mathcal{T}(\bx)_{t_{2i}'-}-
\mathcal{T}(\bx)_{t_{2i-1}'}=0$, 
and $\mathcal{R}(\bx)$ is not strictly positive on any
subinterval of  $[t_{2i}',t_{2i+1}']$ of length $\delta$ for
$i=0,1,\ldots,m$. Let $B_{\delta,\sigma}$ denote 
  the set of functions $(\bx,\by)\in\D_C([0,1],\R^2)$ 
 such that $\bx\in B'_{\delta,\sigma}$, $\by$ is non-decreasing 
with $\by_0=0$, and
$\abs{\by_{t_{2i}'}-\by_{t_{2i-1}'}-r_{\sigma(i)}}\le\delta$  for
$i=1,2\ldots,m$, where the $t_i'$ are associated with $\bx$, 
and let $B_\delta$ be the union of the
$B_{\delta,\sigma}$ over all permutations $\sigma$. 
By the construction of $Q^n$ and $E^n$, if there exists a connected
component of size $l$ of the random graph with $k$ excess edges, 
then there exist integers
 $k_1$ and $k_2$ ranging in $\{0,1,\ldots,n\}$ such
that $k_2-k_1=l$, $Q^n_{k_1}=Q^n_{k_2}=0$, 
$Q^n_i\ge1$ for $i=k_1+1,\ldots,k_2-1$, and
$E^n_{k_2}-E^n_{k_1}=k$. 
Also,  $\Phi^n$ does not increase on
$[k_1,k_2-1]$ and
$\Phi^n_n$ equals the  number of the connected
components of $\mathcal{G}(n,c_n/n)$.
 Therefore, recalling \eqref{eq:33} and \eqref{eq:33a}, we
 have that  $B^n_{\delta}\subset
\{(\overline{S}^n+\overline{\epsilon}^n,\overline{E}^n)
\subset B_{\delta}\}$.
Noting that  $B_{\delta}$ and its closure 
 in $\D_C([0,1],\R^2)$ have the same intersection with $\C([0,1],\R^2)$,
we have  by Corollary~\ref{co:excess} and Lemma~\ref{le:eps}
that
\begin{multline}\label{eq:51'}
  \limsup_{n\to\infty}\frac{1}{n}\log
  \mathbf{P}(B^n_{\delta})\le
\limsup_{n\to\infty}\frac{1}{n}\log
  \mathbf{P}\bl((\overline{S}^n+\overline{\epsilon}^n,\overline{E}^n)
\subset B_{\delta}\br)\\\le
-\inf_{(\bx,\by)\in B_{\delta}\cap\C([0,1],\R^2)}\bl(I^S(\bx)+I^E_\bx(\by)\br).
\end{multline}
Let $B'_\sigma$  denote the set
of functions $\bx\in\C([0,1],\R)$ with $\bx_0=0$
 such that $\mathcal{T}(\bx)_1=a$ and 
  there exist points $0=t_0\le t_1\le t_2\le\ldots\le t_{2m}\le
  t_{2m+1}=1$  with
$t_{2i}-t_{2i-1}=u_{\sigma(i)}$ for $i=1,2\ldots,m$ for which
$\mathcal{R}(\bx)_{t_{2i-1}}=
\mathcal{R}(\bx)_{t_{2i}}=0$, $\mathcal{T}(\bx)_{t_{2i}}=
\mathcal{T}(\bx)_{t_{2i-1}}$, and
$\mathcal{R}(\bx)$ equals zero on the intervals  
$[t_{2i},t_{2i+1}]$  for $i=0,1,\ldots,m$.
Let $\hat{B}_\sigma$ denote the set of functions
$(\bx,\by)\in\C([0,1],\R^2)$ such that $\bx\in B'_\sigma$, 
$\by$ is non-decreasing with $\by_0=0$ and
$\by_{t_{2i}}-\by_{t_{2i-1}}=r_{\sigma(i)}$  for
$i=1,2\ldots,m$ and the $t_i$ associated
with $\bx$.
Since $\cap_{\delta>0}B_{\delta}\cap\C([0,1],\R^2)=\cup_\sigma\hat{B}_\sigma$,
 we have by \eqref{eq:51'}
\begin{equation*}
\limsup_{\delta\to0}  \limsup_{n\to\infty}\frac{1}{n}\log
  \mathbf{P}(B^n_{\delta})\le
-\inf_\sigma\inf_{(\bx,\by)\in \hat{B}_\sigma}\bl(I^S(\bx)+I^E_\bx(\by)\br).
\end{equation*} 
As the function $\pi$ is convex and $\pi(1)=0$, it follows 
by the form of $I^E_\bx(\by)$ in Corollary~\ref{co:excess}
that the infimum of $I^E_\bx(\by)$ over
$\by$ such that $(\bx,\by)\in \hat{B}_\sigma$, where $\bx\in
B'_\sigma$ is fixed as well as the points $t_i$, 
is attained at $\hat{\by}$ defined by 
$\dot{\hat{\by}}_t=r_{\sigma(i)}
\mathcal{R}(\bx)_t/\int_{t_{2i-1}}^{t_{2i}}
\mathcal{R}(\bx)_s\,ds$ for $t\in[t_{2i-1},t_{2i}]$, where
$i=1,\ldots,m$, and $\dot{\hat{\by}}_t=c\mathcal{R}(\bx)_t$ 
elsewhere, and is equal to 
$\sum_{i=1}^m\pi\bl(r_{\sigma(i)}/(c\int_{t_{2i-1}}^{t_{2i}}
\mathcal{R}(\bx)_s\,ds)\br)c\int_{t_{2i-1}}^{t_{2i}}
\mathcal{R}(\bx)_s\,ds$. 
We can thus write
\begin{equation}
  \label{eq:14}
\limsup_{\delta\to0}  \limsup_{n\to\infty}\frac{1}{n}\log
  \mathbf{P}(B^n_{\delta})\le
-\inf_\sigma\inf_{\bx\in B'_\sigma}\biggl(I^S(\bx)
+\sum_{i=1}^m\pi\biggl(\frac{r_{\sigma(i)}}{%
 c\int_{t_{2i-1}}^{t_{2i}}
\mathcal{R}(\bx)_s\,ds}\biggr)\,c\int_{t_{2i-1}}^{t_{2i}}
\mathcal{R}(\bx)_s\,ds\biggr).
\end{equation} 
We now evaluate the infimum over $B'_\sigma$.
For $\bx\in B'_\sigma$ with $I^S(\bx)<\infty$, let
$\phi=(\phi_t,\,t\in[0,1])=\mathcal{T}(\bx)$. 
The condition $\dot{\bx}_t\ge-1$ a.e. 
implies that $\dot{\phi}_t\le 1$ a.e.
 The function $\phi$ does not increase on the
intervals $[t_{2i-1},t_{2i}],\,i=1,2,\ldots,m$, so $a=\phi_1=
\sum_{i=0}^m \int_{t_{2i}}^{t_{2i+1}}\dot{\phi}_t\,dt\le1- \sum_{i=1}^mu_i$
implying that $I^S(\bx)=\infty$ for  $\bx\in  B'_\sigma$ if
$\sum_{i=1}^m u_i>1-a$. This proves the second
limit in the statement of the lemma. In the rest of the argument we
assume that $\sum_{i=1}^m u_i\le1-a$.
We have on using that $\dot{\bx}_t\ge-1$ a.e. 
\begin{multline}
  \label{eq:15}
  \inf_{\bx\in B'_\sigma}\biggl(I^S(\bx)
+\sum_{i=1}^m\pi\Bl(\frac{r_{\sigma(i)}}{c\int_{t_{2i-1}}^{t_{2i}}
\mathcal{R}(\bx)_s\,ds}\Br)\,c\int_{t_{2i-1}}^{t_{2i}}
\mathcal{R}(\bx)_s\,ds\biggr)\\
=\inf_{\substack{w_i\in[0,u_i^2/2),\\i=1,2,\ldots,m}}
\biggl(  \inf_{\bx\in B'_\sigma(w_1,\ldots,w_m)}I^S(\bx)
+\sum_{i=1}^m\pi\Bl(\frac{r_i}{cw_i}\Br)
\,cw_i\biggr),
\end{multline}
where  $B'_\sigma(w_1,\ldots,w_m)=
\{\bx\in B'_\sigma:\,\int_{t_{2i-1}}^{t_{2i}}
\mathcal{R}(\bx)_s\,ds=w_{\sigma(i)},\,i=1,\ldots,m\}$.
We next prove that
  \begin{multline}
    \label{eq:1410}
\inf_{\bx\in B'_\sigma(w_1,\ldots,w_m)}I^S(\bx)
=\sum_{i=1}^m\sup_{\rho\in\R_+}\bl(K_{\rho}(u_i)+
(\rho-c)w_i\br)
+L_c\br((1-2a)\vee \sum_{i=1}^m u_i\bl)
\\+\frac{c}{2}\bl(1-(1-2a)\vee\sum_{i=1}^m u_i\br)^2
\;\pi\biggl(\frac{2\bl(1-a-(1-2a)\vee\sum_{i=1}^m u_i\br)}{%
 c\bl(1-(1-2a)\vee\sum_{i=1}^m u_i\br)^2}\biggr).
  \end{multline}
Since for $\bx\in B'_\sigma$ we have that
$\mathcal{R}(\bx)_{t_{2i-1}}=0$ and 
$\mathcal{T}(\bx)_{t_{2i}}=\mathcal{T}(\bx)_{t_{2i-1}}$ for
$i=1,2,\ldots,m$,  it follows
that $\mathcal{R}(\bx)_t=\bx_t-\bx_{t_{2i-1}}$ for
$t\in[t_{2i-1},t_{2i}]$. 
Hence, in view of  the form of $I^S$ in Theorem~\ref{the:1}
and Lemma~\ref{le:var}, if we change
 $\bx\in B'_\sigma(w_1,\ldots,w_m)$ 
with $I^S(\bx)<\infty$  on intervals $[t_{2i-1},t_{2i}]$  to
$\bl(\bx_{t_{2i-1}}+\tilde{\bx}_p(t_{2i-1},t_{2i}),
\,p\in[t_{2i-1},t_{2i}]\br)$, 
where $\tilde{\bx}_p(t_{2i-1},t_{2i})$ is defined in the statement 
 of Lemma~\ref{le:var}, this  will not increase the value of $I^S(\bx)$. The
altered function $\bx$ will still belong to
$B'_\sigma(w_1,\ldots,w_m)$ 
(note that $\phi$ is not affected by this modification of
$\bx$). Since $\bx_t+\phi_t=0$ on $\cup_{i=0}^m [t_{2i},t_{2i+1}]$,
the function
$\phi$ and the intervals $[t_{2i-1},t_{2i}]$ uniquely determine the
modified function $\bx$. We may thus optimise over $\phi$ and the 
$[t_{2i},t_{2i+1}]$, and assume    in view of Lemma~\ref{le:min},
Lemma~\ref{le:var}, Theorem~\ref{the:1},   and the
fact that $\dot{\phi}_t=0$ a.e. on $\cup_{i=1}^m[t_{2i-1},t_{2i}] $
that $\bx$ is such that
\begin{multline*}
  I^S(\bx)=\sum_{i=1}^m 
  \bl(\sup_{\rho\in\R_+}\bl(K_{\rho}(t_{2i}-t_{2i-1})+
(\rho-c)w_{\sigma(i)} \br)
+L_c(t_{2i})-L_c(t_{2i-1})\br)\\ +
\int_0^1\ind\bl(t\in\cup_{i=0}^m [t_{2i},t_{2i+1}]\br)
\pi\Bl(\frac{1-\dot{\phi}_t}{c(1-t)}\Br)\,c(1-t)\,dt
\\=\sum_{i=1}^m 
\sup_{\rho\in\R_+}\bl(  K_{\rho}(u_i)+(\rho-c)w_i\br)
+\int_0^1\pi\Bl(\frac{1-\dot{\phi}_t}{c(1-t)}\Br)\,c(1-t)\,dt\,,
\end{multline*}
where for the latter equality we used the definition of $L_c$ in 
\eqref{eq:l}.
An application of Lemma~\ref{le:optim} yields \eqref{eq:1410}.

Now,  a minimax argument
(cf., e.g., Aubin and Ekeland \citeyear{AubEke84}) shows that
  \begin{equation}
    \label{eq:131}
\inf_{w_i\in[0,u_i^2/2)}
\Bl(\sup_{\rho\in\R_+}\bl(K_{\rho}(u_i)+(\rho-c)w_i\br)+
\pi\Bl(\frac{r_i}{cw_i}\Br)cw_i\Br) =
\sup_{\rho\in\R_+}\Bl(K_{\rho}(u_i)+r_i\log\frac{\rho}{c}\Br).
\end{equation}
Thus,  by \eqref{eq:14}, \eqref{eq:15},  \eqref{eq:1410}, and \eqref{eq:131}, 
if  $\sum_{i=1}^mu_i\le 1-a$, then
\begin{multline*}
  \limsup_{\delta\to0}\limsup_{n\to\infty}\frac{1}{n}\log
      \mathbf{P}(B_{\delta}^n)\le-\Biggl[
\sum_{i=1}^m\sup_{\rho\in\R_+}
\Bl( K_{\rho}(u_i)+r_i\log\frac{\rho}{c}\Br)
+L_c\br((1-2a)\vee \sum_{i=1}^m u_i\bl)
\\+
\frac{c}{2}\bl(1-(1-2a)\vee\sum_{i=1}^m u_i\br)^2
\;\pi\biggl(\frac{2\bl(1-a-(1-2a)\vee\sum_{i=1}^m u_i\br)}{%
 c\bl(1-(1-2a)\vee\sum_{i=1}^m u_i\br)^2}\biggr)\Biggr].
\end{multline*}
We now establish the lower bound: if $\sum_{i=1}^mu_i\le 1-a$, then
\begin{multline}
  \label{eq:23}
      \liminf_{\delta\to0}\liminf_{n\to\infty}\frac{1}{n}\log
      \mathbf{P}\bl(B_{\delta}^n\br)\ge-\Biggl[
\sum_{i=1}^m\sup_{\rho\in\R_+}
\Bl( K_{\rho}(u_i)+r_i\log\frac{\rho}{c}\Br)
+L_c\br((1-2a)\vee \sum_{i=1}^m u_i\bl)
\\+
\frac{c}{2}\bl(1-(1-2a)\vee\sum_{i=1}^m u_i\br)^2
\;\pi\biggl(\frac{2\bl(1-a-(1-2a)\vee\sum_{i=1}^m u_i\br)}{%
 c\bl(1-(1-2a)\vee\sum_{i=1}^m u_i\br)^2}\biggr)\Biggr].
\end{multline}
Let $(\tilde{w}_i,\tilde{\rho}_i)$ denote 
the saddle point  of
 the function on the left-hand side of \eqref{eq:131} so that
 \begin{equation}
   \label{eq:22}
   K_{\tilde{\rho}_i}(u_i)+(\tilde{\rho}_i-c)\tilde{w}_i
+\pi\Bl(\frac{r_i}{c\tilde{w}_i}\Br)c\tilde{w}_i =
\sup_{\rho\in\R_+}\Bl(K_{\rho}(u_i)+r_i\log\frac{\rho}{c}\Br).
 \end{equation}
Calculations show that $\tilde{\rho}_i$ and $\tilde{w}_i$ are
specified by the equalities
\begin{equation}
  \label{eq:26}
      \frac{\tilde{\rho}_iu_i}{1-e^{-\tilde{\rho}_iu_i}}=
1+\frac{r_i}{u_i}+\frac{\tilde{\rho}_iu_i}{2},\;\;\tilde{w}_i=
\frac{r_i}{\tilde{\rho}_i}
\end{equation}
with $\tilde{\rho}_i=\tilde{w}_i=0$ if $r_i=0$.
 Let $s_0=0$, $s_i=\sum_{j=1}^i
u_j,\,j=1,\ldots,m$.
Motivated by the form of the optimal trajectory in
Lemma~\ref{le:var}, the  definition of $\hat{\by}$ above, and the
 definitions of $\tilde{\phi}$ and
$\hat{\phi}$ in the proof of Lemma~\ref{le:optim}, we  define  for 
$\eta\in(0,\min_{i=1,2,\ldots,m}u_i)$ 
an absolutely continuous function $\breve{\bx}^\eta$  by
$  \breve{\bx}^\eta_0=0$,
\begin{equation*}
\dot{\breve{\bx}}^\eta_t=
-1+\frac{u_i-\eta}{1-e^{-(\tilde{\rho}_i\vee\eta)(u_i-\eta)}}\,
(\tilde{\rho}_i\vee\eta)
\,e^{-(\tilde{\rho}_i\vee\eta)(t-s_{i-1})}
\end{equation*}
for $t\in(s_{i-1},s_i)$ and $i=1,\ldots,m$, 
$\dot{\breve{\bx}}^\eta_t=-\eta$
for $t\in\bl(s_m,s_m\vee(1-2a)\br)$, and
\begin{equation*}
\dot{\breve{\bx}}^\eta_t=-1+2\,\frac{1-s_m\vee(1-2a)-a}{(1-s_m\vee(1-2a))^2}\,
(1-t)
\end{equation*}
for $t\in\bl(s_m\vee(1-2a),1\br)$, 
and we define an absolutely continuous  function $\breve{\by}^\eta$ by
$\dot{\breve{\by}}^\eta_t=
r_i\mathcal{R}(\breve{\bx}^\eta)_t/\bl(\int_{s_{i-1}}^{s_i-\eta}
\mathcal{R}(\breve{\bx}^\eta)_s\,ds\br)$
for $ t\in(s_{i-1},s_i-\eta),\;\;i=1,\ldots,$ and
$\dot{\breve{\by}}^\eta_t=c\mathcal{R}(\bx)_t$ elsewhere.

Let us fix arbitrary $\delta\in(0,\min_{i=1,2,\ldots,m}u_i)$. 
For $\epsilon>0$, 
let $\breve{B}_{\epsilon,\eta}$ denote
the $\epsilon$-neighbourhood of $(\breve{\bx}^\eta,\breve{\by}^\eta)$ in
$\D_C([0,1],\R^2)$.
It follows from the definitions of $\breve{\bx}^\eta$,
$\breve{\by}^\eta$, and  the operator
$\mathcal{R}$ that if
$\epsilon$ and $\eta$ are small enough, then
for arbitrary $(\bx,\by)\in \breve{B}_{\epsilon,\eta}$
with $\bx_0=\by_0=0$ and $\by$  non-decreasing
there exist disjoint 
segments $(\tilde{s}_{i-1},\tilde{s}_i),\,i=1,2,\ldots,m$ 
with $\abs{\tilde{s}_i-\tilde{s}_{i-1}-u_i}<\delta$ such that
 the function $\mathcal{R}(\bx)$ is positive on these segments and
 equals zero at the endpoints,
 the other intervals where $\mathcal{R}(\bx)$ is positive
are of lengths  less
than $\delta$,  and 
$\abs{\by_{\tilde{s}_i}-\by_{\tilde{s}_{i-1}}-r_i}<\delta$. 
Furthermore, it may be assumed that
  $\mathcal{T}(\bx)_1\in(a-\delta,a+\delta)$.
We therefore have by \eqref{eq:33} and \eqref{eq:33a} that
$\{(\overline{S}^n+\overline{\epsilon}^n,\overline{E}^n)
\subset \breve{B}_{\epsilon,\eta}\}
\subset B^n_\delta$ for all small enough  $\epsilon$ and
$\eta$. As the set
$\breve{B}_{\epsilon,\eta}$ is open in
$\D_C([0,1],\R^2)$,
in view of  Lemma~\ref{le:eps} and Corollary~\ref{co:excess}
\begin{multline}
      \label{eq:44}
  \liminf_{n\to\infty}\frac{1}{n}\log
  \mathbf{P}(B^n_\delta)\ge\liminf_{n\to\infty}\frac{1}{n}\log
\mathbf{P}\bl((\overline{S}^n+\overline{\epsilon}^n,\overline{E}^n)
\subset \breve{B}_{\epsilon,\eta}\br)\\\ge
-\inf_{(\bx,\by)\in \breve{B}_{\epsilon,\eta}}(I^S(\bx)+I^E_\bx(\by))
\ge  -\bl(I^S(\breve{\bx}^\eta)+
I^E_{\breve{\bx}^\eta}(\breve{\by}^\eta)\br).
\end{multline}
By the definitions of $\breve{\bx}^\eta$ and $\breve{\by}^\eta$,
\eqref{eq:22}, \eqref{eq:26},
 the form of $I^S$ in Theorem~\ref{the:1},  the
form of $I^E$ in Corollary~\ref{co:excess},  part 1 of
 Lemma~\ref{le:var}, and part 1 of Lemma~\ref{le:optim}, we have
that $I^S(\breve{\bx}^\eta)+I^E_{\breve{\bx}^\eta}(\breve{\by}^\eta)$
 converges as $\eta\to0$ to the sum on
the right of \eqref{eq:23} which together with \eqref{eq:44} concludes
the proof of \eqref{eq:23}.
\end{proof}
\begin{proof}[Proof of Theorem~\ref{the:jointld}]
 We check that the sequence 
$(\alpha^n/n,\overline{U}^n,\overline{R}^n),\, n\in\N,$ is
  exponentially tight (of order $n$) 
in $[0,1]\times \mathbb{S}_1\times\mathbb{S}$.
By Lemma~\ref{le:comp}, the subsets  of $[0,1]\times
\mathbb{S}_1\times\mathbb{S}$ of elements $(a,\bu,\bfr)$, 
where $\bfr=(r_1,r_2,\ldots)$,
with the property that
$\sum_{i=1}^\infty r_i\le B$ for some  $B>0$ and 
$r_i\to0$ as $i\to\infty$ uniformly,
 are compact. Therefore, it suffices to check that
 \begin{align}
   \label{eq:129}
&\lim_{B\to\infty}\limsup_{n\to\infty}\mathbf{P}\Bl(\sum_{i=1}^\infty
\overline{R}^n_i>B\Br)^{1/n}=0,\\
&\lim_{i\to\infty}\limsup_{n\to\infty}\mathbf{P}\bl(
\sup_{j=i,i+1,\ldots}\overline{R}^n_j>\eta\br)^{1/n}=0,
\,\eta>0.   \label{eq:130}
 \end{align}
The first limit follows by exponential tightness 
of the $\overline{E}^n_1$ valid in view of Corollary~\ref{co:excess}
and the fact that $\sum_{i=1}^\infty\overline{R}^n_i= \overline{E}^n_1$.
For the second limit, we note that $\overline{R}^n_i$ equals the
increment of $\overline{E}^n_t$ over a time interval of length
$\overline{U}^n_i$, so for $\delta>0$
\begin{equation}
  \label{eq:118}
    \bigcup_{i=1}^\infty\{\overline{R}^n_{i}>
\eta,\,\overline{U}^n_{i}\le\delta\}\subset 
\bl\{\sup_{\substack{s,t\in[0,1]:%
\\\abs{s-t}\le\delta}}
\abs{\overline{E}^n_t-\overline{E}^n_s}>\eta\br\}.
\end{equation}
 Since $u_i\le 1/i$ for an element $\bu=(u_1,u_2,\ldots)$ of
$\mathbb{S}_1$, we have that 
  \begin{equation*}
\limsup_{i\to\infty}\limsup_{n\to\infty}\mathbf{P}\bl(
\sup_{j=i,i+1,\ldots}\overline{R}^n_j>\eta\br)^{1/n}\le
\limsup_{n\to\infty}\mathbf{P}\bl(\sup_{\substack{s,t\in[0,1]:%
\\\abs{s-t}\le\delta}}
\abs{\overline{E}^n_t-\overline{E}^n_s}>\eta\br)^{1/n}.
\end{equation*}
Therefore, \eqref{eq:130} follows on using
 that by $\C$-exponential tightness of the $\overline{E}^n$
\begin{equation}
  \label{eq:66}
\lim_{\delta\to0}
  \limsup_{n\to\infty}\mathbf{P}\bl(\sup_{\substack{s,t\in[0,1]:%
\\\abs{s-t}\le\delta}}
\abs{\overline{E}^n_t-\overline{E}^n_s}>\eta\br)^{1/n}=0.
\end{equation}

It thus remains to check that
  \begin{multline*}
    \lim_{\epsilon\to0}\limsup_{n\to\infty}\frac{1}{n}\log
      \mathbf{P}\Bl(d\Bl(\bl(\frac{\alpha^n}{n}
,\overline{U}^n,\overline{R}^n\br),
(a,\bu,\bfr)\Br)\le\epsilon\Br)\\= 
\lim_{\epsilon\to0}\liminf_{n\to\infty}\frac{1}{n}\log
 \mathbf{P}\Bl(d\Bl(\bl(\frac{\alpha^n}{n},\overline{U}^n,\overline{R}^n\br),
(a,\bu,\bfr)\Br)\le\epsilon\Br)=-I_c^{\alpha,U,R}(a,\bu,\bfr),
  \end{multline*}
where $d$ is a product metric on $[0,1]\times\mathbb{S}_1\times\mathbb{S}$ and 
$(a,\bu,\bfr)\in[0,1]\times \mathbb{S}_1\times\mathbb{S}$.
Let  $\bu=(u_1,u_2,\ldots)$ and $\bfr=(r_1,r_2,\ldots)$.
If all the $u_i>0$, then given  $\delta>0$,    for
all small enough $\epsilon>0$ and all large enough $m$
\begin{equation}
  \label{eq:57}
  \Bl\{d\Bl(\bl(\frac{\alpha^n}{n},\overline{U}^n,\overline{R}^n\br),
(a,\bu,\bfr)\Br)\le\epsilon\Br\}\subset B^n_{\delta}(a;\{u_i,r_i\}_{i=1}^m).
\end{equation}
If  $u_1>0$ and $u_i=0$ for all large $i$, then \eqref{eq:57} holds
for $m$ that is  the greatest index $i$ with $u_i>0$.
If $u_1=0$, then  we have the inclusion
\begin{equation*}
\Bl\{d\Bl(\bl(\frac{\alpha^n}{n},\overline{U}^n,\overline{R}^n\br),
(a,\bu,\bfr)\Br)\le\epsilon\Br\}\subset \tilde{B}^n_{\delta}(a).  
\end{equation*}
Therefore,  Lemma~\ref{le:konechngia} and the form of
 $I_c^{\alpha,U,R}(a,\bu,\bfr)$ imply that, provided $r_i=0$ when $u_i=0$,
 \begin{equation}
   \label{eq:132}
     \limsup_{\epsilon\to0}\limsup_{n\to\infty}\frac{1}{n}\log
\mathbf{P}\Bl(d\Bl(\bl(\frac{\alpha^n}{n},\overline{U}^n,\overline{R}^n\br),
(a,\bu,\bfr)\Br)\le\epsilon\Br)
\le-I_c^{\alpha,U,R}(a,\bu,\bfr).
\end{equation}
If for some $i$ we have that $u_i=0$ and $r_i>0$, then 
by \eqref{eq:118} and \eqref{eq:66} the left-hand side of
\eqref{eq:132} equals $-\infty$, so the required inequality holds as well.

For the lower bound
\begin{equation}
  \label{eq:133}
  \liminf_{\epsilon\to0}\liminf_{n\to\infty}\frac{1}{n}\log
\mathbf{P}\Bl(d\Bl(\bl(\frac{\alpha^n}{n},\overline{U}^n,\overline{R}^n\br),
(a,\bu,\bfr)\Br)\le\epsilon\Br)\ge-I_c^{\alpha,U,R}(a,\bu,\bfr)
  \end{equation}
we may assume that $r_i=0$ when $u_i=0$.
Let us be given $\epsilon>0$ and $B>0$.
If all the $u_i$ are positive, then for  all small enough $\delta>0$,
$\eta>0$, and large enough $m$
 we have the inclusion
 \begin{equation*}
B^n_{\delta}(a;\{u_i,r_i\}_{i=1}^m)\subset
\Bl\{d\bl(\bl(\frac{\alpha^n}{n},\overline{U}^n,\overline{R}^n\br),
(a,\bu,\bfr)\br)\le \epsilon\Br\}\cup\Bl\{\sum_{i=1}^\infty
\overline{R}^n_i> B\Br\}\cup\Bl\{\sup_{i=m+1,\ldots}\overline{R}^n_i>
\eta\Br\}.   
 \end{equation*}
To see the latter we use the inequality
$\sum_{i=m+1}^\infty \chi(u'_i/\epsilon)\le
\sup_{i=m+1,\ldots}\bl(\chi(u'_i/\epsilon)/u'_i\br)\sum_{i=m+1}^\infty
u_i'$ for $(u_1',u_2',\ldots)\in\mathbb{S}$ 
and  the convergence $\chi(x)/x\to0$ as $x\to0$.
Lemma~\ref{le:konechngia}, \eqref{eq:129}, and \eqref{eq:130} imply
\eqref{eq:133}. 
If $u_1>0$ and not all the $u_i$ are positive, then by a similar argument
\begin{equation*}
B^n_{\delta}(a;\{u_i,r_i\}_{i=1}^m)\subset
\Bl\{d\bl(\bl(\frac{\alpha^n}{n},\overline{U}^n,\overline{R}^n\br),
(a,\bu,\bfr)\br)\le \epsilon\Br\}\cup\Bl\{\sum_{i=1}^\infty
\overline{R}^n_i> B\Br\}\cup\bigcup_{i=m+1}^\infty\{\overline{R}^n_{i}>
\eta,\,\overline{U}^n_{i}\le\delta\}, 
\end{equation*}
 where $m$ is the greatest index $i$ with
$u_i>0$.
If $u_1=0$, then
\begin{equation*}
\tilde{B}^n_{\delta}(a)\subset
\Bl\{d\bl(\bl(\frac{\alpha^n}{n},\overline{U}^n,\overline{R}^n\br),
(a,\bu,\bfr)\br)\le \epsilon\Br\}\cup\Bl\{\sum_{i=1}^\infty
\overline{R}^n_i> B\Br\}\cup
\bigcup_{i=1}^\infty\{\overline{R}^n_{i}>
\eta,\,\overline{U}^n_{i}\le\delta\}.  
\end{equation*}
In either case, \eqref{eq:133}  follows by 
Lemma~\ref{le:konechngia}, \eqref{eq:129},  \eqref{eq:118}, and
 \eqref{eq:66}.
  \end{proof}
Corollaries \ref{the:4} and \ref{the:5} follow by an application
  of the contraction
principle. In some more detail, the infima of
$I_c^{\alpha,U,R}(a,\bu,\bfr)$ and   $I_c^{\alpha,U}(a,\bu)$
over $a\in[0,1]$ are attained  at $a^\ast=1/(2c)$ if $\sum_{i=1}^\infty
u_i< 1-1/c$ and at $a^\ast=1-\sum_{i=1}^\infty
u_i-c(1-\sum_{i=1}^\infty u_i)^2/2$ if $\sum_{i=1}^\infty u_i\ge
  1-1/c$; the infimum of $I_c^{\alpha,U,R}(a,\bu,\bfr)$ over
  $\bfr\in\mathbb{S}$ 
is found by a minimax argument (it is actually attained at
  $\bfr^\ast=(r^\ast_1,r^\ast_2,\ldots)$ with
  $r^\ast_i=cu_i^2/(1-\exp(-cu_i))-cu_i^2/2-u_i$), 
cf., Aubin and Ekeland \citeyear{AubEke84}.
The expression for $I_c^\alpha(a)$ is obtained on noting that
subadditivity of $K_c(u)$ in $u$ implies that $\sum_{i=1}^\infty K_c(u_i)\ge 
K_c(\sum_{i=1}^\infty u_i)$, so  one should minimise
  $I^{\alpha,U}_c(a,\bu)$  with respect to
$\sum_{i=1}^\infty u_i$, and that $K_c(u)$ is monotonically
decreasing in $u$, so 
the infimum can be taken over $\sum_{i=1}^\infty u_i\ge
1-2a$. We  provide  more detail as to the proofs of
Corollaries~\ref{the:2} and \ref{co:oc}.
\begin{proof}[Proof of Corollary~\ref{the:2}]
 Let $A_\delta(u_1,\ldots,u_m)$ for $\delta\in(0,\min_{i=1,\ldots,m}u_i/2)$
 denote the subset of $\mathbb{S}_1$  of vectors 
$\tilde{\bu}=(\tilde{u}_1,\tilde{u}_2,\ldots)$ such that 
there exist distinct 
$j_i\in\{1,2,\ldots,\lfloor 2/u_i\rfloor\}$ with
$\abs{\tilde{u}_{j_i}-u_i}<\delta$ for $i=1,2,\ldots, m$.
Let a set $A(u_1,\ldots,u_m)$ be defined as 
 the set of
$\tilde{\bu}=(\tilde{u}_1,\tilde{u}_2,\ldots)\in \mathbb{S}_1$ such that 
$\tilde{u}_{j_i}=u_i,\,i=1,2,\ldots,m,$ for some
$j_1,\ldots,j_m$.
Since $A(u_1,\ldots,u_m)$ equals the 
intersection of the closures of the $A_\delta(u_1,\ldots,u_m)$ over
$\delta>0$,
the sets $A_\delta(u_1,\ldots,u_m)$ are open in $\mathbb{S}_1$, and
$A^n_\delta(u_1,\ldots,u_m)=\{\overline{U}^n\in
A_\delta(u_1,\ldots,u_m)\}$, we have by 
Corollary~\ref{the:5} and the definition of the LDP 
\begin{multline*}
    \lim_{\delta\to0}\limsup_{n\to\infty}\frac{1}{n}\log
      \mathbf{P}\bl(A_{\delta}^n(u_1,\ldots,u_m)\br)= 
\lim_{\delta\to0}\liminf_{n\to\infty}\frac{1}{n}\log
      \mathbf{P}\bl(A_{\delta}^n(u_1,\ldots,u_m)\br)\\=
-\inf_{\bu\in A(u_1,\ldots,u_m)}I_c^U(\bu).
\end{multline*}
We evaluate the latter infimum. Since $I_c^U(\bu)$ is invariant with
respect to permutations of the entries of $\bu$, we may replace
$\bu$ with its permutation that has $u_1,\ldots,u_m$ as the first
$m$ entries. By subadditivity of $K_c(u)$ in $u$ we have that 
$\sum_{i=m+1}^{\infty}K_c(u_i)\ge 
K_c\bl(\sum_{i=m+1}^\infty u_i\br)$, so it is optimal to assume
that $u_{m+2}=u_{m+3}=\ldots=0$. We thus need to find  optimal
$u_{m+1}$. If $\sum_{i=1}^m u_i\ge 1-1/c$, then 
$I_c^U(\bu)=\sum_{i=1}^{m+1}K_c(u_i)+L_c\bl(\sum_{i=1}^{m+1}u_i\br)$.
By Lemma~\ref{le:min} 
$K_c(u_{m+1})+L_c\bl(\sum_{i=1}^{m+1}u_i\br)>L_c\bl(\sum_{i=1}^{m}u_i\br)$
for any $u_{m+1}>0$, so it is optimal to take $u_{m+1}=0$, accordingly
$\inf_{\bu\in A(u_1,\ldots,u_m)}I_c^U(\bu)=\sum_{i=1}^{m}K_c(u_i)
+L_c\bl(\sum_{i=1}^{m}u_i\br)$.
If $\sum_{i=1}^m u_i< 1-1/c$, then  Lemma~\ref{le:min} implies that 
for $u^\ast>0$ such that 
$u^\ast/(1-\exp(-cu^\ast))=1-\sum_{i=1}^mu_i$ we have
$K_c(u^\ast)+L_c\bl(\sum_{i=1}^{m}u_i+u^\ast\br)=
L_c\bl(\sum_{i=1}^{m}u_i\br)$. Also  $\sum_{i=1}^m
u_i+u^\ast>1-1/c$, so the
 choice of $u^\ast$ as $u_{m+1}$ yields the value of the action functional
$\sum_{i=1}^{m}K_c(u_i)+
K_c(u^\ast)+L_c\bl(\sum_{i=1}^{m}u_i+u^\ast\br)=
\sum_{i=1}^{m}K_c(u_i)+L_c\bl(\sum_{i=1}^{m}u_i\br)$.
If $u_{m+1}\not=u^\ast$ and is such that  $\sum_{i=1}^{m+1}u_i\ge1-1/c$,
then
$I_c^U(\bu)=\sum_{i=1}^{m+1}K_c(u_i)+L_c\bl(\sum_{i=1}^{m+1}u_i\br)$,
which is greater than
$\sum_{i=1}^{m}K_c(u_i)+L_c\bl(\sum_{i=1}^{m}u_i\br)$ by Lemma~\ref{le:min}.
Finally, if $u_{m+1}$ is such that $\sum_{i=1}^{m+1}u_i<1-1/c$, then
with the use of  Lemma~\ref{le:min}
$I_c^U(\bu)=\sum_{i=1}^{m+1}K_c(u_i)+L_c\bl(1-1/c\br)>
\sum_{i=1}^{m+1}K_c(u_i)+L_c\bl(\sum_{i=1}^{m+1}u_i\br)\ge
\sum_{i=1}^{m}K_c(u_i)+L_c\bl(\sum_{i=1}^{m}u_i\br)$.
 Therefore, $u^\ast$
is the optimal value of $u_{m+1}$. Thus, $\inf_{\bu\in
  A(u_1,\ldots,u_m)}I_c^U(\bu)=\sum_{i=1}^mK_c(u_i)
+L_c\bl(\sum_{i=1}^mu_i\br)$ and it is attained at a unique point
$\bu^\ast$ given by
$\bu^\ast=(u_1,u_2,\ldots,u_m,0,0,\ldots)$ if $\sum_{i=1}^mu_i\ge 1-1/c$
and 
$\bu^\ast=(u_1,u_2,\ldots,u_m,u^\ast,0,0,\ldots)$ 
if $\sum_{i=1}^mu_i<1-1/c$.
We also have by the form of $I^{U,R}_c$ in Corollary~\ref{the:4}
 that the infimum of
$I^{U,R}_c(\bu^\ast,\bfr)$ over $\bfr$ equals $I^{U}_c(\bu^\ast)$
and is attained at the unique point
$\bfr^\ast=(r^\ast_1,\ldots,r^\ast_m,0,0,\ldots)$ if
$\sum_{i=1}^mu_i\ge 1-1/c$ and 
$\bfr^\ast=(r^\ast_1,\ldots,r^\ast_m,r^\ast,0,0,\ldots)$ if
$\sum_{i=1}^mu_i< 1-1/c$. 
Therefore, letting $\tilde{d}$ 
denote a metric on $\mathbb{S}_1\times\mathbb{S}$,
\begin{align*}
        \lim_{\delta\to0}\limsup_{n\to\infty}\frac{1}{n}\log
      \mathbf{P}\bl(A_\delta^n(u_1,\ldots,u_m)\br)&=
\lim_{\delta\to0}\liminf_{n\to\infty}\frac{1}{n}\log
      \mathbf{P}\bl(A_\delta^n(u_1,\ldots,u_m)\br)\\
=\lim_{\eta\to0}\limsup_{n\to\infty}\frac{1}{n}\log
      \mathbf{P}\bl(\tilde{d}\bl((\overline{U}^n,\overline{R}^n),
(\bu^\ast,\bfr^\ast)\br)<\eta\br)&= 
\lim_{\eta\to0}\liminf_{n\to\infty}\frac{1}{n}\log
\mathbf{P}\bl(\tilde{d}\bl((\overline{U}^n,\overline{R}^n),
(\bu^\ast,\bfr^\ast)\br)<\eta\br)\\
&=-\Bl(\sum_{i=1}^mK_c(u_i)+
L_c\bl(\sum_{i=1}^mu_i\br)\Br).
\end{align*}
In addition, 
$\liminf_{\delta\to0} \lim_{n\to\infty}
\mathbf{P}\bl(\{\tilde{d}\bl((\overline{U}^n,\overline{R}^n),
(\bu^\ast,\bfr^\ast)\br)<\eta\}|
A_{\delta}^n(u_1,\ldots,u_m)\br)=1
$ for $\eta>0$ as 
in  Freidlin and Wentzell \citeyear[Theorem 3.4 of Chapter 3]{wf2}. 
The proof is completed by noting that 
 $\{\tilde{d}\bl((\overline{U}^n,\overline{R}^n),
(\bu^\ast,\bfr^\ast)\br)<\eta\}\subset
\tilde{A}_{\delta,\epsilon}^n(u_1,\ldots,u_m)
\subset A_\delta^n(u_1,\ldots,u_m)$ 
for all small enough $\eta>0$.
\end{proof}
\begin{proof}[Proof of Corollary~\ref{co:oc}]
    By Theorem~\ref{the:jointld} and the contraction principle
  \begin{equation}
    \label{eq:59}
    I_c^{\alpha,\beta,\gamma}(a,u,r)=\inf_{(\bu,\bfr)\in
    O(u,r)}I_c^{\alpha,U,R}(a,\bu,\bfr),
  \end{equation}
where $O(u,r)=\{(\bu,\bfr)\in
\mathbb{S}_1\times\mathbb{S}:\,u_1=u,\,r_1=r\}$. The assertion of the corollary
for $u=0$ follows. Let us
assume now that $u>0$. The infimum of $\sup_{\rho\in\R_+}\bl( K_{\rho}(x)+
r\log(\rho/c)\br)$ over $r\in\R_+$ equals $K_c(x)$, therefore, it suffices
 to minimise over $u_2,u_3,\ldots$ the function
  \begin{equation*}\sum_{i=1}^\infty K_c(u_i)+
L_c\bl((1-2a)\vee\sum_{i=1}^\infty  u_i\br)
+\frac{c}{2}\bl(1-(1-2a)\vee\sum_{i=1}^\infty u_i\br)^2
\;\pi\biggl(\frac{2\bl(1-a-(1-2a)\vee\sum_{i=1}^\infty u_i\br)}{%
 c\bl(1-(1-2a)\vee\sum_{i=1}^\infty u_i\br)^2}\biggr)
\end{equation*}
By the fact that $K_c(x)<0$ for $x>0$  and is
decreasing in $x$ (Lemma~\ref{le:min}), we can
assume that in an optimal configuration $\sum_{i=1}^\infty u_i\ge
1-2a$. Next, since $K_c(x)$ is concave in $x$, $K_c(0)=0$ and $u_i\le
u$, we have that
\begin{equation}
  \label{eq:5}
  \sum_{i=1}^\infty K_c(u_i)\ge
  \Bigl\lfloor\frac{\sum_{i=1}^\infty u_i}{u}\Bigr\rfloor K_c(u)
+K_c\Bl(\sum_{i=1}^\infty u_i
-u\Bigl\lfloor\frac{\sum_{i=1}^\infty u_i}{u}\Bigr\rfloor\Br).
\end{equation}
Hence, by Theorem~\ref{the:jointld}
\begin{multline}
      \label{eq:20}
    I_c^{\alpha,\beta,\gamma}(a,u,r)=\sup_{\rho\in\R_+}\Bl( K_{\rho}(u)+
r\log\frac{\rho}{c}\Br)-K_c(u)+\inf_{\tau\in[(1-2a)\vee u,1-a]}\biggl(
 \Bigl\lfloor\frac{\tau}{u}\Bigr\rfloor K_c(u)\\
+K_c\Bl(\tau-u\Bigl\lfloor\frac{\tau}{u}\Bigr\rfloor\Br)
+L_c( \tau)+\frac{c}{2}(1-\tau)^2
\;\pi\Bl(\frac{2(1-a-\tau)}{c(1-\tau)^2}\Br)\biggr),
\end{multline}
as required. Part 1 has been proved.

We prove part 2.
By the contraction principle the sequence $(\beta^n/n,\gamma^n/n),\,n\in\N,$
obeys the LDP for scale $n$ with action functional
$I^{\beta,\gamma}_c(u,r)=\inf_{a\in[0,1]}I_c^{\alpha,\beta,\gamma}(a,u,r)$, 
which yields
the assertion of part 2 for $(u,r)=(0,0)$. Let $u>0$.
The infimum of the right-most term on the right of \eqref{eq:20} over 
$a\in[(1-\tau)/2,1-\tau]$ is attained at $(1-\tau)/2$ if $\tau<1-1/c$
and at $1-\tau-c(1-\tau)^2/2$ if $\tau\ge 1-1/c$ with respective values
$c(1-\tau)^2/2\;\pi\bl(1/( c(1-\tau))\br)$ and $0$.
If $\tau<1-1/c$, then by Lemma~\ref{le:min} there exists 
$\tau^\ast\in(0,1-\tau)$
such that $\tau+\tau^\ast> 1-1/c$ and $L_c(\tau)=K_c(\tau^\ast)+
L_c(\tau+\tau^\ast)$. Therefore, in analogy with \eqref{eq:5}
\begin{equation*}
  \Bigl\lfloor\frac{\tau}{u}\Bigr\rfloor K_c(u)
+K_c\Bl(\tau-u\Bigl\lfloor\frac{\tau}{u}\Bigr\rfloor\Br)
+L_c( \tau)\ge \Bigl\lfloor\frac{\tau+\tau^\ast}{u}\Bigr\rfloor K_c(u)
+K_c\Bl(\tau+\tau^\ast-u\Bigl\lfloor\frac{\tau+\tau^\ast}{u}\Bigr\rfloor\Br)
+L_c( \tau+\tau^\ast),
\end{equation*}
which implies that we may disregard the domain $\tau<1-1/c$. Hence,
\eqref{eq:20} yields
\begin{equation}
\label{eq:41}
I_c^{\beta,\gamma}(u,r)=\sup_{\rho\in\R_+}\Bl( K_{\rho}(u)+
r\log\frac{\rho}{c}\Br)-K_c(u)+\inf_{\tau\in
[(1-1/c)\vee u,1]}\biggl(
  \Bigl\lfloor\frac{\tau}{u}\Bigr\rfloor K_c(u)
+K_c\Bl(\tau-u\Bigl\lfloor\frac{\tau}{u}\Bigr\rfloor\Br)
+L_c( \tau)\biggr).
\end{equation}
If $u\ge1-1/c$, then 
for $\tau\ge  u$  by Lemma~\ref{le:min}
$
\lfloor\tau/u\rfloor K_c(u)+K_c(\tau-\lfloor\tau/u\rfloor u)
+L_c( \tau)\ge K_c(u)+K_c(\tau-u)+L_c(\tau)\ge   K_c(u)+L_c( u)$,
so 
\begin{equation}
  \label{eq:24}
I_c^{\beta,\gamma}(u)=\sup_{\rho\in\R_+}\Bl( K_{\rho}(u)+
r\log\frac{\rho}{c}\Br)+L_c(u).  
\end{equation}
Let us now assume that  $u<1-1/c$, so $c>1$. If $\tau\ge
\lfloor(1-1/c)/u\rfloor u+u$, then by the fact that
$\lfloor(1-1/c)/u\rfloor u+u>1-1/c$ and Lemma~\ref{le:min} 
\begin{multline*}
  \Bl\lfloor\frac{\tau}{u}\Br\rfloor K_c(u)+
K_c\Bl(\tau-\Bl\lfloor\frac{\tau}{u}\Br\rfloor u\Br)+L_c(\tau)
\ge
\biggl(\Bl\lfloor\frac{1}{u}\Bl(1-\frac{1}{c}\Br)\Br\rfloor+1\biggr) K_c(u)+
K_c\Br(\tau-\Bl\lfloor\frac{1}{u}\Bl(1-\frac{1}{c}\Br)\Br\rfloor
u-u\Br)+L_c(\tau)\\\ge
\biggl(\Bl\lfloor\frac{1}{u}\Bl(1-\frac{1}{c}\Br)\Br\rfloor+1\biggr)
K_c(u)+L_c\Br(\Bl\lfloor\frac{1}{u}\Bl(1-\frac{1}{c}\Br)\Br\rfloor
u+u\Br) ,
\end{multline*}
so by \eqref{eq:41}
\begin{equation}
  \label{eq:40}
  I_c^\beta(u)=\sup_{\rho\in\R_+}\Bl( K_{\rho}(u)+
r\log\frac{\rho}{c}\Br)-K_c(u)+\inf_{\tau\in
[1-1/c,\lfloor(1-1/c)/u\rfloor u+u]}\biggl(
  \Bigl\lfloor\frac{\tau}{u}\Bigr\rfloor K_c(u)
+K_c\Bl(\tau-u\Bigl\lfloor\frac{\tau}{u}\Bigr\rfloor\Br)
+L_c( \tau)\biggr).
\end{equation}
By subadditivity of $K_c(x)$ in $x$, for $\tau\ge1-1/c$
\begin{equation}
  \label{eq:58}
  \Bl\lfloor\frac{\tau}{u}\Br\rfloor K_c(u)+
K_c\Bl(\tau-\Bl\lfloor\frac{\tau}{u}\Br\rfloor u\Br)
\ge
\Bl\lfloor\frac{1}{u}\Bl(1-\frac{1}{c}\Br)\Br\rfloor K_c(u)+
K_c\Br(\tau-\Bl\lfloor\frac{1}{u}\Bl(1-\frac{1}{c}\Br)\Br\rfloor u\Br).
\end{equation}
By Lemma~\ref{le:min} and the definition of $\hat{u}$ 
for $\tau\in[1-1/c,\lfloor(1-1/c)/u\rfloor u+u]$
\begin{equation*}
K_c\Br(\tau-\Bl\lfloor\frac{1}{u}\Bl(1-\frac{1}{c}\Br)\Br\rfloor u\Br)
+L_c( \tau)\ge
K_c(\hat{u}\wedge u)
+L_c\Bl(\Bl\lfloor\frac{1}{u} \Bl(1-\frac{1}{c}\Br)\Br\rfloor u
+\hat{u}\wedge u \Br),
\end{equation*}
which implies by  \eqref{eq:58} that the minimum in 
  \eqref{eq:40} is attained at 
$\hat{\tau}=\lfloor(1-1/c)/u\rfloor u+\hat{u}\wedge u$
 completing the proof of part 2.

Part 3 follows by minimising $I_c^{\beta,\gamma}(u,r)$ over $r\in\R_+$.
\end{proof}

\section{Normal and moderate deviations for the largest component}
\label{sec:norm-moder-devi}
In this section we prove Theorems~\ref{the:7} and \ref{the:6}.
We start by establishing a law-of-large-numbers result.
Let  
\begin{align}
  \label{eq:69}
    \overline{M}^n_t&=\frac{1}{n}
\sum_{i=1}^{\lfloor nt\rfloor}
\sum_{j=1}^{n-Q^n_{i-1}-(i-1)}\bl(\xi^n_{ij}-\frac{c_n}{n}\br),\;
t\in[0,1],\\
  \label{eq:29}
    \overline{L}^n_t&=\frac{1}{n}
\sum_{i=1}^{\lfloor nt\rfloor}
\sum_{j=1}^{Q^n_{i-1}-1}\bl(\zeta^n_{ij}-\frac{c_n}{n}\br),\;
t\in[0,1],
\end{align}
so that by \eqref{eq:2}, \eqref{eq:3},  \eqref{eq:1a}, and \eqref{eq:27}
  \begin{align}
    \label{eq:72}
      \overline{Q}^n_{t}&=
\int_0^{\lfloor nt\rfloor/n}\Bl(c_n\Bl(1-\overline{Q}^n_{s}-
\frac{\lfloor ns\rfloor}{n}\Br)-1\Br)\,ds
+\overline{\epsilon}^n_{t}+\overline{M}^n_t
+\overline{\Phi}^n_t\,,\\
    \label{eq:52}
    \overline{E}^n_{t}&=c_n
\int_0^{\lfloor
  nt\rfloor/n}\overline{Q}^n_{s}\,ds
+\overline{L}^n_t-\frac{c_n}{n}\int_0^{\lfloor  nt\rfloor/n}
\ind(\overline{Q}^n_s>0)\,ds\,.
  \end{align}
The processes
$\overline{M}^n=(\overline{M}^n_t,\,t\in[0,1])$ and
$\overline{L}^n=(\overline{L}^n_t,\,t\in[0,1])$   are
orthogonal square integrable martingales relative to the filtration 
$(\mathcal{F}^n_t,\,t\in[0,1])$ with respective predictable quadratic
characteristics
\begin{align}
  \label{eq:73}
  \langle\overline{M}^n\rangle_t&=\frac{c_n}{n}\bl(1-\frac{c_n}{n}\br) 
\int_{0}^{\lfloor nt\rfloor/n}
\Bl(1-\overline{Q}^n_s-\frac{\lfloor ns\rfloor}{n}\Br)\,ds\,,\\
  \label{eq:70}
    \langle\overline{L}^n\rangle_t&=\frac{c_n}{n}\bl(1-\frac{c_n}{n}\br) 
\int_{0}^{\lfloor nt\rfloor/n}\Bl(\overline{Q}^n_{s}-\frac{1}{n}\Br)^+\,ds\,.
\end{align}
Let functions 
$\overline{q}=(\overline{q}_t,\,t\in[0,1])$,
$\overline{\phi}=(\overline{\phi}_t,\,t\in[0,1])$, and 
$\overline{e}=(\overline{e}_t,\,t\in[0,1])$
be  defined by
\begin{align}
  \label{eq:42}
\overline{q}_t&=
  \begin{cases}
1-t-e^{-ct}&\,\text{ if }t\in[0,\beta]\,,\\
    0&\,\text{ otherwise}\,,
  \end{cases}\\
  \label{eq:31}
      \overline{\phi}_t&=
\begin{cases}
\displaystyle
\frac{c}{2}(t^2-\beta^2)-(c-1)(t-\beta)&\,\text{ if }t\in[\beta,1]\,,\\
0&\,\text{ otherwise}\,,
\end{cases}
  \end{align}
and 
\begin{equation}
  \label{eq:76}
  \overline{e}_t=e^{-c(t\wedge\beta)}
-1+c(t\wedge\beta)-\frac{c(t\wedge\beta)^2}{2}\,.
\end{equation}
Equivalently, the pair $(\overline{q},\overline{\phi})$ can be defined
as the solution to the Skorohod problem 
\begin{equation}
  \label{eq:74}
  \overline{q}_t=\int_0^t\bl(c(1-\overline{q}_s-s)-1\br)\,ds+
\overline{\phi}_t\;\text{ and }\;
\overline{\phi}_t=\int_0^t\ind(\overline{q}_s=0)\,d\overline{\phi}_s.
\end{equation}
We note that
\begin{align}
  \label{eq:32}  
    \overline{q}_t=
\int_0^t\bl(c(1-\overline{q}_s-s)-1\br)\,ds\;\text{ for
}t\in[0,\beta]&
\;\text{ and }\;\overline{e}_t=c\int_0^t\overline{q}_s\,ds
\;\text{ for }t\in[0,1]\,.
\end{align}

\begin{lemma}
  \label{le:lln}Let $c_n\to c>0$ as $n\to\infty$.
Then the processes $\overline{Q}^n$,
  $\overline{\Phi}^n$, and $\overline{E}^n$
 converge in probability uniformly on $[0,1]$ to
the functions  $\overline{q}$, $\overline{\phi}$, and $\overline{e}$
 respectively.
\end{lemma}
\begin{proof}
By \eqref{eq:73}, \eqref{eq:70}, 
and Doob's inequality the
$\overline{M}^n$ and $\overline{L}^n$
 converge to $0$ in probability uniformly over $[0,1]$ as
$n\to\infty$. Also, the $\overline{\epsilon}^n$ converge in probability
to $0$ uniformly on $[0,1]$ by Lemma~\ref{le:eps}. 
Now, a standard tightness argument applied to \eqref{eq:72} and
 \eqref{eq:52}  shows that the
sequence $(\overline{Q}^n,\overline{\Phi}^n,\overline{E}^n),\,
n\in\N,$ is $\C$-tight in $\D_C([0,1],\R^3)$,
where a limit point 
$(\tilde{q},\tilde{\phi},\tilde{e})$
is such that 
$  \tilde{q}_t=\int_0^t\bl(c(1-\tilde{q}_s-s)-1\br)\,ds+\tilde{\phi}_t,
$
$\tilde{\phi}$  is
non-decreasing with $\tilde{\phi}_t=\int_0^t\ind(\tilde{q}_s=0)
\,d\tilde{\phi}_s$, and $\tilde{e}_t=c\int_0^t\tilde{q}_s\,ds$. 
Hence, $(\tilde{q},\tilde{\phi},\tilde{e})=
(\overline{q},\overline{\phi},\overline{e})$
concluding the proof.
\end{proof}
\begin{remark}
  The convergences $\overline{Q}^n\overset{\mathbf{P}}{\to}
\overline{q}$ and $\overline{\Phi}^n\overset{\mathbf{P}}{\to}
\overline{\phi}$  also follow from
  Remark~\ref{re:1} since the action functionals
$I^Q$ and $I^\Phi$ are equal 
to $0$ at $\overline{q}$ and $\overline{\phi}$, respectively.
\end{remark}

We now prove a diffusion limit theorem, which will lead to
the proof of Theorem~\ref{the:7}. Let us define
processes $M^n=(M^n_t,\,t\in[0,1])$, $L^n=(L^n_t,\,t\in[0,1])$, 
$X^n=(X^n_t,\,t\in[0,1])$, 
$Y^n=(Y^n_t,\,t\in[0,1])$, and $Z^n=(Z^n_t,\,t\in[0,1])$   by
the respective equalities 
 $M^n_t=\sqrt{n}\:\overline{M}^n_t$,  $L^n_t=\sqrt{n}\:\overline{L}^n_t$,
$  X^n_t=\sqrt{n}(\overline{Q}^n_t-\overline{q}_t)$,
 $  Y^n_t=\sqrt{n}(\overline{\Phi}^n_t-\overline{\phi}_t)$, and 
$  Z^n_t=\sqrt{n}(\overline{E}^n_t-\overline{e}_t)$. 
By \eqref{eq:72}, \eqref{eq:52}, \eqref{eq:42}, 
\eqref{eq:74}, and \eqref{eq:32} these processes satisfy the equations
\begin{align}
  \label{eq:77a}
  X^n_{t}&=-c_n\int_0^tX^n_s\,ds+
\sqrt{n}(c_n-c)\int_0^t\sigma_s^2\,ds
+M^n_t+\tilde{\epsilon}^n_t+Y^n_t,\\
  \label{eq:80}
Z^n_t&=c_n\int_0^tX^n_s\,ds+\sqrt{n}(c_n-c)\int_0^t\overline{q}_s\,ds+
L^n_t+\tilde{\delta}^n_t,
  \end{align}
where 
\begin{align}
  \label{eq:6.4b}
\sigma^2_t&=
\begin{cases}
\displaystyle  e^{-ct}&\text{ if }t\in[0,\beta]\,,\\
\displaystyle 1-t&\text{ if }t\in[\beta,1]\,,
\end{cases}\\
\label{eq:6.9a}
\tilde{\epsilon}^n_t&=\sqrt{n}\overline{\epsilon}^n_t+
\sqrt{n}\int_t^{\lfloor
  nt\rfloor/n}\Bl(c_n\Bl(1-\overline{Q}^n_{s}-\frac{\lfloor
  ns\rfloor}{n}\Br)-1\Br)\,ds+\sqrt{n}\int_0^t\Bl(s-\frac{\lfloor
  ns\rfloor}{n}\Br)\,ds,\\
  \label{eq:81}
\tilde{\delta}^n_t&=\sqrt{n}\,c_n\int_t^{\lfloor
  nt\rfloor/n}\overline{Q}^n_s\,ds
-\frac{c_n}{\sqrt{n}}\int_0^{\lfloor  nt\rfloor/n}
\ind(\overline{Q}^n_s>0)\,ds\,.
\end{align}
We note that 
Lemma~\ref{le:eps} implies
 that if $c_n\to c$ as $n\to\infty$, then 
for arbitrary $\eta>0$
\begin{align}
  \label{eq:tiep}
\sup_{t\in[0,1]}\abs{\tilde{\epsilon}^n_t}&\overset{\mathbf{P}}{\to}0,
\\  \intertext{also}
  \label{eq:tide}
   \sup_{t\in[0,1]}\abs{\tilde{\delta}^n_t}&\le\frac{c_n}{\sqrt{n}}.
\end{align}
Let
$W^{(1)}=(W^{(1)}_t,\,t\in[0,1])$ and $W^{(2)}=(W^{(2)}_t,\,t\in[0,1])$ 
be independent  Wiener processes, and processes
$H=(H_t,\,t\in[0,1])$ and $Z=(Z_t,\,t\in[0,1])$
 be specified  by the equations
\begin{align}
\label{eq:6.11a}
  H_t&=-c\int_0^{t\wedge\beta}H_s\,ds+\theta\int_0^t\sigma^2_s\,ds+
\sqrt{c}\int_0^t\sigma_s\,dW^{(1)}_s,\\
Z_t&=c\int_0^{t\wedge\beta}H_s\,ds+\theta\int_0^t\overline{q}_s\,ds+
\sqrt{c}\int_0^t\sqrt{\overline{q}_s}\,dW^{(2)}_s\,.
  \label{eq:90}
  \end{align}
We also define  processes $M=(M_t,\,t\in[0,1])$ and $L=(L_t,\,t\in[0,1])$ by
 $M_t=\sqrt{c} \int_0^t\sigma_s\,dW_s^{(1)}$ and  $L_t=\sqrt{c}
 \int_0^t\sqrt{\overline{q}_s}\,dW_s^{(2)}$.
\begin{lemma}
  \label{the:diff}Let $\sqrt{n}(c_n-c)\to\theta\in\R$ as $n\to\infty$,
  where $c>0$. Then
  \begin{equation*}
\lim_{B\to\infty}\limsup_{n\to\infty}
\mathbf{P}(\sup_{t\in[0,1]}\abs{X^n_t}>B)=0\,.
  \end{equation*}
Also the following holds.
  \begin{enumerate}
  \item If $\beta>0$, then
 for   $\delta\in(0,\beta\wedge(1-\beta))$ the processes $M^n$, $L^n$,
$(X^n_t,\,t\in[0,\beta-\delta])$, 
$(Y^n_t,\,t\in[\beta+\delta,1])$, and $(Z^n_t,\,t\in[0,1])$
 jointly converge in distribution in
$\D_C([0,1],\R^2)\times\D_C([0,\beta-\delta],\R)\times
\D_C([\beta+\delta,1],\R)\times\D_C([0,1],\R)$
 to the respective processes
  $M$, $L$, $(H_t,\,t\in[0,\beta-\delta])$, 
$(-H_t,\,t\in[\beta+\delta,1])$, and $Z$.
In addition, 
$    \lim_{n\to\infty}
\mathbf{P}(\sup_{t\in[0,\beta-\delta]}\abs{Y^n_t}>\delta)=0\,.
$
\item If $\beta=0$, then the processes 
$Y^n$  converge in distribution in
$\D_C([0,1],\R)$
to the process $-H$.
  \end{enumerate}
\end{lemma}
\begin{proof}
We start by proving that 
the processes $(M^n,L^n)$ converge in distribution
in $\D_C([0,1],\R^2)$ to the process $(M,L)$.
The processes $M^n$ and $L^n$   are orthogonal square
integrable martingales relative to the filtration
$(\mathcal{F}^n_t,\,t\in[0,1])$, whose respective
 predictable quadratic characteristics
$n\langle\overline{M}^n\rangle_t$ and
$n\langle\overline{L}^n\rangle_t$   converge in
probability as $n\to\infty$  to $c\int_0^t\sigma^2_s\,ds$
and $c\int_0^t\overline{q}_s\,ds$ respectively
in view of \eqref{eq:73}, \eqref{eq:70}, \eqref{eq:42}, \eqref{eq:6.4b},
 and  Lemma~\ref{le:lln}.
The predictable measure of jumps of $(M^n,L^n)$ is given by
\begin{equation*}
  \tilde{\nu}^n([0,t],\Gamma\times\Gamma')=
\sum_{k=0}^{\lfloor nt\rfloor-1}
\tilde{F}^n\Bl(1-\frac{Q^n_{k}}{n}-\frac{k}{n},\,\Gamma\setminus\{0\}
\Br)\tilde{F}^n\Bl(\Bl(\frac{Q^n_{k}}{n}-\frac{1}{n}\Br)^+,
\,\Gamma'\setminus\{0\}\Br),\,\Gamma,\Gamma'\in\mathcal{B}(\R),
\end{equation*}
where 
\begin{equation*}
  \tilde{F}^n(s,\,\Gamma'')=
\mathbf{P}\Bl(\frac{1}{\sqrt{n}}\sum_{j=1}^{\lfloor ns\rfloor}
\bl(\xi^n_{1j}-\frac{c_n}{n}\br)\in \Gamma''\Br),\,s\in[0,1],\,
\Gamma''\in\mathcal{B}(\R).
\end{equation*}
Therefore, for $\epsilon>0$ and $n$ large enough
\begin{multline*}
  \int_0^1\int_{\R^2}
  \abs{x}^2\,\ind(\abs{x}>\epsilon)\,\tilde{\nu}^n(ds,dx)
\le \frac{1}{\epsilon^2}\int_0^1\int_{\R^2}
  \abs{x}^4\,\,\tilde{\nu}^n(ds,dx)\le 
\frac{2}{\epsilon^2}
\sum_{k=1}^{n}\int_{\R} \abs{x}^4\,
\tilde{F}^n\Bl(1-\frac{Q^n_{k-1}}{n}-\frac{k-1}{n},\,dx\Br)\\+
\frac{2}{\epsilon^2}
\sum_{k=1}^{n}\int_{\R} \abs{x}^4\,
\tilde{F}^n\Bl(\Bl(\frac{Q^n_{k-1}}{n}-\frac{1}{n}\Br)^+,\,dx\Br)
\le\frac{4(2c_n+3c_n^2)}{n^2\epsilon^2},
\end{multline*}
which converges to $0$ as $n\to\infty$. Therefore,
extending the $(M^n,L^n)$ to processes with trajectories in $\D(\R_+,\R^2)$ by
setting $(M^n_t,L^n_t)=(M^n_1,L^n_1),\,t\ge1$, we  see by Jacod and
Shiryaev \citeyear[Theorem VIII.3.22]{jacshir} that these processes converge in
distribution to the extension of $(M,L)$ defined as
$(M_t,L_t)=(M_1,L_1),\,t\ge1$.
Since the projection $p_1$ from $\D(\R_+,\R^2)$ to $\D_C([0,1],\R^2)$ is
continuous at continuous functions from $\D(\R_+,\R^2)$,  we conclude that 
the (non-extended) processes $(M^n,L^n)$ converge in distribution
in $\D_C([0,1],\R^2)$ to the process $(M,L)$.

By \eqref{eq:72}, \eqref{eq:74} and Lipshitz continuity of reflection
for $r\in[0,1]$
\begin{equation*}
    \abs{\overline{Q}^n_r-\overline{q}_r}\le
2\sup_{t\in[0,r]}\Bl|\int_0^{\lfloor nt\rfloor/n}
\Bl(c_n\Bl(1-\overline{Q}^n_{s}-\frac{\lfloor
  ns\rfloor}{n}\Br)-1\Br)\,ds
+\overline{\epsilon}^n_{t}+\overline{M}^n_t 
-\int_0^t\bl(c(1-\overline{q}_s-s)-1\br)\,ds
\Br|,
\end{equation*}
so the definitions of $X^n_t$ and $M^n_t$, 
\eqref{eq:42}, \eqref{eq:6.4b}, and \eqref{eq:6.9a}  yield 
  \begin{equation}
    \label{eq:30}
\abs{X^n_{t}}\le
2  c_n\int_0^t\abs{X_s^n}\,ds
+2\sup_{s\in[0,t]}\abs{M^n_s}+2\sqrt{n}\abs{c_n-c}\int_0^t\sigma_s^2\,ds+
2\sup_{s\in[0,1]}\abs{\tilde{\epsilon}^n_s},\,t\in[0,1].
\end{equation}
In view of  $\C$-tightness  of the $M^n$, the convergence 
$\sqrt{n}(c_n-c)\to\theta$, \eqref{eq:tiep},
  and Gronwall's inequality,  \eqref{eq:30}
 yields the asymptotic boundedness in
 probability of the 
$\sup_{t\in[0,1]}\abs{X^n_t}$  asserted 
 in the first display of the statement  of the lemma.
This implies by \eqref{eq:80}, the convergence
$\sqrt{n}(c_n-c)\to\theta$, \eqref{eq:tide} and $\C$-tightness of the
$L^n$ that the sequence $Z^n,\,n\in\N$, is $\C$-tight in $\D([0,1],\R)$.

We next show that for arbitrary $\delta\in(0,1-\beta)$
\begin{equation}
  \label{eq:82}
  \lim_{n\to\infty}
\mathbf{P}\bl(\sup_{t\in[\beta+\delta,1]}\abs{X^n_t}>\delta\br)=0.
\end{equation}
On recalling the definition of $Y^n_t$, 
we write \eqref{eq:77a} in the following form
\begin{equation}
  \label{eq:77}
  X^n_{t}=-c_n\int_0^tX^n_s\,ds+
\sqrt{n}(c_n-c)\int_0^t\sigma_s^2\,ds
+M^n_t+\tilde{\epsilon}^n_t-\sqrt{n}\overline{\phi}_t
+\sqrt{n}\: \overline{\Phi}^n_t.
\end{equation}
Since $X^n_t=\sqrt{n}\overline{Q}^n_t$ for $t\in[\beta,1]$,
$\overline{\phi}_\beta=0$,  and
$\overline{\Phi}^n_t$ increases only when $\overline{Q}^n_t=0$,
\eqref{eq:77} implies   that $(X^n_t,\,t\in[\beta,1])$
 is the reflection of the process 
$\bl(X^n_\beta-c_n\int_\beta^tX^n_s\,ds+
\sqrt{n}(c_n-c)\int_\beta^t\sigma_s^2\,ds +
(M^n_t-M^n_\beta)-\sqrt{n}\overline{\phi}_t
+(\tilde{\epsilon}^n_t-\tilde{\epsilon}^n_\beta),\,
t\in[\beta,1]\br)$, so by $X^n_s$ being non-negative on $[\beta,1]$ it is not
greater than the reflection of
$\bl(X^n_\beta+\sqrt{n}(c_n-c)\int_\beta^t\sigma_s^2\,ds+(M^n_t-M^n_\beta)
-\sqrt{n}\overline{\phi}_t+(\tilde{\epsilon}^n_t-\tilde{\epsilon}^n_\beta),
\,t\in[\beta,1]\br)$.
Therefore, 
\begin{multline}
  \label{eq:87}
  X^n_t\le 
\sup_{s\in[\beta,t]}\bl(\sqrt{n}(c_n-c)\int_s^t\sigma_p^2\,dp+(M^n_t-M^n_s)
+\sqrt{n}(\overline{\phi}_s-\overline{\phi}_t)+(\tilde{\epsilon}^n_t
-\tilde{\epsilon}^n_s)\br)\\ \vee
\br(X^n_\beta+\sqrt{n}(c_n-c)\int_\beta^t\sigma_s^2\,ds+(M^n_t-M^n_\beta)
-\sqrt{n}\overline{\phi}_t+(\tilde{\epsilon}^n_t-\tilde{\epsilon}^n_\beta)\br),
\;t\in[\beta,1].
\end{multline}
Hence, for  $t\ge\beta+\delta$ and $\eta\in(0,\delta)$,
\begin{multline}
  \label{eq:89}
  X^n_t  \le\bl(\abs{\sqrt{n}(c_n-c)}
\int_\beta^1\sigma_s^2\,ds
+2 \sup_{s\in[\beta,1]}\abs{M^n_s}
+2 \sup_{s\in[\beta,1]}\abs{\tilde{\epsilon}^n_s}+X^n_\beta
+\sqrt{n}(\overline{\phi}_{t-\eta}
-\overline{\phi}_t)\br)\\
\vee \sup_{s\in[t-\eta,t]}(\abs{\sqrt{n}(c_n-c)}\int_s^t\sigma_p^2\,dp
+\abs{M^n_t-M^n_s}+\abs{\tilde{\epsilon}^n_t-\tilde{\epsilon}^n_s}).
\end{multline}
Limit \eqref{eq:82} follows by \eqref{eq:89}, \eqref{eq:tiep},  
$\C$-tightness of  the $M^n$, asymptotic boundedness in probability of
the $\sup_{t\in[0,1]}\abs{X^n_t}$,
 the convergence $\sqrt{n}(c_n-c)\to \theta$, and 
convergence of $\sup_{t\in[\beta+\delta,1]}
\sqrt{n}(\overline{\phi}_{t-\eta}-\overline{\phi}_t)$ to $-\infty$
as $n\to\infty$.
Now,  \eqref{eq:82} implies by
 \eqref{eq:77a},   \eqref{eq:tiep},  the convergence 
$\sqrt{n}(c_n-c)\to\theta$, asymptotic boundedness in probability of
the $\sup_{t\in[0,1]}\abs{X^n_t}$, and $\C$-tightness of the $M^n$   that
the processes $Y^n$ restricted to $[\beta+\delta,1]$ are
$\C$-tight in $\D([\beta+\delta,1],\R)$.

Let us now assume that $\beta>0$. By \eqref{eq:77}, the definition of $X^n_t$,
 and the definition of the reflection mapping  for $t\in[0,1]$
\begin{equation}
  \label{eq:79}
  \sqrt{n}\: \overline{\Phi}^n_t=-\inf_{s\in[0,t]}
\Bl(-c_n\int_0^sX^n_p\,dp+\sqrt{n}(c_n-c)\int_0^s\sigma_p^2\,dp+M^n_s+
\tilde{\epsilon}^n_s+\sqrt{n}\overline{q}_s-\sqrt{n}\overline{\phi}_s
\Br)\wedge0.
\end{equation}
Convergence in distribution of the $M^n$ to a continuous-path process
implies that for $\delta>0$
$\lim_{\eta\to0}  \limsup_{n\to\infty}\mathbf{P}\bl(\sup_{t\in[0,\eta]}\abs{M^n_t}
>\delta\br)=0$.
Therefore,  given $\delta\in(0,\beta)$, we derive from \eqref{eq:79},
  taking into consideration the convergences
$\sqrt{n}(c_n-c)\to\theta$ and
$\sqrt{n}\inf_{t\in[\eta,\beta-\delta]}\overline{q}_t\to\infty$  as
$n\to\infty$, where  $\eta\in(0,\beta-\delta)$, the fact that 
$\overline{\phi}_t=0$ for $t\in[0,\beta]$, \eqref{eq:tiep},
 and asymptotic boundedness in probability of the 
$\sup_{t\in[0,1]}\abs{X^n_t}$ and $\sup_{t\in[0,1]}\abs{M^n_t}$ that
\begin{equation}
  \label{eq:78}
    \lim_{n\to\infty}\mathbf{P}\bl(\sup_{t\in[0,\beta-\delta]}
\abs{Y^n_t}>\delta \br)=0.
\end{equation}
Putting together  \eqref{eq:77a}, \eqref{eq:tiep},   \eqref{eq:78}, the
convergence $\sqrt{n}(c_n-c)\to\theta$, asymptotic boundedness 
in probability of the $\sup_{t\in[0,1]}\abs{X^n_t}$, and
$\C$-tightness of the $M^n$, we conclude 
  that the  $X^n$ restricted to $[0,\beta-\delta]$ are $\C$-tight in 
$\D_C([0,\beta-\delta],\R)$.

We have thus established that for $\beta>0$ and
$\delta\in(0,\beta\wedge(1-\beta))$  the processes 
$M^n$, $L^n$, $X^n$ restricted to $[0,\beta-\delta]$, 
$Y^n$ restricted to $[\beta+\delta,1]$, and $Z^n$ 
 are $\C$-tight in the
associated function spaces, so they are jointly tight as random
elements with values in the product space.  Convergence in
distribution in $\D_C([0,1],\R^2)\times\D_C([0,\beta-\delta],\R)\times
\D_C([\beta+\delta,1],\R)\times\D_C([0,1],\R)$ of   the 
$\bl(M^n,L^n,(X^n_t,\,t\in[0,\beta-\delta]),
(Y^n_t,\,t\in[\beta+\delta,1]),Z^n\br)$
 to
$\bl(M,L,(H_t,\,t\in[0,\beta-\delta]),
(-H_t,\,t\in[\beta+\delta,1]),Z\br)$ now follows
by  \eqref{eq:77a}, \eqref{eq:80}, \eqref{eq:tiep},  \eqref{eq:tide},
 \eqref{eq:6.11a}, \eqref{eq:90}, 
\eqref{eq:82}, \eqref{eq:78}, the convergence
$\sqrt{n}(c_n-c)\to\theta$, 
 convergence in distribution of the $(M^n,L^n)$ to $(M,L)$,
and uniqueness of the solution $(H,Z)$ to  \eqref{eq:6.11a} and \eqref{eq:90}.

Let us now assume that $\beta=0$. 
Inequality \eqref{eq:30} in view of asymptotic boundedness in
probability of the $\sup_{t\in[0,1]}\abs{X^n_t}$, $\C$-tightness of
the $M^n$, limits \eqref{eq:tiep}, \eqref{eq:82}, and 
$\sqrt{n}(c_n-c)\to\theta$ yields the limit
$\lim_{\eta\to0}\limsup_{n\to\infty}
\mathbf{P}\bl(\sup_{t\in[0,\eta]}\abs{X^n_t}>\delta\br)=0$
for $\delta>0$, so by  \eqref{eq:82}
$\lim_{n\to\infty}
\mathbf{P}\bl(\sup_{t\in[0,1]}\abs{X^n_t}>\delta\br)=0$.
Therefore, by  \eqref{eq:77a}, the convergence
$\sqrt{n}(c_n-c)\to\theta$,  and convergence in
distribution of the $M^n$ to $M$,  the  $Y^n$ 
 converge in
distribution  in $\D_C([0,1],\R)$ to $-H$.
\end{proof}
\begin{remark}
A slight modification of the proof
allows one to strengthen the assertion of the  lemma for $\beta>0$ 
to the joint
  convergence in distribution in
$\D_C([0,1],\R^2)\times\D_C([0,\beta-\delta],\R)^2\times
\D_C([\beta+\delta,1],\R)^2\times\R^2\times\D_C([0,1],\R)$ of
the $M^n$, $L^n$,
 $(X^n_t,\,t\in[0,\beta-\delta])$, $(Y^n_t,\,t\in[0,\beta-\delta])$,
$(X^n_t,\,t\in[\beta+\delta,1])$,  
$(Y^n_t,\,t\in[\beta+\delta,1])$, $X^n_\beta$,  $Y^n_\beta$, and $Z^n$
  to the respective random elements
  $M$, $L$, $(X_t,\,t\in[0,\beta-\delta])$, $(Y_t,\,t\in[0,\beta-\delta])$,
$(X_t,\,t\in[\beta+\delta,1])$, 
$(Y_t,\,t\in[\beta+\delta,1])$, $X_\beta$,  $Y_\beta$, and $Z$, where 
  \begin{align*}
  X_t&=
\left\{  \begin{aligned}
\displaystyle
   & H_t\,&&\text{ for
    }t\in[0,\beta),\\
&H_\beta\vee0\,&&\text{ for }t=\beta,\\
&0\,&&\text{ for }t\in(\beta,1],
  \end{aligned}\right.
&\text{ and }&&Y_t=\left\{
  \begin{aligned}
&0\,&&\text{ for }t\in[0,\beta),\\
&(-H_\beta)\vee0\,&&\text{ for }t=\beta,\\
\displaystyle&-H_t\,&&\text{ for
    }t\in(\beta,1].
  \end{aligned}\right.
\end{align*} 
We thus have convergence in
  distribution with unmatched
  jumps in the limit process mentioned in the introduction.  
\end{remark}
\begin{proof}[Proof of Theorem~\ref{the:7}]
Let $c>1$, so $\beta>0$.
  We prove that as $n\to\infty$
  \begin{equation}
    \label{eq:84}
  \Bl(\sqrt{n}\Bl(\frac{\alpha^n}{n}-\alpha\Br),\,\sqrt{n}\Bl(\frac{\beta^n}{n}
-\beta\Br),\,\sqrt{n}\Bl(\frac{\gamma^n}{n}-\gamma\Br)\Br) \overset{d}{\to}
\Bl(-H_1,\,\frac{H_\beta}{1-c(1-\beta)},\,Z_\beta\Br),
  \end{equation}
which implies the assertion of part 2 of the theorem.

Let $\tau^n$ be the last time $t$ before $\beta/2$ when 
$\overline{Q}^n_t=0$ and 
$\tilde{\beta}^n$ be the first time $t$ not before  $\beta/2$ when
$\overline{Q}^n_t=0$. By 
Lemma~\ref{le:lln} and \eqref{eq:42}
$\overline{Q}^n_t>0$ for $t\in[\delta,\beta-\delta]$
with probability tending to $1$ as $n\to\infty$ for arbitrary
$\delta\in(0,\beta/2)$, so $\mathbf{P}(\tau^n\le\delta)\to1$ 
and $\mathbf{P}(\tilde{\beta}^n\ge
\beta-\delta)\to1$. Also, noting that $\overline{\Phi}^n_t=
\overline{\Phi}^n_{\tau^n}$ for $t\in(\tau^n,\tilde{\beta}^n)$,
  Lemma~\ref{le:lln}, and  \eqref{eq:31}
\begin{equation*}
\limsup_{n\to\infty}\mathbf{P}(\tilde{\beta}^n>\beta+\delta)
\le\limsup_{n\to\infty}\mathbf{P}\br(\overline{\Phi}^n_{\tau^n}
=\overline{\Phi}^n_{\beta+\delta}\bl)\le
\ind(0=\overline{\phi}_{\beta+\delta})=0\,,
\end{equation*}
so as $n\to\infty$
\begin{equation}
  \label{eq:50}
\tilde{\beta}^n\overset{\mathbf{P}}{\to}\beta.
\end{equation}
Similarly, the event that  there exists an excursion of
 $\overline{Q}^n$ of duration greater than $\eta$, where
 $\eta\in(0,1-\beta)$,
which  ends at some time after $\beta+\eta$, 
 is contained in the event
$\{\inf_{t\in[\beta,1-\eta]}(\overline{\Phi}^n_{t+\eta}-
\overline{\Phi}^n_t)=0\}$. 
Lemma~\ref{le:lln} and the fact that
$\overline{\phi}_t$ is strictly increasing on $[\beta,1]$ in view of 
\eqref{eq:31}  imply  that the
probability of the latter event tends to $0$ as $n\to\infty$.
As the sizes of the connected components of $\mathcal{G}(n,c_n/n)$ are 
equal to $n$ multiplied by the excursion lengths of $\overline{Q}^n$,
we see that with probability tending to $1$ as $n\to\infty$ the
largest component ``starts'' at $n\tau^n$ and ``ends'' at
$n\tilde{\beta}^n$, so 
\begin{align}
  \label{eq:51}
\mathbf{P}\Bl(\frac{\beta^n}{n}=\tilde{\beta}^n-\tau^n  \Br)&\to1,\\
  \label{eq:93}
  \mathbf{P}\Bl(\frac{\gamma^n}{n}=\overline{E}^n_{\tilde{\beta}^n}-
\overline{E}^n_{\tau^n } \Br)&\to1\,.
\end{align}
By \eqref{eq:77a} and the facts that
 $X^n_{\tau^n}=-\sqrt{n}\overline{q}_{\tau^n}$ and 
 $X^n_{\tilde{\beta}^n}=-\sqrt{n}\overline{q}_{\tilde{\beta}^n}$,
    \begin{align}
     \label{eq:39}
-\sqrt{n}\overline{q}_{\tau^n}
& =-c_n\int_{0}^{\tau^n}X^n_s\,ds+
\sqrt{n}(c_n-c)\int_0^{\tau^n}\sigma_s^2\,ds
+M^n_{\tau^n}+\tilde{\epsilon}^n_{\tau^n}+Y^n_{\tau^n},\\
   \label{eq:38}
-\sqrt{n}\overline{q}_{\tilde{\beta}^n}
&=-c_n\int_{0}^{\tilde{\beta}^n}X^n_s\,ds+
\sqrt{n}(c_n-c)\int_0^{\tilde{\beta}^n}\sigma_s^2\,ds
+M^n_{\tilde{\beta}^n}+
\tilde{\epsilon}^n_{\tilde{\beta}^n}
+\sqrt{n}\overline{\Phi}^n_{\tilde{\beta}^n}
-\sqrt{n}\overline{\phi}_{\tilde{\beta}^n}.
\end{align}
Since $\tau^n\overset{\mathbf{P}}{\to}0$, the right-hand side 
of \eqref{eq:39} converges in probability to zero by \eqref{eq:tiep}
and Lemma~\ref{the:diff}, so
$\sqrt{n}\overline{q}_{\tau^n}\overset{\mathbf{P}}{\to}0$ and, consequently, by
 \eqref{eq:32} and the fact that $c>1$
\begin{equation}
  \label{eq:49}
  \sqrt{n}\tau^n\overset{\mathbf{P}}{\to}0.
\end{equation}
Since $\overline{\Phi}^n_{\tilde{\beta}^n}=\overline{\Phi}^n_{\tau^n}
+1/n$ (see \eqref{eq:36}),
$\sqrt{n}\overline{\Phi}^n_{\tau^n}\overset{\mathbf{P}}{\to}0$,
and $\overline{q}_{\tilde{\beta}^n}-\overline{\phi}_{\tilde{\beta}^n}=
\int_0^{\tilde{\beta}^n}\bl(c(1-\overline{q}_s-s)-1\br)\,ds=
\int_\beta^{\tilde{\beta}^n}\bl(c(1-\overline{q}_s-s)-1\br)\,ds$ (see
\eqref{eq:74}),  we derive
from \eqref{eq:38} on  using \eqref{eq:50}, \eqref{eq:tiep}, and 
 Lemma~\ref{the:diff} that 
\begin{equation}
  \label{eq:86}
\sqrt{n}\int_\beta^{\tilde{\beta}^n}\bl(c(1-\overline{q}_s-s)-1\br)\,ds
  -c_n\int_{0}^{\tilde{\beta}^n}X^n_s\,ds
+\sqrt{n}(c_n-c)\int_0^{\tilde{\beta}^n}\sigma_s^2\,ds  
+M^n_{\tilde{\beta}^n}
\overset{\mathbf{P}}{\to}0.
\end{equation}
Since  $\alpha=\overline{\phi}_1$ (see \eqref{eq:31}) and 
$\alpha^n=\Phi^n_n$, we also have that
\begin{equation}
  \label{eq:91}
\sqrt{n}\Bl(\frac{\alpha^n}{n}-\alpha\Br)=Y^n_1.
\end{equation}
Convergence \eqref{eq:84} follows by \eqref{eq:50},  
 \eqref{eq:51}, \eqref{eq:93}, \eqref{eq:49}, \eqref{eq:86}, \eqref{eq:91},
the observation that $\gamma=\overline{e}_\beta$ (see \eqref{eq:76}),
asymptotic boundedness in probability  of the
$\sup_{t\in[0,1]}\abs{X^n_t}$, the convergence
 $\sqrt{n}(c_n-c)\to\theta$, the joint convergence in distribution
$\bl(M^n,Y^n_1,(X^n_s,\,s\in[0,\beta-\delta]),Z^n\br)
\overset{d}{\to}
\bl(M,-H_1,(H_s,\,s\in[0,\beta-\delta]),Z\br)$ 
in  $\D_C([0,1],\R)\times\R\times\D_C([0,\beta-\delta],\R)\times\D_C([0,1],\R)$
valid by Lemma~\ref{the:diff}, 
and the continuous mapping theorem.

If $c\le 1$ the
$Y^n_1$ converge in distribution to $-H_1$ by part 2 of
Lemma~\ref{the:diff}, which completes the proof of part 1.
\end{proof}
We now prove Theorem~\ref{the:6}.
As mentioned above, the proof is along the lines of
 the proof of Theorem~\ref{the:7}, so we begin with an
idempotent analogue of Lemma~\ref{the:diff}. 
We recall that $b_n,\, n\in\N,$ is a real-valued sequence such that
$b_n\to\infty$ and $b_n/\sqrt{n}\to0$ as $n\to\infty$, and 
introduce processes  $\hat{M}^n=(\hat{M}^n_t,\,t\in[0,1])$, 
$\hat{L}^n=(\hat{L}^n_t,\,t\in[0,1])$, 
$\hat{X}^n=(\hat{X}^n_t,\,t\in[0,1])$,
 $\hat{Y}^n=(\hat{Y}^n_t,\,t\in[0,1])$, 
 and $\hat{Z}^n=(\hat{Z}^n_t,\,t\in[0,1])$
 by the respective equalities $\hat{M}^n_t=M^n_t/b_n$,
$\hat{L}^n_t=L^n_t/b_n$,
$  \hat{X}^n_t=X^n_t/b_n$,
$  \hat{Y}^n_t=Y^n_t/b_n$, and 
$  \hat{Z}^n_t=Z^n_t/b_n$.
Dividing \eqref{eq:77a} and \eqref{eq:80} through by $b_n$ yields
for $t\in[0,1]$
\begin{align}
  \label{eq:77h}
  \hat{X}^n_{t}&=-c_n\int_0^t\hat{X}^n_s\,ds+
\frac{\sqrt{n}}{b_n}\,(c_n-c)\int_0^t\sigma_s^2\,ds
+\hat{M}^n_t+\hat{\epsilon}^n_t+\hat{Y}^n_t\,,\\
  \label{eq:80h}
\hat{Z}^n_t&=c_n\int_0^t\hat{X}^n_s\,ds
+\frac{\sqrt{n}}{b_n}\,(c_n-c)\int_0^t\overline{q}_s\,ds+
\hat{L}^n_t+\hat{\delta}^n_t,
  \end{align}
where 
\begin{align}
  \label{eq:hep}
\hat{\epsilon}^n_t=\frac{\tilde{\epsilon}^n_t}{b_n},\quad
\hat{\delta}^n_t=\frac{\tilde{\delta}^n_t}{b_n}.
\end{align}
We note that by \eqref{eq:6.9a}, \eqref{eq:81}, 
\eqref{eq:hep},  and Lemma~\ref{le:eps},
\begin{align}
\label{eq:100}
  \sup_{t\in[0,1]}\abs{\hat{\epsilon}^n_t}
&\overset{\mathbf{P}^{1/b_n^2}}{\to}0,\\
\intertext{provided $c_n\to c$ as $n\to\infty$,
and}  \label{eq:100a}
\sup_{t\in[0,1]}\abs{\hat{\delta}^n_t}&\le\frac{c_n}{b_n\sqrt{n}}.
\end{align}
Let
$\hat{W}^{(1)}=(\hat{W}^{(1)}_t,\,t\in[0,1])$ and
$\hat{W}^{(2)}=(\hat{W}^{(2)}_t,\,t\in[0,1])$  be independent 
 idempotent  Wiener processes on an idempotent probability space
 $(\Upsilon,\mathbf{\Pi})$ adapted to a complete $\tau$-flow $\mathbf{A}$,
idempotent processes  $\hat{M}=(\hat{M}_t,\,t\in[0,1])$
and $\hat{L}=(\hat{L}_t,\,t\in[0,1])$ be defined by 
$\hat{M}_t=\sqrt{c} \int_0^t\sigma_s\,\dot{\hat{W}}^{(1)}_s\,ds$
and $\hat{L}_t=\sqrt{c}
 \int_0^t\sqrt{\overline{q}_s}\,\dot{\hat{W}}_s^{(2)}\,ds$, 
respectively, an idempotent
 process
$\hat{H}=(\hat{H}_t,\,t\in[0,1])$ be the Luzin strong  solution of the
equation 
\begin{equation}
  \label{eq:113}
    \hat{H}_t=
-c\int_0^{t\wedge \beta}\hat{H}_s\,ds
+\hat{\theta}\int_0^t\sigma^2_s\,ds+
\sqrt{c}\int_0^t\sigma_s\dot{\hat{W}}^{(1)}_s\,  ds,
\end{equation}
and an idempotent process $\hat{Z}=(\hat{Z}_t,\,t\in[0,1])$ be given by
\begin{equation}
  \hat{Z}_t=c\int_0^{t\wedge\beta}\hat{H}_s\,ds
+\hat{\theta}\int_0^t\overline{q}_s\,ds+
\sqrt{c}\int_0^t\sqrt{\overline{q}_s}\,\dot{\hat{W}}^{(2)}_s\,ds.
  \label{eq:90h}
\end{equation}
 \begin{lemma}
  \label{the:idemdiff} 
Let $(\sqrt{n}/b_n)(c_n-c)\to\hat{\theta}\in\R$ as $n\to\infty$, where
$c>0$, $b_n\to\infty$ and $b_n/\sqrt{n}\to0$.
Then for arbitrary $\eta>0$
\begin{align*}
\lim_{n\to\infty}\mathbf{P}(\sup_{t\in[0,1]}
\abs{\overline{Q}^n_t-\overline{q}_t}>\eta)^{1/b_n^2}&=0,\\
\lim_{n\to\infty}\mathbf{P}(\sup_{t\in[0,1]}
\abs{\overline{\Phi}^n_t-\overline{\phi}_t}>\eta)^{1/b_n^2}&=0,\\  
\intertext{and}    \lim_{B\to\infty}\limsup_{n\to\infty}
\mathbf{P}(\sup_{t\in[0,1]}\abs{\hat{X}^n_t}>B)^{1/b_n^2}&=0.
\end{align*}
Also the following holds.
\begin{enumerate}
\item 
If $\beta>0$, then 
for   $\delta\in(0,\beta\wedge(1-\beta))$ the stochastic processes
 $\hat{M}^n$, $\hat{L}^n$, $(\hat{X}^n_t,\,t\in[0,\beta-\delta])$, 
$(\hat{Y}^n_t,\,t\in[\beta+\delta,1])$, and 
 $(\hat{Z}^n_t,\,t\in[0,1])$
 jointly LD converge in distribution at rate $b_n^2$ in
$\D_C([0,1],\R^2)\times\D_C([0,\beta-\delta],\R)\times
\D_C([\beta+\delta,1],\R)\times\D_C([0,1],\R)$
 to the respective idempotent processes
  $\hat{M}$, $\hat{L}$, $(\hat{H}_t,\,t\in[0,\beta-\delta])$, 
$(-\hat{H}_t,\,t\in[\beta+\delta,1])$, 
 and $\hat{Z}$.
In addition,
$    \lim_{n\to\infty}\mathbf{P}(\sup_{t\in[0,\beta-\delta]}
\abs{\hat{Y}^n_t}>\delta)^{1/b_n^2}=0.
$
  \item
If $\beta=0$, then 
 the stochastic processes
  $\hat{Y}^n$
 LD converge in distribution at rate $b_n^2$ in
$\D_C([0,1],\R)$
 to the  idempotent process
   $-\hat{H}$.
\end{enumerate}
\end{lemma}
\begin{proof}
We have by \eqref{eq:69} and \eqref{eq:29}
\begin{align*}
    \hat{M}^n_t&=\frac{1}{b_n \sqrt{n}}
\sum_{i=1}^{\lfloor nt\rfloor}
\sum_{j=1}^{n-Q^n_{i-1}-(i-1)}\bl(\xi^n_{ij}-\frac{c_n}{n}\br),\;
t\in[0,1],\\
    \hat{L}^n_t&=\frac{1}{b_n\sqrt{n}}
\sum_{i=1}^{\lfloor nt\rfloor}
\sum_{j=1}^{Q^n_{i-1}-1}\bl(\zeta^n_{ij}-\frac{c_n}{n}\br),\;
t\in[0,1].
\end{align*}
Therefore, the $\mathbf{F}^n$-predictable measure of jumps
of $(\hat{M}^n,\hat{L}^n)$ 
has the form
\begin{equation}
  \label{eq:95}
    \hat{\nu}^n([0,t],\Gamma\times\Gamma')=\sum_{k=0}^{\lfloor
    nt\rfloor-1}\hat{F}^{n}\Bl(1-\overline{Q}^n_{k/n}-
\frac{k}{n},\Gamma\setminus\{0\}\Br)
\hat{F}^n\Bl( \Bl(\overline{Q}^n_{k/n}-\frac{1}{n}\Br)^+,\,
\Gamma'\setminus\{0\}\Br),\;\Gamma,\Gamma'\in\mathcal{B}(\R),
\end{equation}
where
  \begin{equation}
    \label{eq:85}
  \hat{F}^n(s,\Gamma'')=\mathbf{P}\Bl(\frac{1}{b_n\sqrt{n}}\sum_{j=1}^{\lfloor
  ns\rfloor}\Bl(\xi^n_{1j}-\frac{c_n}{n}\Br)
\in \Gamma''\Br),\,s\in[0,1],\,\Gamma''\in\mathcal{B}(\R).
\end{equation} 
Accordingly, the stochastic exponential
$(\hat{\mathcal{E}}^n_t(\lambda),\,t\in[0,1])$, 
where $\lambda\in\R$,
associated with $\hat{M}^n$ is given by
\begin{multline*}
\log  \hat{\mathcal{E}}^n_t(\lambda)=
\sum_{k=1}^{\lfloor nt\rfloor }
\log\Bl(1+\int_{\R}\bl(e^{\lambda
  x}-1\br)\,\hat{\nu}^n\bl(\bl\{\frac{k}{n}\br\},dx\times\R\br)\Br)\\
=n\log\Bl(\mathbf{E}\exp\Bl(\frac{\lambda}{b_n\sqrt{n}}
\bl(\xi^n_{11}-\frac{c_n}{n}\br)\Br)
\Br)\sum_{k=0}^{\lfloor nt\rfloor-1 }
 \Bl(1-\overline{Q}^n_{k/n}-\frac{k}{n}\Br).
\end{multline*}
Since, for $B>0$, by Doob's inequality
\begin{multline*}
  \mathbf{P}(\sup_{t\in[0,1]}\abs{\hat{M}^n_t}>B)^{1/b_n^2}\le
e^{-B}\Bl(\bl(\mathbf{E}\,e^{b_n^2\hat{M}^n_1}\br)^{1/b_n^2}
+\bl(\mathbf{E}\,e^{-b_n^2\hat{M}^n_1}\br)^{1/b_n^2}\Br)\\
\le
e^{-B}\Bl(\bl(\mathbf{E}\,\hat{\mathcal{E}}^n_1(2b_n^2)\br)^{1/(2b_n^2)}
+\bl(\mathbf{E}\,\hat{\mathcal{E}}^n_1(-2b_n^2)\br)^{1/(2b_n^2)}\Br)
\end{multline*}
and $(n/b_n)^2\log \mathbf{E}\exp\bl(\pm(2b_n/\sqrt{n})%
\bl(\xi^n_{11}-c_n/n\br)\br)\to 2 c$ as $n\to\infty$, we conclude that
\begin{equation}
  \label{eq:92}
\lim_{B\to\infty}
\limsup_{n\to\infty}
\mathbf{P}(\sup_{t\in[0,1]}\abs{\hat{M}^n_t}>B)^{1/b_n^2}=0.
\end{equation}
Dividing  \eqref{eq:30} through by $b_n$ and recalling \eqref{eq:hep} yields
\begin{equation}
  \label{eq:45}
  \abs{\hat{X}^n_{t}}\le
2  c_n\int_0^t\abs{\hat{X}_s^n}\,ds+
2\frac{\sqrt{n}}{b_n}\abs{c_n-c}\int_0^t\sigma_s^2\,ds
+2\sup_{s\in[0,t]}\abs{\hat{M}^n_s}+
2\sup_{s\in[0,1]}\abs{\hat{\epsilon}^n_s},\,t\in[0,1].
\end{equation}
Applying Gronwall's inequality to \eqref{eq:45}, we have 
by  \eqref{eq:100}, \eqref{eq:92}, and the
convergence $(\sqrt{n}/b_n)(c_n-c)\to\hat{\theta}$ that
\begin{equation}
  \label{eq:88}
  \lim_{B\to\infty}
\limsup_{n\to\infty}\mathbf{P}(\sup_{t\in[0,1]}\abs{\hat{X}^n_t}>B)^{1/b_n^2}=0
\end{equation}
proving the third display   in the statement of the lemma.
As a consequence of \eqref{eq:88}, the definition of 
$\hat{X}^n_t$, and the convergence $\sqrt{n}/b_n\to\infty$
\begin{equation}
  \label{eq:107}
\lim_{n\to\infty}
\mathbf{P}(\sup_{t\in[0,1]}\abs{\overline{Q}^n_t-\overline{q}_t}>
\eta)^{1/b_n^2}=0\,,
\end{equation}
and then by \eqref{eq:72}, \eqref{eq:74}, \eqref{eq:100}, and
\eqref{eq:92}
\begin{equation*}
\lim_{n\to\infty}
\mathbf{P}(\sup_{t\in[0,1]}\abs{\overline{\Phi}^n_t-\overline{\phi}_t}>
\eta)^{1/b_n^2}=0
\end{equation*}
 for arbitrary $\eta>0$, proving the other claimed
 super-exponential convergences in probability.

We now prove that  the $(\hat{M}^n,\hat{L}^n)$ LD
converge in distribution at rate $b_n^2$ 
to $(\hat{M},\hat{L})$ in $\D_C([0,1],\R^2)$. 
This is accomplished by
checking the conditions of Corollary~4.3.13 in Puhalskii \citeyear{Puh01}.
Extending $\hat{M}^n$ and $\hat{L}^n$ to  processes defined on $\R_+$ by
letting $\hat{M}^n_t=\hat{M}^n_1$ and $\hat{L}^n_t=\hat{L}^n_1$ for 
$t\ge 1$, 
we have by \eqref{eq:73} and \eqref{eq:70}
that  $\hat{M}^n$ and $\hat{L}^n$ are orthogonal
$\mathbf{F}^n$-square integrable martingales with respective
 $\mathbf{F}^n$-predictable quadratic characteristics
\begin{align*}
  \langle\hat{M}^n\rangle_t&=\frac{c_n}{b_n^2}\,\Bl(1-\frac{c_n}{n}\Br) 
\int_0^{\lfloor n(t\wedge1)\rfloor/n}
\Bl(1-\overline{Q}^n_{s}-\frac{\lfloor ns\rfloor}{n}\Br)\,ds,\\
    \langle\hat{L}^n\rangle_t&=\frac{c_n}{b_n^2}\Bl(1-\frac{c_n}{n}\Br) 
\int_{0}^{\lfloor nt\rfloor/n}
\Bl(\overline{Q}^n_{s}-\frac{1}{n}\Br)^+\,ds,
\end{align*}
so by \eqref{eq:42}, \eqref{eq:6.4b}, and \eqref{eq:107} for $\epsilon>0$
\begin{align*}
  \lim_{n\to\infty}\mathbf{P}\Bl(\abs{b_n^2\langle\hat{M}^n\rangle_t-
c\int_0^{t\wedge1}\sigma^2_s\,ds}>\epsilon\Br)^{1/b_n^2}&=0,\\
  \lim_{n\to\infty}\mathbf{P}\Bl(\abs{b_n^2\langle\hat{L}^n\rangle_t-
c\int_0^{t\wedge1}\overline{q}_s\,ds}>\epsilon\Br)^{1/b_n^2}&=0,
\end{align*}
checking condition $(C_0')$ of the corollary. The processes
$(\hat{M}^n,\hat{L}^n)$ 
satisfy the Cram\'er condition by  \eqref{eq:95} and \eqref{eq:85}.
We check condition $(L_e)$:
\begin{equation}
    \label{eq:98}
      \lim_{n\to\infty}\mathbf{P}\Bl(\frac{1}{b_n^2}\int_0^1\int_{\R^2}
e^{\lambda b_n^2\abs{x}}
\ind(b_n^2\abs{x}>\epsilon)\hat{\nu}^n(ds,dx)>\eta\Br)^{1/b_n^2}=0,
\;\lambda>0,\,\epsilon>0,\,\eta>0.
\end{equation}
We have for $n$ large enough by \eqref{eq:95} and \eqref{eq:85}
\begin{multline*}
\frac{1}{b_n^2}\int_0^1\int_{\R^2} 
e^{\lambda b_n^2\abs{x}}
\ind(b_n^2\abs{x}>\epsilon)\hat{\nu}^n(ds,dx)
\le \frac{e^{-\epsilon\sqrt{n}/b_n}}{b_n^2}
\int_0^1\int_{\R^2} e^{(\lambda+\epsilon) b_n\sqrt{n}\abs{x}}
\,\hat{\nu}^n(ds,dx)
\\\le\frac{e^{-\epsilon\sqrt{n}/b_n}}{2b_n^2}
\sum_{k=0}^{n-1}\biggl(\int_\R e^{2(\lambda+\epsilon)b_n\sqrt{n}\abs{x}}
\hat{F}^{n}\Bl(1-\overline{Q}^n_{k/n}-
\frac{k}{n},\,dx\Br)+
\int_\R e^{2(\lambda+\epsilon)b_n\sqrt{n}\abs{x}}
\hat{F}^n\Bl(\Bl(\overline{Q}^n_{k/n}-\frac{1}{n}\Br)^+,\,dx
\Br)\biggr)\\
\le e^{-\epsilon\sqrt{n}/b_n}\frac{n}{b_n^2}\,
e^{c_n\bl(\exp(2(\lambda+\epsilon))-1+2(\lambda+\epsilon)\br)}.
\end{multline*}
 Since the latter expression converges to $0$ as
$n\to\infty$,
convergence \eqref{eq:98} holds. Conditions $(0)$ and $(\sup B')$ of the
corollary trivially hold. Thus, the extended $(\hat{M}^n,\hat{L}^n)$
LD converge in distribution in $\D(\R_+,\R^2)$ at rate $b_n^2$ 
 to $(\hat{M},\hat{L})$. Since the projection $p_1$ from $\D(\R_+,\R^2)$ 
to $\D_C([0,1],\R^2)$ is
continuous at continuous functions from $\D(\R_+,\R^2)$,  we conclude by
the contraction principle that 
the processes $(\hat{M}^n,\hat{L}^n)$ 
LD converge in distribution at rate $b_n^2$
in $\D_C([0,1],\R^2)$ to the idempotent process $(\hat{M},\hat{L})$.
As a byproduct of $\C$-exponential 
tightness of the $\hat{L}^n$, we deduce by
 \eqref{eq:88},  \eqref{eq:80h}, the convergence
$(\sqrt{n}/b_n)(c_n-c)\to\hat{\theta}$, and \eqref{eq:100a}  
that the sequence $\hat{Z}^n,\,n\in\N$, is 
$\C$-exponentially tight in $\D([0,1],\R)$.

We next show that for arbitrary $\delta>0$
\begin{equation}
  \label{eq:103}
  \lim_{n\to\infty}
\mathbf{P}%
\bl(\sup_{t\in[\beta+\delta,1]}\abs{\hat{X}^n_t}>\delta\br)^{1/b_n^2}=0.
\end{equation}
Dividing  \eqref{eq:89} through by $b_n$ yields
 for  $t\ge\beta+\delta$ and $\eta\in(0,\delta)$
\begin{multline*}
  \hat{X}^n_t  \le\Bl(\frac{\sqrt{n}}{b_n}\,\abs{c_n-c}
\int_\beta^1\sigma_s^2\,ds
+2 \sup_{s\in[\beta,1]}\abs{\hat{M}^n_s}
+2 \sup_{s\in[\beta,1]}\abs{\hat{\epsilon}^n_s}+\hat{X}^n_\beta
+\frac{\sqrt{n}}{b_n}(\overline{\phi}_{t-\eta}
-\overline{\phi}_t)\Br)\\
\vee \sup_{s\in[t-\eta,t]}(\frac{\sqrt{n}}{b_n}\,\abs{c_n-c}
\int_s^t\sigma_p^2\,dp
+\abs{\hat{M}^n_t-\hat{M}^n_s}
+\abs{\hat{\epsilon}^n_t-\hat{\epsilon}^n_s}).
\end{multline*}
Convergence \eqref{eq:103} follows  if we
recall that the $\hat{M}^n$ are $\C$-exponentially tight of order $b_n^2$, 
$(\sqrt{n}/b_n)(c_n-c)\to \hat{\theta}$,
\eqref{eq:100} and \eqref{eq:88}   hold,
and use that $\sup_{t\in[\beta+\delta,1]}
(\overline{\phi}_{t-\eta}-\overline{\phi}_t)<0$. 
Consequently, 
by  \eqref{eq:77h}, \eqref{eq:100},  \eqref{eq:88}, \eqref{eq:103},
 $\C$-exponential tightness of the $\hat{M}^n$,
 and the convergence $(\sqrt{n}/b_n)(c_n-c)\to \hat{\theta}$  
the processes $\hat{Y}^n$ restricted to $[\beta+\delta,1]$ are
$\C$-exponentially tight of order $b_n^2$.

Next, let us assume that $\beta>0$. Representation \eqref{eq:79}
implies that for $t\in[0,1]$
\begin{equation}
  \label{eq:43}
\frac{\sqrt{n}}{b_n}  \overline{\Phi}^n_t=-\inf_{s\in[0,t]}
\Bl(-c_n\int_0^s\hat{X}^n_p\,dp+\frac{\sqrt{n}}{b_n}
(c_n-c)\int_0^s\sigma_p^2\,dp
+\hat{M}^n_s+
\hat{\epsilon}^n_s+\frac{\sqrt{n}}{b_n}\overline{q}_s-
\frac{\sqrt{n}}{b_n}\overline{\phi}_s \Br)\wedge0
\end{equation}
In view of LD convergence in distribution at rate $b_n^2$ 
of the $\hat{M}^n$ to a continuous-path idempotent process 
 $
\lim_{\eta\to0}  
\limsup_{n\to\infty}\mathbf{P}\bl(\sup_{t\in[0,\eta]}\abs{\hat{M}^n_t}
>\delta\br)^{1/b_n^2}=0$ for $\delta>0$.
Therefore,  given $\delta\in(0,\beta)$, we derive from \eqref{eq:43},
  taking into consideration the convergences
$(\sqrt{n}/b_n)(c_n-c)\to\hat{\theta}$ and
$(\sqrt{n}/b_n)\inf_{t\in[\eta,\beta-\delta]}\overline{q}_t\to\infty$  as
$n\to\infty$, where  $\eta\in(0,\beta-\delta)$, the fact that 
$\overline{\phi}_t=0$ for $t\in[0,\beta]$, \eqref{eq:100}, \eqref{eq:92}, 
 \eqref{eq:88}, and  $\C$-exponential tightness of the $\hat{M}^n$
  that for  $\delta\in(0,\beta)$
 \begin{equation}
  \label{eq:102}
   \lim_{n\to\infty}\mathbf{P}\bl(\sup_{t\in[0,\beta-\delta]}
  \abs{Y^n_t}>\delta \br)^{1/b_n^2}=0.
\end{equation}
Putting together   \eqref{eq:77h}, \eqref{eq:100}, \eqref{eq:88},
 \eqref{eq:102}, 
the convergence $(\sqrt{n}/b_n)(c_n-c)\to\hat{\theta}$, 
and 
LD convergence in distribution at rate $b_n^2$ 
of the $\hat{M}^n$ to $\hat{M}$, we conclude 
  that the sequence of laws of the
$\hat{X}^n$ restricted to $[0,\beta-\delta]$ is 
$\C$-exponentially tight of order $b_n^2$ in $\D([0,\beta-\delta],\R)$.

We have thus established that for $\beta>0$ and 
 $ \delta\in(0,\beta\wedge(1-\beta))$ 
the processes 
$\hat{M}^n$, $\hat{L}^n$, 
$\hat{X}^n$ restricted to $[0,\beta-\delta]$,
$\hat{Y}^n$ restricted to $[\beta+\delta,1]$, and $\hat{Z}^n$
 are $\C$-exponentially tight of order $b_n^2$ in the
associated function spaces, so they are jointly exponentially tight of
order $b_n^2$ as random
elements with values in the product space. Now, 
LD convergence in distribution at rate $b_n^2$ in
 $\D_C([0,1],\R^2)\times\D_C([0,\beta-\delta],\R)\times
\D_C([\beta+\delta,1],\R)\times\D_C([0,1],\R)$ of the 
$\bl(\hat{M}^n,\hat{L}^n,(\hat{X}^n_t,\,t\in[0,\beta-\delta]),
(\hat{Y}^n_t,\,t\in[\beta+\delta,1]),\hat{Z}^n\br)$
  to $\bl(\hat{M},\hat{L},(\hat{H}_t,\,t\in[0,\beta-\delta]),
(-\hat{H}_t,\,t\in[\beta+\delta,1]),\hat{Z}\br)$
 follows 
by  \eqref{eq:77h}, \eqref{eq:80h}, \eqref{eq:100}, \eqref{eq:100a},
\eqref{eq:113}, \eqref{eq:90h}, 
\eqref{eq:103},    \eqref{eq:102}, the convergence
 $(\sqrt{n}/b_n)(c_n-c)\to\hat{\theta}$,  
 LD convergence in distribution of the $(\hat{M}^n,\hat{L}^n)$ 
to $(\hat{M},\hat{L})$, and strong uniqueness of the solution 
$(\hat{H},\hat{L})$ of \eqref{eq:113} and \eqref{eq:90h}.

Let us now assume that $\beta=0$. In
view of limits \eqref{eq:100}, \eqref{eq:88}, the convergence
 $(\sqrt{n}/b_n)(c_n-c)\to\hat{\theta}$,   and LD convergence in
distribution at rate $b_n^2$ of the $\hat{M}^n$ to $\hat{M}$, we have 
by \eqref{eq:45}  the convergence
$\lim_{\eta\to0}\limsup_{n\to\infty}
\mathbf{P}\bl(\sup_{t\in[0,\eta]}\abs{\hat{X}^n_t}>\delta\br)^{1/b_n^2}=0$
for $\delta>0$,
so by \eqref{eq:103}
$    \lim_{n\to\infty}\mathbf{P}\bl(\sup_{t\in[0,1]}
\abs{\hat{X}^n_t}>\delta\br)^{1/b_n^2}    =0.$
Therefore, by 
  \eqref{eq:77h}, the convergence
 $(\sqrt{n}/b_n)(c_n-c)\to\hat{\theta}$,   and LD convergence in
distribution at rate $b_n^2$ 
of the $\hat{M}^n$ to $\hat{M}$  the  $\hat{Y}^n$ 
LD converge in
distribution at rate $b_n^2$ in $\D_C([0,1],\R)$ to $-\hat{H}$.
\end{proof}

\begin{remark}
A slight modification of the proof shows that
for $\beta>0$ and   $\delta\in(0,\beta\wedge(1-\beta))$  
the random elements  $\hat{M}^n$, $\hat{L}^n$, 
$(\hat{X}^n_t,\,t\in[0,\beta-\delta])$, 
$(\hat{Y}^n_t,\,t\in[0,\beta-\delta])$,
$(\hat{X}^n_t,\,t\in[\beta+\delta,1])$,  
$(\hat{Y}^n_t,\,t\in[\beta+\delta,1])$, 
$\hat{X}^n_\beta$,  $\hat{Y}^n_\beta$, and $\hat{Z}$
 jointly LD converge in distribution at rate $b_n^2$ in
$\D_C([0,1],\R^2)\times\D_C([0,\beta-\delta],\R)^2\times
\D_C([\beta+\delta,1],\R)^2\times\R^2\times\D_C([0,1],\R)$
 to the respective idempotent  elements
  $\hat{M}$, $\hat{L}$, $(\hat{X}_t,\,t\in[0,\beta-\delta])$, 
$(\hat{Y}_t,\,t\in[0,\beta-\delta])$,
$(\hat{X}_t,\,t\in[\beta+\delta,1])$, 
$(\hat{Y}_t,\,t\in[\beta+\delta,1])$, $\hat{X}_\beta$, 
$\hat{Y}_\beta$, and $\hat{Z}$,
where
idempotent processes
 $\hat{X}=(\hat{X}_t,\,t\in[0,1])$ 
and $\hat{Y}=(\hat{Y}_t,\,t\in[0,1])$ are defined by
  \begin{equation*}
\begin{aligned} 
  \hat{X}_t&=\left\{
    \begin{aligned}
      \displaystyle
    &\hat{H}_t\,&&\text{ for
    }t\in[0,\beta),\\
&\hat{H}_\beta\vee0\,&&\text{ for }t=\beta,\\
&0\,&&\text{ for }t\in(\beta,1],
    \end{aligned}\right.
\end{aligned}
\quad\text{ and }\quad
\begin{aligned} 
\hat{Y}_t&=&\left\{
  \begin{aligned}
&    0\,&&\text{ for }t\in[0,\beta),\\
&(-\hat{H}_\beta)\vee0\,&&\text{ for }t=\beta,\\
\displaystyle &-\hat{H}_t&&\text{ for
    }t\in(\beta,1].
  \end{aligned}
\right.
\end{aligned}
\end{equation*}
\end{remark}
\begin{proof}[Proof of Theorem~\ref{the:6}]
The proof replicates the proof of Theorem~\ref{the:7}.
  We  begin by proving that in analogy with 
\eqref{eq:84} if $c>1$, then  as $n\to\infty$
  \begin{equation}
      \label{eq:108}
  \Bl(\frac{\sqrt{n}}{b_n}\Bl(\frac{\alpha^n}{n}-\alpha\Br),
\,\frac{\sqrt{n}}{b_n}\Bl(\frac{\beta^n}{n}-\beta\Br),
\,\frac{\sqrt{n}}{b_n}\Bl(\frac{\gamma^n}{n}-\gamma\Br)\Br) 
\xrightarrow[b_n^2]{ld}
\Bl(-\hat{H}_1,\,\frac{\hat{H}_\beta}{1-c(1-\beta)},\,\hat{Z}_\beta\Br).
  \end{equation}
As in the proof of Theorem~\ref{the:7}, we 
let $\tau^n$ be the last time $t$ before $\beta/2$ when 
$\overline{Q}^n_t=0$ and 
$\tilde{\beta}^n$ be the first time $t$ not before  $\beta/2$ when
$\overline{Q}^n_t=0$.
The  argument of the proof of Theorem~\ref{the:7} with the 
 super-exponential limits in probability of Lemma~\ref{the:idemdiff}
 used in place of Lemma~\ref{le:lln} 
implies  that under the hypotheses as $n\to\infty$
\begin{equation}
  \label{eq:109}
\tau^n\overset{\mathbf{P}^{1/b_n^2}}{\to}0,\quad
\tilde{\beta}^n\overset{\mathbf{P}^{1/b_n^2}}{\to}\beta,\quad
\mathbf{P}\Bl(\frac{\beta^n}{n}\not=
\tilde{\beta}^n-\tau^n  \Br)^{1/b_n^2}\to0,\quad
  \mathbf{P}\Bl(\frac{\gamma^n}{n}\not=\overline{E}^n_{\tilde{\beta}^n}-
\overline{E}^n_{\tau^n } \Br)^{1/b_n^2}\to0.
\end{equation}
By   \eqref{eq:39} and \eqref{eq:38} with the use of \eqref{eq:hep}
    \begin{align}
              \label{eq:53}
      -\frac{\sqrt{n}}{b_n}\,\overline{q}_{\tau^n}&
 =-c_n\int_{0}^{\tau^n}\hat{X}^n_s\,ds+
\frac{\sqrt{n}}{b_n}(c_n-c)\int_0^{\tau^n}\sigma_s^2\,ds
+\hat{M}^n_{\tau^n}+\hat{\epsilon}^n_{\tau^n}+
\hat{Y}^n_{\tau^n},\\
  \label{eq:54}
-\frac{\sqrt{n}}{b_n}\,\overline{q}_{\tilde{\beta}^n}&
=-c_n\int_{0}^{\tilde{\beta}^n}\hat{X}^n_s\,ds+
\frac{\sqrt{n}}{b_n}\,(c_n-c)\int_0^{\tilde{\beta}^n}\sigma_s^2\,ds
+\hat{M}^n_{\tilde{\beta}^n}+\hat{\epsilon}^n_{\tilde{\beta}^n}+
\frac{\sqrt{n}}{b_n}\overline{\Phi}^n_{\tilde{\beta}^n}
-\frac{\sqrt{n}}{b_n}\overline{\phi}_{\tilde{\beta}^n}.
\end{align}
 The left-most convergence in
\eqref{eq:109} implies by Lemma~\ref{the:idemdiff}, \eqref{eq:100},
and the convergence
 $(\sqrt{n}/b_n)(c_n-c)\to\hat{\theta}$   
that the right-hand side of \eqref{eq:53}
converges super-exponentially in probability at rate $b_n^2$ 
to $0$, which yields the convergence 
\begin{equation}
  \label{eq:55}
    \frac{\sqrt{n}}{b_n}\tau^n\overset{\mathbf{P}^{1/b_n^2}}{\to}0.
\end{equation}
Next,  \eqref{eq:109}, \eqref{eq:54}, and Lemma~\ref{the:idemdiff}
imply by an argument along the lines of the one 
used for deriving \eqref{eq:86} that 
\begin{equation}
  \label{eq:60}
  \frac{\sqrt{n}}{b_n}
\int_\beta^{\tilde{\beta}^n}\bl(c(1-\overline{q}_s-s)-1\br)\,ds
  -c_n\int_{0}^{\tilde{\beta}^n}\hat{X}^n_s\,ds+
\frac{\sqrt{n}}{b_n}\,(c_n-c)\int_0^{\tilde{\beta}^n}\sigma_s^2\,ds
+\hat{M}^n_{\tilde{\beta}^n}
\overset{\mathbf{P}^{1/b_n^2}}{\to}0.
\end{equation}
Also by the definition of $\hat{Y}^n$ and \eqref{eq:91}
  \begin{equation}
    \label{eq:61}
\frac{\sqrt{n}}{b_n}\Bl(\frac{\alpha^n}{n}-\alpha\Br)=\hat{Y}^n_1.
\end{equation}
Convergence \eqref{eq:108} follows by   \eqref{eq:60},
  \eqref{eq:61}, the convergence
  $(\sqrt{n}/b_n)(c_n-c)\to\hat{\theta}$,
  the joint LD convergence in distribution 
$\bl(\hat{M}^n,\hat{Y}^n_1,(\hat{X}^n_s,\,s\in[0,\beta-\delta]),\hat{Z}^n\br)
\xrightarrow[b_n^2]{ld}
\bl(\hat{M},-\hat{H}_1,(\hat{H}_s,\,s\in[0,\beta-\delta]),\hat{Z}\br)$ 
in $\D_C([0,1],\R)\times \R\times\D_C([0,\beta-\delta],\R)\times
\D_C([0,1],\R)$, 
the third super-exponential convergence in probability
in the statement of Lemma~\ref{the:idemdiff},
the last three
convergences in \eqref{eq:109}, \eqref{eq:55}, and the contraction principle.

If $c\le 1$, then
 the $Y^n_1$ LD converge in distribution to $-\hat{H}_1$ by part 2 of
Lemma~\ref{the:idemdiff}.
 
We complete the proof by showing 
 that the right-hand side of \eqref{eq:108} is idempotent
Gaussian with parameters $(\mu,\Sigma)$, i.e.,
\begin{equation}
  \label{eq:114}
\mathbf{S}\exp\Bl(-\lambda_1
\hat{H}_1+\lambda_2
\frac{\hat{H}_\beta}{1-c(1-\beta)}+\lambda_3\hat{Z}_\beta\Br)=
\exp\Bl(\lambda^T\mu+\frac{1}{2}\lambda^T\Sigma\lambda\Br),
\end{equation}
where $\lambda=(\lambda_1,\lambda_2,\lambda_3)^T\in\R^3$ and $\mathbf{S}$ denotes
idempotent expectation with respect to $\mathbf{\Pi}$.
By \eqref{eq:113}, \eqref{eq:90h}, \eqref{eq:42}, and \eqref{eq:6.4b}
\begin{align*}
\hat{H}_\beta&=
\hat{\theta}\beta e^{-\beta c}+
\sqrt{c}e^{-\beta c}\int_0^\beta e^{cs/2}\dot{\hat{W}}^{(1)}_s\,ds,\\
\hat{Z}_\beta&=
\frac{\hat{\theta}\beta^2}{2}
+\sqrt{c}\int_0^\beta\bl(1-e^{c(s-\beta)}\br)e^{-cs/2}\dot{\hat{W}}^{(1)}_s\,
ds+\sqrt{c}\int_0^\beta\sqrt{\overline{q}_s}\,\dot{\hat{W}}^{(2)}_s\,ds.
\end{align*}
On noting that by \eqref{eq:113} and \eqref{eq:6.4b}
$\hat{H}_1=\hat{H}_\beta+\hat{\theta}\int_\beta^{1}(1-s)\,ds+
\sqrt{c}\int_\beta^1\sqrt{1-s}\,\dot{\hat{W}}^{(1)}_s\,ds$, 
 $\hat{W}^{(1)}$ and
 $\hat{W}^{(2)}$ are independent, we can write using Lemma~\ref{le:wieind}
\begin{multline}
\label{eq:116}
  \mathbf{S}\exp\Bl(-\lambda_1\hat{H}_1+\lambda_2
\frac{\hat{H}_\beta}{1-c(1-\beta)}+\lambda_3\hat{Z}_\beta\Br)=
\exp\Bl(-\frac{\lambda_1\hat{\theta}(1-\beta)^2}{2}+
\Bl(\frac{\lambda_2}{1-c(1-\beta)}-\lambda_1\Br)\hat{\theta}\beta
e^{-\beta c}+\lambda_3\frac{\hat{\theta}\beta^2}{2} \Br)\\
\mathbf{S}\exp\Bl(\Bl(\frac{\lambda_2}{1-c(1-\beta)}-\lambda_1-\lambda_3\Br)
\sqrt{c}e^{-\beta c}\int_0^\beta e^{cs/2}\dot{\hat{W}}^{(1)}_s\,ds
+\lambda_3\sqrt{c}\int_0^\beta
e^{-cs/2}\dot{\hat{W}}^{(1)}_s\,ds\Br)\\
\;\mathbf{S}\exp\bl(-\lambda_1\sqrt{c}\int_\beta^1\sqrt{1-s}\,\dot{\hat{W}}^{(1)}_s\,
ds\br)\;\mathbf{S}\exp\bl(\lambda_3\sqrt{c}\int_0^\beta\sqrt{\overline{q}_s}
\,\dot{\hat{W}}^{(2)}_s\,ds \br).
\end{multline}
Lemma~\ref{le:wieind} also yields
\begin{equation}
  \begin{split}\label{eq:99}
    \mathbf{S}\exp\Bl(\sqrt{c}
\int_0^\beta\Bl(\Bl(\frac{\lambda_2}{1-c(1-\beta)}-\lambda_1-\lambda_3\Br)
 e^{cs/2-\beta c}+\lambda_3e^{-cs/2}\Br)
\dot{\hat{W}}^{(1)}_s\,ds\Br)\\
=\exp\Bl(\frac{c}{2}
\int_0^\beta\Bl(\Bl(\frac{\lambda_2}{1-c(1-\beta)}-\lambda_1-\lambda_3\Br)
e^{cs/2-\beta c}+\lambda_3e^{-cs/2}\Br)^2 \,ds\Br),
  \end{split}
\end{equation}
\begin{align}
  \label{eq:101}
  \mathbf{S}\exp\bl(-\lambda_1\sqrt{c}\int_\beta^1\sqrt{1-s}\,\dot{\hat{W}}^{(1)}_s\,
ds\br)&=\exp\Bl(\frac{c\lambda^2_1}{2}\int_\beta^1(1-s)\,ds\Br),\\
  \label{eq:104}
\mathbf{S}\exp\bl(\lambda_3\sqrt{c}\int_0^\beta\sqrt{\overline{q}_s}
\,\dot{\hat{W}}^{(2)}_s\,ds \br)&
=\exp\Bl(\frac{c\lambda_3^2}{2}\int_0^\beta\overline{q}_s\,ds \Br).
\end{align}
Equality \eqref{eq:114} follows on substituting 
\eqref{eq:99}, \eqref{eq:101}, and \eqref{eq:104}
into \eqref{eq:116} and recalling \eqref{eq:42}.
\end{proof}
\begin{remark}
  Equality \eqref{eq:114} admits also a direct proof by solving the
  variational problem on the left.
\end{remark}
\section{The critical random graph}
\label{sec:crigra}
In this section, we prove Theorem~\ref{the:crit}, so the 
  notation of the theorem is adopted.
We denote $\tilde{S}^n_t=S^n_{\lfloor n^{2/3}t\rfloor
\,\wedge  n}/n^{1/3}$,  
$\tilde{E}^n_t=E^n_{\lfloor n^{2/3}t\rfloor\wedge n}$,
 $\tilde{Q}^n_t=Q^n_{\lfloor n^{2/3}t\rfloor\wedge n}/n^{1/3}$,
$\breve{S}^n_t=S^n_{\lfloor
  (nb_n)^{2/3}t\rfloor\wedge n}/(n^{1/3}b_n^{4/3})$, 
 $\breve{E}^n_t=E^n_{\lfloor (nb_n)^{2/3}t\rfloor\wedge
  n}/b_n^2$, and 
 $\breve{Q}^n_t=Q^n_{\lfloor (nb_n)^{2/3}t\rfloor\wedge n}/(n^{1/3}b_n^{4/3})$
  for $t\in\R_+$, and introduce  processes
$\tilde{S}^n=(\tilde{S}^n_t,\,t\in\R_+)$,
$\tilde{E}^n=(\tilde{E}^n_t,\,t\in\R_+)$,
$\tilde{Q}^n=(\tilde{Q}^n_t,\,t\in\R_+)$,
$\breve{S}^n=(\breve{S}^n_t,\,t\in\R_+)$,
$\breve{E}^n=(\breve{E}^n_t,\,t\in\R_+)$,  and 
$\breve{Q}^n=(\breve{Q}^n_t,\,t\in\R_+)$.  Let stochastic processes
$\tilde{S}=(\tilde{S}_{t},\,t\in\R_+)$ and 
$\tilde{E}=(\tilde{E}_{t},\,t\in\R_+)$ be defined by the respective equalities
$\tilde{S}_t=W_t+\tilde{\theta} t-t^2/2$ and $\tilde{E}_t=
N_{\int_0^t\tilde{X}_s\,ds}$.
Let idempotent processes
$\breve{S}=(\breve{S}_{t},\,t\in\R_+)$ and 
$\breve{E}=(\breve{E}_{t},\,t\in\R_+)$ be defined by the respective equalities
$\breve{S}_t=\breve{W}_t+\breve{\theta} t-t^2/2$ and $\breve{E}_t=
\breve{N}_{\int_0^t\mathcal{R}(\breve{S})_p\,dp}$, where 
$\breve{W}=(\breve{W}_t,\,t\in\R_+)$ and $\breve{N}=(\breve{N}_t,\,t\in\R_+)$
are independent  Wiener and  Poisson idempotent
processes, respectively. The first assertion of part 1 of the next
lemma is in the theme of  Aldous \citeyear[eq. (31)]{Ald97}.
\begin{lemma}
  \label{le:crit}
  \begin{enumerate}
  \item If $n^{1/3}(c_n-1)\to\tilde{\theta}\in\R$ as $n\to\infty$, then the
   $(\tilde{S}^n,\tilde{E}^n)$  converge in distribution
  in $\D(\R_+,\R^2)$ as $n\to\infty$ to $(\tilde{S},\tilde{E})$.
If $\sqrt{n}(c_n-1)\to\theta\in\R$ as $n\to\infty$, then the 
$\bl(\sqrt{n}(\alpha^n/n-1/2),\tilde{S}^n,\tilde{E}^n\br)$ converge in
  distribution in $\R\times\D(\R_+,\R^2)$ to
$(\tilde{\alpha},\tilde{S},\tilde{E})$, where $(\tilde{S},\tilde{E})$
  correspond to $\tilde{\theta}=0$ and are independent of $\tilde{\alpha}$.
\item
If $(n^{1/3}/b_n^{2/3})(c_n-1)\to\breve{\theta}\in\R $ as $n\to\infty$,
then the   $(\breve{S}^n,\breve{E}^n)$  LD 
converge in distribution in $\D_C(\R_+,\R^2)$ at rate $b_n^2$
 to $(\breve{S},\breve{E})$.
If $(\sqrt{n}/b_n)(c_n-1)\to\hat{\theta}\in\R$ as $n\to\infty$, then the 
$\bl((\sqrt{n}/b_n)(\alpha^n/n-1/2),\breve{S}^n,\breve{E}^n\br)$ LD converge in
  distribution at rate $b_n^2$ in $\R\times\D_C(\R_+,\R^2)$ to
$(\breve{\alpha},\breve{S},\breve{E})$, where $(\breve{S},\breve{E})$
  correspond to $\breve{\theta}=0$, $\breve{\alpha}$ is idempotent
  Gaussian with parameters $(-\hat{\theta}/2,1/2)$ and is independent of 
$(\breve{S},\breve{E})$.
  \end{enumerate}
\end{lemma}
\begin{proof}
We begin with the proof of part 1.
  By \eqref{eq:3}
  \begin{equation}
    \label{eq:83}
    \tilde{S}^n_t=\tilde{M}^n_t+
n^{1/3}(c_n-1)\frac{\lfloor n^{2/3}t\rfloor\wedge n}{n^{2/3}}-
c_n\int_0^{\lfloor n^{2/3}t\rfloor\wedge n/n^{2/3}}
\frac{\lfloor n^{2/3}s\rfloor}{n^{2/3}}\,ds
-\frac{c_n}{n^{1/3}}\int_0^{\lfloor n^{2/3}t\rfloor\wedge n/n^{2/3}}
\tilde{Q}^n_s\,ds\,,
  \end{equation}
where
\begin{equation}
      \label{eq:94}
\tilde{M}^n_t=\frac{1}{n^{1/3}}\sum_{i=1}^{\lfloor n^{2/3}t\rfloor\wedge n}
\sum_{j=1}^{n-Q^n_{i-1}-(i-1)}\Bl(\xi^n_{ij}-\frac{c_n}{n}\Br).
\end{equation}
Let   $\tilde{\mathcal{F}}^n_t,\,t\in\R_+,$ denote
the $\sigma$-algebras generated by the 
$\xi^n_{ij},\,\zeta^n_{ij},\,
i=1,2,\ldots, \lfloor n^{2/3}t\rfloor\wedge n,\,j\in\mathbb{N}$, 
completed with sets of $\mathbf{P}$-measure zero. Then
$\tilde{M}^n=(\tilde{M}^n_t,\,t\in\R_+)$  
is a square-integrable martingale relative to the   filtration
$\tilde{\F}^n=(\tilde{\mathcal{F}}^n_t,\,t\in\R_+)$ with predictable
quadratic characteristic
\begin{equation}
  \label{eq:120}
    \langle
  \tilde{M}^n\rangle_t=\frac{1}{n^{2/3}}\frac{c_n}{n}\Bl(1-\frac{c_n}{n}\Br)
\sum_{i=1}^{\lfloor n^{2/3}t\rfloor\wedge n}
(n-Q^n_{i-1}-(i-1)).
\end{equation}
By 
Lemma~\ref{le:lln} $\langle  \tilde{M}^n\rangle_t\overset{\mathbf{P}}{\to}t$
as $n\to\infty$.
The predictable measure of jumps of $\tilde{M}^n$ is given by
\begin{equation*}
  \check{\nu}^n([0,t],\Gamma)=
\sum_{k=0}^{\lfloor n^{2/3}t\rfloor\wedge n-1}
\check{F}^n\Bl(1-\frac{Q^n_{k}}{n}-\frac{k}{n},\,\Gamma\setminus\{0\}
\Br),\,\Gamma\in\mathcal{B}(\R),
\end{equation*}
where 
\begin{equation*}
\check{F}^n(s,\,\Gamma')=
\mathbf{P}\Bl(\frac{1}{n^{1/3}}\sum_{j=1}^{\lfloor ns\rfloor}
\bl(\xi^n_{1j}-\frac{c_n}{n}\br)\in \Gamma'\Br),\,s\in\R_+,
\,\Gamma'\in\mathcal{B}(\R).
\end{equation*}
Therefore, for $\epsilon>0$ and $n$ large enough
\begin{multline*}
  \int_0^t\int_{\R}
  \abs{x}^2\,\ind(\abs{x}>\epsilon)\,\check{\nu}^n(ds,dx)
\le \frac{1}{\epsilon^2}
\sum_{k=1}^{\lfloor n^{2/3}t\rfloor\wedge n}\int_{\R} \abs{x}^4\,
\breve{F}^n\Bl(1-\frac{Q^n_{k-1}}{n}-\frac{k-1}{n},\,dx\Br)
\le\frac{(2c_n+3c_n^2)t}{n^{2/3}\epsilon^2},
\end{multline*}
which converges to $0$ as $n\to\infty$.    Consequently,  by Liptser and
Shiryaev \citeyear[Theorem 7.1.4]{lipshir} the processes $\tilde{M}^n$
converge in distribution in $\D(\R_+,\R)$ to the process $W$ as
$n\to\infty$.
Hence, the processes $\tilde{S}'^n=(\tilde{S}'^n_t,\,t\in\R_+)$, where 
\begin{equation*}
  \tilde{S}'^n_t=\tilde{M}^n_t+
n^{1/3}(c_n-1)\frac{\lfloor n^{2/3}t\rfloor\wedge n}{n^{2/3}}-
c_n\int_0^{\lfloor n^{2/3}t\rfloor\wedge n/n^{2/3}}
\frac{\lfloor n^{2/3}s\rfloor}{n^{2/3}}\,ds\,,
\end{equation*}
converge in distribution to the process
$\tilde{S}$.

Let $\tilde{\tilde{\epsilon}}^n=(\tilde{\tilde{\epsilon}}^n_t,\,t\in\R_+)$ be
  defined by $\tilde{\tilde{\epsilon}}^n_t=
\epsilon^n_{\lfloor n^{2/3}t\rfloor\wedge n}/n^{1/3} $.
According to \eqref{eq:33} and \eqref{eq:1a},
\begin{equation}
  \label{eq:106}
\tilde{Q}^n
=\mathcal{R}(\tilde{S}^n+\tilde{\tilde{\epsilon}}^n).
\end{equation}
Besides, by Lemma~\ref{le:eps}
\begin{equation}
  \label{eq:110}
\sup_{s\in\R_+}\abs{\tilde{\tilde{\epsilon}}^n_s}\overset{\mathbf{P}}{\to}0\;
\text{ as }n\to\infty.
\end{equation}
Since the difference
$\tilde{S}'^n_t- \tilde{S}^n_t$ is non-negative and non-decreasing in
$t$,     it follows by \eqref{eq:106} 
that the values of the process $\tilde{Q}^n$ are not greater
than the corresponding values of the reflection of $\tilde{S}'^n+
\tilde{\tilde{\epsilon}}^n$. On using that the 
$\sup_{r\in[0,t]}\abs{\tilde{S}'^n_r}$ are asymptotically bounded in
probability and that  \eqref{eq:110} holds, we conclude that the
$\sup_{s\in[0,t]}\tilde{Q}^n_s$ 
are asymptotically bounded in probability, so the right-most term of
\eqref{eq:83} tends in probability to $0$ uniformly over bounded
intervals as $n\to\infty$ implying that
  the $\tilde{S}^n$ converge in distribution to
$\tilde{S}$.  

Next, according to \eqref{eq:27}
\begin{equation}
  \label{eq:34}
  \tilde{E}^n_t=\sum_{i=1}^{\lfloor n^{2/3}t\rfloor\wedge    n}
\sum_{j=1}^{Q^n_{i-1}-1}\zeta^n_{ij}\,,\,t\in\R_+\,.
\end{equation}
   Given a sequence $\bx^n,\,n\in\N$,   of elements of $\D(\R_+,\R)$, let 
\begin{equation*}
\tilde{E}'^n_t=\sum_{i=1}^{\lfloor n^{2/3}t\rfloor\wedge n}
\;\;\sum_{j=1}^{\lfloor n^{1/3}\mathcal{R}(\bx^n)_{(i-1)/n^{2/3}}\rfloor-1}
\zeta^n_{ij},\,t\in\R_+.  
\end{equation*}
The $\tilde{E}'^n=(\tilde{E}'^n_t,\,t\in\R_+)$ are jump processes with
 $\tilde{\mathbf{F}}^n$-predictable measures of jumps
$\tilde{\nu}'^n([0,t],\Gamma)=
\sum_{i=0}^{\lfloor n^{2/3}t\rfloor\wedge
  n-1}\tilde{F}'^n\bl(\mathcal{R}(\bx^n)_{i/n^{2/3}},
\Gamma\setminus\{0\}\br),\,\Gamma\in\mathcal{B}(\R)$, where 
$\tilde{F}'^n(y,\,\Gamma')=
\mathbf{P}\bl(\sum_{j=1}^{\lfloor n^{1/3}y\rfloor-1}
\zeta^n_{1j}\in \Gamma'\br),\,\Gamma'\in\mathcal{B}(\R)$.
Theorem~VII.3.7 in Jacod and Shiryaev \citeyear{jacshir} implies that 
 if $\bx^n\to \bx$ as $n\to\infty$ in $\D(\R_+,\R)$, then the
sequence $\tilde{E}'^n,\,n\in\N$, converges in distribution 
in $\D(\R_+,\R)$ to a
compound Poisson process with compensator $\int_0^t
 \mathcal{R}(\bx)_s\,ds$. On noting that, in view of independence of 
$\tilde{S}^n$ and the $\zeta^n_{ij}$, \eqref{eq:106} and \eqref{eq:34},
the $\tilde{E}'^n$ are
 distributed according to  the regular conditional
 distributions of $\tilde{E}^n$ given that
 $\tilde{S}^n+\tilde{\tilde{\epsilon}}^n
=\bx^n$, we  conclude by  \eqref{eq:106},  \eqref{eq:110}, 
and \eqref{eq:34} that
 the $(\tilde{S}^n,\tilde{E}^n)$ jointly converge in distribution in
 $\D(\R_+,\R^2)$ to $(\tilde{S},\tilde{E})$ as $n\to\infty$. 
The first assertion of part 1
 has been proved.

For the second assertion, let in analogy with   \eqref{eq:69}
  for $\eta>0$
\begin{equation}
  \label{eq:124}
  \tilde{M}^{n,\eta}_t=
\frac{\ind(t\ge\eta)}{n}\sum_{i=\lfloor n\eta\rfloor+1}^{\lfloor nt\rfloor}
\sum_{j=1}^{n-\lfloor n\overline{q}_{(i-1)/n}\rfloor
-(i-1)}\bl(\xi^n_{ij}-\frac{c_n}{n}\br),\;
t\in[0,1],
\end{equation}
and $ \tilde{Q}^{n,\eta}=(\tilde{Q}^{n,\eta}_{t},\,t\in[0,1])$ be
defined in analogy with \eqref{eq:72} 
by the condition that it is the reflection of the process
$\int_0^{t}\bl(c_n(1-\tilde{Q}^{n,\eta}_{s}-
s)-1\br)\,ds
+\tilde{M}^{n,\eta}_t$, i.e., $\tilde{Q}^{n,\eta}_{t}\ge0$ and 
\begin{equation}
  \label{eq:125}
  \tilde{Q}^{n,\eta}_{t}=\int_0^{t}
\bl(c_n(1-\tilde{Q}^{n,\eta}_{s}-s)-1\br)\,ds
+\tilde{M}^{n,\eta}_t+\tilde{\Phi}^{n,\eta}_t,
\end{equation}
where $
\tilde{\Phi}^{n,\eta}=(\tilde{\Phi}^{n,\eta}_{t},\,t\in[0,1])$ is
non-decreasing with
$  \tilde{\Phi}^{n,\eta}_{t}=\int_0^t\ind(\tilde{Q}^{n,\eta}_{s}=0)\,
d\tilde{\Phi}^{n,\eta}_{s}$.
(For existence of $ \tilde{Q}^{n,\eta}$, one can first prove 
that a solution exists  between the jumps of $ \tilde{M}^{n,\eta}$ by
using the method of successive approximations and 
making use of Lipshitz continuity of the reflection mapping 
and Gronwall's inequality, and then account for the jumps by
introducing, if necessary, jumps in $\tilde{\Phi}^{n,\eta}$. Strong
uniqueness for $ \tilde{Q}^{n,\eta}$ follows by 
Lipshitz continuity of the reflection mapping 
and Gronwall's inequality too.)
By \eqref{eq:69}, \eqref{eq:124} and the convergence of the 
$\overline{Q}^n$ to $\overline{q}$ (Lemma~\ref{le:lln}) for $\tilde{\eta}>0$
\begin{equation*}
  \label{eq:126}
  \lim_{\eta\to0}\limsup_{n\to\infty}
\mathbf{P}\bl(\sup_{t\in[0,1]}\sqrt{n}\abs{\tilde{M}^{n,\eta}_t-
\overline{M}^{n}_t}>\tilde{\eta}\br)=0,
\end{equation*}
which implies by \eqref{eq:72}, \eqref{eq:125}, Lemma~\ref{le:eps}, 
Lipshitz continuity of
the reflection mapping, and Gronwall's inequality that
\begin{equation*}
  \lim_{\eta\to0}\limsup_{n\to\infty}
\mathbf{P}\bl(\sup_{t\in[0,1]}\sqrt{n}\abs{\tilde{Q}^{n,\eta}_t-
\overline{Q}^{n}_t}>\tilde{\eta}\br)=0,
\end{equation*}
and consequently
\begin{equation}
  \label{eq:127}
  \lim_{\eta\to0}\limsup_{n\to\infty}
\mathbf{P}\bl(\sqrt{n}\abs{\tilde{\Phi}^{n,\eta}_1-
\overline{\Phi}^{n}_1}>\tilde{\eta}\br)=0.
\end{equation}
Since $\tilde{\Phi}^{n,\eta}$ is independent of 
the $\xi^n_{ij},\,i=1,2,\ldots,\lfloor n\eta\rfloor,\,j\in\N$,  
and $\zeta^n_{ij},\,i\in\N,\,j\in\N$,  and the
$(\tilde{S}^n_t,\tilde{E}^n_t)$ are measurable functions of 
$\xi^n_{ij},\,i=1,2,\ldots,\lfloor n^{2/3}t\rfloor\wedge n,\,j\in\N$,  
and $\zeta^n_{ij},\,i=1,2,\ldots,\lfloor n^{2/3}t\rfloor\wedge n,\,j\in\N$,
it follows that $\tilde{\Phi}^{n,\eta}_1$ and
finite-dimensional distributions of the $(\tilde{S}^n,\tilde{E}^n)$
are independent for all large $n$, which yields by \eqref{eq:127}
 the asymptotic independence of 
$\sqrt{n}(\overline{\Phi}^{n}_1-
\overline{\phi}_1)$ and finite-dimensional distributions of the
$(\tilde{S}^n,\tilde{E}^n)$. 
The proof of part 1 is over.

The proof of part 2 is similar. In
analogy with \eqref{eq:83} and \eqref{eq:94}
\begin{multline}
  \label{eq:111}
      \breve{S}^n_t=\breve{M}^n_t+\frac{n^{1/3}}{b_n^{2/3}}\,(c_n-1)\,
\frac{\lfloor (nb_n)^{2/3}t\rfloor\wedge n}{(nb_n)^{2/3}}
-c_n\int_0^{\lfloor(nb_n)^{2/3}t\rfloor\wedge n/(nb_n)^{2/3}}
\frac{\lfloor   (nb_n)^{2/3}s\rfloor}{(nb_n)^{2/3}}
\,ds\\-c_n\,\frac{b_n^{2/3}}{n^{1/3}}
\int_0^{\lfloor(nb_n)^{2/3}t\rfloor\wedge n/(nb_n)^{2/3}}
\breve{Q}^n_s\,ds\,,
\end{multline}
where
\begin{equation}
\label{eq:mbreve}
  \breve{M}^n_t=\frac{1}{n^{1/3}b_n^{4/3}}
\sum_{i=1}^{\lfloor (nb_n)^{2/3}t\rfloor\wedge n}
\;\;\sum_{j=1}^{n-Q^n_{i-1}-(i-1)}\Bl(\xi^n_{ij}-\frac{c_n}{n}\Br).
\end{equation}
Let   $\breve{\mathcal{F}}^n_t,\,t\in\R_+,$ denote
the $\sigma$-algebras generated by the 
$\xi^n_{ij},\,\zeta^n_{ij},\,
i=1,2,\ldots, \lfloor (nb_n)^{2/3}t\rfloor\wedge n,\,j\in\mathbb{N}$, 
completed with sets of $\mathbf{P}$-measure zero. Then
$\breve{M}^n=(\breve{M}^n_t,\,t\in\R_+)$  
is a square-integrable martingale relative to the   filtration
$\breve{\F}^n=(\breve{\mathcal{F}}^n_t,\,t\in\R_+)$ with predictable
quadratic characteristic
\begin{equation}
  \label{eq:117}
    \langle
  \breve{M}^n\rangle_t=\frac{1}{n^{2/3}b_n^{8/3}}
\frac{c_n}{n}\Bl(1-\frac{c_n}{n}\Br)
\sum_{i=1}^{\lfloor (nb_n)^{2/3}t\rfloor\wedge n}
(n-Q^n_{i-1}-(i-1))
\end{equation}
and predictable measure of jumps 
  \begin{equation}
    \label{eq:105}
    \breve{\nu}^n([0,t],\Gamma)=
\sum_{k=0}^{\lfloor (nb_n)^{2/3}t\rfloor\wedge n-1}
\breve{F}^n\Bl(1-\frac{Q^n_{k}}{n}-\frac{k}{n},\,\Gamma\setminus\{0\}
\Br),\,\Gamma\in\mathcal{B}(\R),
\end{equation}
where 
\begin{equation}\label{eq:fbreve}
\breve{F}^n(s,\,\Gamma')=
\mathbf{P}\Bl(\frac{1}{n^{1/3}b_n^{4/3}}\sum_{j=1}^{\lfloor ns\rfloor}
\bl(\xi^n_{1j}-\frac{c_n}{n}\br)\in \Gamma'\Br),\,s\in\R_+,
\,\Gamma'\in\mathcal{B}(\R).
\end{equation}
By  \eqref{eq:117} and the first super-exponential
convergence in probability in Lemma~\ref{the:idemdiff} 
$b_n^2\langle  \breve{M}^n\rangle_t\overset{\mathbf{P}^{1/b_n^2}}{\to}t$
as $n\to\infty$. Next, in analogy with \eqref{eq:98} it is established
that
\begin{equation*}
        \lim_{n\to\infty}\mathbf{P}\Bl(\frac{1}{b_n^2}\int_0^t\int_\R 
e^{\lambda b_n^2\abs{x}}
\ind(b_n^2\abs{x}>\epsilon)\,\breve{\nu}^n(ds,dx)>\eta\Br)^{1/b_n^2}=0,
\;\lambda>0,\,\epsilon>0,\,\eta>0,\,t>0.
\end{equation*}
By Corollary~4.3.13 in Puhalskii \citeyear{Puh01} we thus have that the
$\breve{M}^n$ LD converge in $\D(\R_+,\R)$ at rate $b_n^2$ to the
idempotent process $\breve{W}$  as $n\to\infty$.
Since in analogy with \eqref{eq:106}
$  \breve{Q}^n
=\mathcal{R}(\breve{S}^n+\breve{\epsilon}^n),
$ where  
\begin{equation}
  \label{eq:107ep}
\breve{\epsilon}^n_t=
\frac{\epsilon^n_{\lfloor (nb_n)^{2/3}t\rfloor\wedge n}}{n^{1/3}b_n^{4/3}},  
\end{equation}
and 
$
\sup_{t\in\R_+}\abs{\tilde{\tilde{\epsilon}}^n_t}
\overset{\mathbf{P}^{1/b_n^2}}{\to}0\;
\text{ as }n\to\infty$ by Lemma~\ref{le:eps},
we conclude by an argument replicating the one used in the first part
of the proof that 
 the $\breve{S}^n$ LD converge to $\breve{S}$.
Finally, a ``conditional'' argument modelled on those used in the
proofs of part 1 and Corollary~\ref{co:excess}  shows that 
$(\breve{S}^n,\breve{E}^n)\xrightarrow[b_n^2]{ld}
(\breve{S},\breve{E})$ in $\D(\R_+,\R^2)$.
Convergence in  $\D_C(\R_+,\R^2)$ follows by continuity of 
$(\tilde{S},\tilde{E})$ and $(\breve{S}^n,\breve{E}^n)$ being a random
element of $\D_C(\R_+,\R^2)$. The proof of the second assertion of part 2
is similar to the proof of the second assertion of part 1.
\end{proof}
\begin{proof}[Proof of Theorem~\ref{the:crit}]
We begin with part 1, so we assume that
$n^{1/3}(c_n-1)\to\tilde{\theta}$. 
The below reasoning repeatedly invokes the property
 that for  almost every
trajectory of $\tilde{S}$ the process $\mathcal{T}(\tilde{S})$ 
is increasing in arbitrarily small neighbourhoods to the left of the initial
point and to the right of the terminal point  of an excursion of
$\mathcal{R}(\tilde{S})$; equivalently, the value of $\tilde{S}$ at the
initial point is strictly less than at any point to its left
 and   the infimum of the values of $\tilde{S}$ in
an arbitrary neighbourhood to the right of 
the terminal point  is strictly less than the
value of $\tilde{S}$ at the terminal point. 
(The stated property can be proved by using
 the  decomposition of the Wiener process into excursions, 
see, e.g., Ikeda and Watanabe
\citeyear{IkeWat}.)

We denote $\tilde{U}^n_i=U^n_i/n^{2/3}$ and 
$\tilde{R}^n_i=R^n_i/n^{2/3}$.  Given intervals 
$[\underline{u}_i,\overline{u}_i]$ and 
$[\overline{r}_i,\underline{r}_i]$, where 
$0<\underline{u}_i< \overline{u}_i$ and $0\le\underline{r}_i<
\overline{r}_i$ for $i=1,\ldots,m$, let  
$\overline{B}^{n}$     denote  the event that
there exist $m$ connected 
components of $\mathcal{G}(n,c_n/n)$ of sizes in the intervals
$[n^{2/3}\underline{u}_i,n^{2/3}\overline{u}_i]$ for $i=1,2,\ldots,m$ and
the numbers of the excess edges of these components belong to the
 respective intervals $[n^{2/3}\underline{r}_i,n^{2/3}\overline{r}_i]$.
Let $\overline{B}_{T}$ for $T>0$ denote 
  the set of functions $(\bx,\by)\in\D(\R_+,\R^2)$ 
 with $\bx_0=0$, $\by_0=0$, and $\by$  non-decreasing 
  such that    there exist non-overlapping intervals $[s_i,t_i]$
  with $t_i-s_i\in[\underline{u}_i,\overline{u}_i]$ and $t_i\le T$
 for which $\mathcal{R}(\bx)_{s_i}=
\mathcal{R}(\bx)_{t_{i}}=0$, $\mathcal{T}(\bx)_{t_{i}-}=
\mathcal{T}(\bx)_{s_{i}}$, 
and $\by_{t_{i}}-\by_{s_{i}}\in[\underline{r}_i,
\overline{r}_i]$  for $i=1,2\ldots,m$. Since the connected components
of $\mathcal{G}(n,c_n/n)$ correspond to excursions of  $\tilde{Q}^n$ 
and may occur
either before time $T$ or after it, we have
$\overline{B}^{n}
\subset\{(\tilde{S}^n+\tilde{\tilde{\epsilon}}^n,\tilde{E}^n)\in 
\overline{B}_{T}\}\cup\{\sup_{t\ge T}
(\tilde{S}^n_t+\tilde{\tilde{\epsilon}}^n_t-\tilde{S}^n_{t-\eta}
-\tilde{\tilde{\epsilon}}^n_{t-\eta})>0\} $ for
$\eta\in(0,T\wedge\min_{i=1,2,\ldots,m}\underline{u}_i)$.
Since   the set
$\overline{B}_{T}$ and its   closure (in $\D(\R_+,\R^2)$)
have the same intersection with $\C(\R_+,\R^2)$,
Lemma~\ref{le:crit} implies that 
\begin{equation}
  \label{eq:119}
\limsup_{n\to\infty}
\mathbf{P}\bl(\overline{B}^{n}\br)\le
\mathbf{P}\bl((\tilde{S},\tilde{E})\in\overline{B}_{T}\br)  +
\limsup_{n\to\infty}\mathbf{P}\bl(\sup_{t\ge T}
(\tilde{S}^n_t+\tilde{\tilde{\epsilon}}^n_t-\tilde{S}^n_{t-\eta}
-\tilde{\tilde{\epsilon}}^n_{t-\eta})>0\br).
\end{equation}
We show that
\begin{equation}
  \label{eq:128}
  \lim_{T\to\infty}\limsup_{n\to\infty}\mathbf{P}\bl(\sup_{t\ge T}
(\tilde{S}^n_t+\tilde{\tilde{\epsilon}}^n_t-\tilde{S}^n_{t-\eta}
-\tilde{\tilde{\epsilon}}^n_{t-\eta})>0\br)=0.
\end{equation}
By \eqref{eq:83},  \eqref{eq:120},  \eqref{eq:35}, 
 and Doob's inequality for all $n$ and $T$ 
large enough
\begin{multline}
      \label{eq:121}
  \mathbf{P}\bl(\sup_{t\ge T}
(\tilde{S}^n_t+\tilde{\tilde{\epsilon}}^n_t-\tilde{S}^n_{t-\eta}
-\tilde{\tilde{\epsilon}}^n_{t-\eta})>0\br)\\\le
\sum_{k=0}^\infty  \mathbf{P}\bl(\sup_{t\in[T+k \eta,T+(k+1)\eta]}
(\tilde{M}^n_t-\tilde{M}^n_{t-\eta}
+\tilde{\tilde{\epsilon}}^n_t
-\tilde{\tilde{\epsilon}}^n_{t-\eta})>c_n\eta (T+(k-1)\eta)
-2\eta n^{1/3}\abs{c_n-1}\br)\\
\le
\sum_{k=0}^\infty  \mathbf{P}\bl(2\sup_{s\in[0,\eta]}
(\abs{\tilde{M}^n_{T+(k-1)\eta+s}-\tilde{M}^n_{T+(k-1)\eta}}+
\abs{\tilde{\tilde{\epsilon}}^n_{T+(k-1)\eta+s}
-\tilde{\tilde{\epsilon}}^n_{T+(k-1)\eta}})
\\+\sup_{s\in[0,\eta]}
(\abs{\tilde{M}^n_{T+k\eta+s}-\tilde{M}^n_{T+k\eta}}+
\abs{\tilde{\tilde{\epsilon}}^n_{T+k\eta+s}
-\tilde{\tilde{\epsilon}}^n_{T+k\eta}})
>c_n\eta (T+(k-1)\eta)
-2\eta n^{1/3}\abs{c_n-1}\br)\\
\le
\sum_{k=-1}^\infty  \mathbf{P}\Bl(\sup_{s\in[0,\eta]}
(\abs{\tilde{M}^n_{T+k\eta+s}-\tilde{M}^n_{T+k\eta}}
+\abs{\tilde{\tilde{\epsilon}}^n_{T+k\eta+s}
-\tilde{\tilde{\epsilon}}^n_{T+k\eta}})>\frac{c_n\eta}{3} (T+k\eta)
-\frac{2\eta}{3}n^{1/3}\abs{c_n-1}\Br)\\
+\sum_{k=-1}^\infty  \mathbf{P}\Bl(\sup_{s\in[0,\eta]}
(\abs{\tilde{M}^n_{T+(k+1)\eta+s}-\tilde{M}^n_{T+(k+1)\eta}}
+\abs{\tilde{\tilde{\epsilon}}^n_{T+(k+1)\eta+s}
-\tilde{\tilde{\epsilon}}^n_{T+(k+1)\eta}})>\frac{c_n\eta}{3} (T+k\eta)
-\frac{2\eta}{3}n^{1/3}\abs{c_n-1}\Br)
\\
\le
2\sum_{k=-1}^\infty \frac{\displaystyle
4\mathbf{E}(\langle\tilde{M}^n\rangle_{T+(k+1)\eta}-
\langle\tilde{M}^n\rangle_{T+k\eta})
+\mathbf{E}\sup_{s\in[0,\eta]}
\abs{\tilde{\tilde{\epsilon}}^n_{T+k\eta+s}
-\tilde{\tilde{\epsilon}}^n_{T+k\eta}}^2
}{\Bl(\dfrac{c_n\eta}{3} (T+k\eta)
-\dfrac{2\eta}{3}n^{1/3}\abs{c_n-1}\Br)^2}\\
+2\sum_{k=-1}^\infty \frac{\displaystyle
4\mathbf{E}(\langle\tilde{M}^n\rangle_{T+(k+2)\eta}-
\langle\tilde{M}^n\rangle_{T+(k+1)\eta})
+\mathbf{E}\sup_{s\in[0,\eta]}
\abs{\tilde{\tilde{\epsilon}}^n_{T+(k+1)\eta+s}
-\tilde{\tilde{\epsilon}}^n_{T+(k+1)\eta}}^2
}{\Bl(\dfrac{c_n\eta}{3} (T+k\eta)
-\dfrac{2\eta}{3}n^{1/3}\abs{c_n-1}\Br)^2}
\\\le
4\sum_{k=-1}^\infty \frac{8c_n\eta+\bl((1+c_n)^2+
  c_n\br)n^{-2/3}}{\Bl(\dfrac{c_n\eta}{3} (T+k\eta)
-\dfrac{2\eta}{3}n^{1/3}\abs{c_n-1}\Br)^2}
\le4\sum_{k=-1}^\infty \frac{12^2
(16\eta+1)}{(T+k\eta)^2\eta^2}.
\end{multline}
The latter sum converges to $0$ as $T\to\infty$, so 
\eqref{eq:128} follows.

Denoting
$\overline{B}=\cup_{T>0}\overline{B}_{T}$ we
deduce from \eqref{eq:119} and \eqref{eq:128} that
$  \limsup_{n\to\infty}\mathbf{P}\bl(\overline{B}^{n}\br)\le
 \mathbf{P}\bl((\tilde{S},\tilde{E})\in\overline{B}\br)$.
By the cited property,
 for almost all $\omega\in\Omega$ any interval $[s,t]$  such that 
$\mathcal{R}(\tilde{S})_s(\omega)=\mathcal{R}(\tilde{S})_t(\omega)=0$ and 
$\mathcal{T}(\tilde{S})_s(\omega)=\mathcal{T}(\tilde{S})_t(\omega)$
 is an excursion of $\mathcal{R}(\tilde{S})(\omega)$.
Therefore,  $ \mathbf{P}\bl((\tilde{S},\tilde{E})\in\overline{B}\br)=
P_{\{[\underline{u}_i,\overline{u}_i],
[\overline{r}_i,\underline{r}_i]\}_{i=1}^m}$, where
$P_{\{[\underline{u}_i,\overline{u}_i],
[\overline{r}_i,\underline{r}_i]\}_{i=1}^m}$ denotes the probability that 
there exist $m$ excursions of $\mathcal{R}(\tilde{S})=\tilde{X}$ with
 lengths in the respective intervals $[\underline{u}_i,\overline{u}_i]$  and
the increments of $\tilde{E}$ over these excursions belong to the
 respective intervals $[\underline{r}_i,\overline{r}_i]$. Hence,
 \begin{equation}
   \label{eq:115}
   \limsup_{n\to\infty}\mathbf{P}(\overline{B}^{n})\le
P_{\{[\underline{u}_i,\overline{u}_i],
[\overline{r}_i,\underline{r}_i]\}_{i=1}^m}.
 \end{equation}
Next, let  $\overset{o}{B^n}$   denote  the event that
there exist $m$ connected 
components of $\mathcal{G}(n,c_n/n)$ of sizes in the segments
$(n^{2/3}\underline{u}_i,n^{2/3}\overline{u}_i)$ for $i=1,2,\ldots,m$ and
the numbers of the excess edges of these components belong to the
 respective segments $(n^{2/3}\underline{r}_i,n^{2/3}\overline{r}_i)$.
Let $\overset{o}{B}$  denote 
  the set of functions $(\bx,\by)\in\D(\R_+,\R^2)$ 
  for which    there exist disjoint intervals $[s_i,t_i]$
  with $t_i-s_i\in(\underline{u}_i,\overline{u}_i)$  
 such that $\bx_{s_i}=
\bx_{t_{i}}<\inf_{p\in[0,(s_i-\eta)^+]}\bx_p$ and
$\bx_{t_{i}}>\inf_{p\in[t_i,t_i+\eta]}\bx_p$ 
for arbitrary $\eta>0$,  $\bx_p>\bx_{s_i}$ for
$p\in(s_i,t_i)$, 
and $\by_{t_{i}}-\by_{s_{i}}\in(\underline{r}_i,\overline{r}_i)$  
for $i=1,2\ldots,m$.
 Since continuous functions from
$\overset{o}{B}$ are interior points of 
$\overset{o}{B}$ and
$\{(\tilde{S}^n+\tilde{\tilde{\epsilon}}^n,\tilde{E}^n)\in 
\overset{o}{B}\}
\subset  \overset{o}{B^n}$, by Lemma~\ref{le:crit}
$\liminf_{n\to\infty}\mathbf{P}\bl(\overset{o}{B^n}\br)\ge 
\mathbf{P}\bl((\tilde{S},\tilde{E})\in 
\overset{o}{B}\br)$. If a sample event
$\omega\in\Omega$ is such that $\tilde{X}(\omega)$  has 
$m$ excursions of  lengths in
the respective segments $(\underline{u}_i,\overline{u}_i)$  and
the increments of $\tilde{E}(\omega)$ over these excursions belong to the
 respective segments $(\underline{r}_i,\overline{r}_i)$, then 
by the cited property 
$(\tilde{X}(\omega),\tilde{E}(\omega))\in
\overset{o}{B}$ with probability 1.
Therefore, denoting the probability of the set of these $\omega$ as
 $P_{\{(\underline{u}_i,\overline{u}_i),
(\overline{r}_i,\underline{r}_i)\}_{i=1}^m}$, we deduce that
\begin{equation}
  \label{eq:112}
     \liminf_{n\to\infty}\mathbf{P}\bl(\overset{o}{B^n}\br)\ge 
P_{\{(\underline{u}_i,\overline{u}_i),
(\overline{r}_i,\underline{r}_i)\}_{i=1}^m}.
 \end{equation}
The asserted in part 1 of the theorem convergence of 
$(\tilde{U}^n,\tilde{R}^n)$ follows by \eqref{eq:115},
 \eqref{eq:112}, and the observation that the right-hand sides of these
 inequalities coincide.
The assertion of the theorem for the case
 $\sqrt{n}(c_n-1)\to\theta$ follows by a similar argument with the use
 of part 1 of Theorem~\ref{the:7} and the second assertion of
 part 1 of Lemma~\ref{le:crit}.

The proof of part 2 is obtained by combining the approaches 
 of the proofs of part 1 and Theorem~\ref{the:jointld}.
We firstly note 
that the action functional $\breve{I}^{S,E}(\bx,\by)$ associated with
 $(\breve{S},\breve{E})$ is of the form
$\breve{I}^{S,E}(\bx,\by)=\int_0^\infty(\dot{\bx}_t-\breve{\theta}+t)^2\,dt/2
 +\int_0^\infty
\pi\bl(\dot{\by}_t/\mathcal{R}(\bx)_t\br)\,\mathcal{R}(\bx)_t\,dt$ if
 $\bx$ and $\by$ are absolutely continuous with $\bx_0=\by_0=0$ and
 $\by$ non-decreasing, and $\breve{I}^{S,E}(\bx,\by)=\infty$ otherwise.
Then the proof is carried out along the lines of 
the proof of Theorem~\ref{the:jointld}, where 
the proof of an analogue of Lemma~\ref{le:konechngia} uses
  parts 2 of Lemmas~\ref{le:var} and 
\ref{le:optim}  instead of  respective parts 1 of these lemmas.
In addition, the proof of  an  analogue of \eqref{eq:51'}, as in the argument
 just given, uses the convergence
 \begin{equation}\label{eq:modhvo}
   \lim_{T\to\infty}
 \limsup_{n\to\infty}\mathbf{P}\bl(\sup_{t\ge T}
(\breve{S}^n_t+\breve{\epsilon}^n_t-\breve{S}^n_{t-\eta}
-\breve{\epsilon}^n_{t-\eta})>0\br)^{1/b_n^2}=0,\;\eta>0.
\end{equation}
We omit most of
 the details  and only show the latter.
Arguing as in \eqref{eq:121} 
\begin{multline}\label{eq:hvomod}
  \mathbf{P}\bl(\sup_{t\ge T}
(\breve{S}^n_t+\breve{\epsilon}^n_t-\breve{S}^n_{t-\eta}
-\breve{\epsilon}^n_{t-\eta})>0\br)^{1/b_n^2}\\
\le
\sum_{k=-1}^\infty  \mathbf{P}\Bl(\sup_{s\in[0,\eta]}
(\abs{\breve{M}^n_{T+k\eta+s}-\breve{M}^n_{T+k\eta}}
+\abs{\breve{\epsilon}^n_{T+k\eta+s}
-\breve{\epsilon}^n_{T+k\eta}})>\frac{c_n\eta}{3} (T+k\eta)
-\frac{2\eta}{3}\frac{n^{1/3}}{b_n^{2/3}}\abs{c_n-1}\Br)^{1/b_n^2}\\
+\sum_{k=-1}^\infty  \mathbf{P}\Bl(\sup_{s\in[0,\eta]}
(\abs{\breve{M}^n_{T+(k+1)\eta+s}-\breve{M}^n_{T+(k+1)\eta}}
+\abs{\breve{\epsilon}^n_{T+(k+1)\eta+s}
-\breve{\epsilon}^n_{T+(k+1)\eta}})>\frac{c_n\eta}{3} (T+k\eta)
-\frac{2\eta}{3}\frac{n^{1/3}}{b_n^{2/3}}\abs{c_n-1}\Br)^{1/b_n^2}\\
\le\biggl(\sum_{i=1}^2 \sup_{t\in\R_+}
\Bl(\mathbf{E}\exp\bl((-1)^ib_n^2
(\breve{M}^n_{t+\eta}-\breve{M}^n_{t})\br)\Br)^{1/b_n^2}
+\Bl(\sup_{t\in\R_+} \mathbf{E}\exp\bl(b_n^2
\sup_{s\in[0,\eta]}\abs{\epsilon^n_{t+s}
-\epsilon^n_{t}}\br)\Br)^{1/b_n^2}\biggr)\\
\sum_{k=-1}^\infty
\exp\Bl(-\Bl(\dfrac{c_n\eta}{6} (T+k\eta)
-\dfrac{\eta}{3}\frac{n^{1/3}}{b_n^{2/3}}\abs{c_n-1}\Br)\Bl).
\end{multline}
Let $\breve{\mathcal{E}}^n_t(\lambda),\,t\in\R_+,\,\lambda\in\R$,
denote the stochastic exponential of $\breve{M}^n$ so that by
\eqref{eq:105} and \eqref{eq:fbreve}
\begin{equation*}
\log  \breve{\mathcal{E}}^n_t(\lambda)=
n\log\mathbf{E}\exp\Bl(\frac{\lambda}{n^{1/3}b_n^{4/3}}
\bl(\xi^n_{11}-\frac{c_n}{n}\br)\Br)
\sum_{k=0}^{\lfloor (nb_n)^{2/3}t\rfloor\wedge n-1 }
 \Bl(1-\frac{{Q}^n_k}{n}-\frac{k}{n}\Br).  
\end{equation*}
Hence, for $t\in\R_+$ and $n$ large enough,
\begin{equation*}
\frac{1}{b_n^2}\log \mathbf{E}\exp\bl(\pm b_n^2
(\breve{M}^n_{t+\eta}-\breve{M}^n_{t})\br)\le
\frac{1}{2b_n^2}\log\mathbf{E}\frac{\breve{\mathcal{E}}^n_{t+\eta}
(\pm 2b_n^2)}{\breve{\mathcal{E}}^n_t(\pm 2b_n^2)}
\le
\frac{n^{5/3}\eta}{b_n^{4/3}}
\Bl(
\log \mathbf{E}\exp\Bl(\pm\frac{2b_n^{2/3}}{n^{1/3}}\,
\bl(\xi^n_{11}-\frac{c_n}{n}\br)\Br)\Br)^+,
\end{equation*}
so since $\log \mathbf{E}\exp\bl(\pm 2b_n^{2/3}(\xi^n_{11}-c_n/n)/n^{1/3}\br)$ 
is asymptotically equivalent to $2 c_nb_n^{4/3}/n^{5/3}$ as
$n\to\infty$, we conclude that
\begin{equation}
\label{eq:122}
\limsup_{n\to\infty}
\sup_{t\in\R_+}\Bl(\mathbf{E}\exp\bl(\pm b_n^2
(\breve{M}^n_{t+\eta}-\breve{M}^n_{t})\br)\Br)^{1/b_n^2}\le
e^{2\eta}.
\end{equation}
Also by \eqref{eq:35} and the definition of $\breve{\epsilon}^n_t$
in \eqref{eq:107ep}
\begin{equation}
  \label{eq:123}
\limsup_{n\to\infty}\sup_{t\in\R_+}
\Bl(\mathbf{E}\exp\bl( b_n^2\sup_{s\in[0,\eta]}
\abs{\breve{\epsilon}^n_{t+s}
-\breve{\epsilon}^n_{t}}\br)\Br)^{1/b_n^2}\le 1.
\end{equation}
Limit \eqref{eq:modhvo} follows by \eqref{eq:hvomod}, \eqref{eq:122},
\eqref{eq:123}, and the convergence 
$\bl(n^{1/3}/b_n^{2/3}\br)(c_n-1)\to\breve{\theta}$.
\end{proof}
Corollary~\ref{co:crit} follows by the contraction principle, in
particular, part 2 is proved in analogy with part 2 of 
Corollary~\ref{co:oc}. (Note that in the expression for 
$\breve{I}^\beta_{\breve{\theta}}$ the role of $K_c(u)$ and $L_c(u)$
are played by the functions $-u^3/24$ and
$\bl((u-\breve{\theta})^3+\breve{\theta}^3\br)/6$, respectively, and
an analogue of Lemma~\ref{le:min} holds with $2(\breve{\theta}-u)$ as
$u^\ast$.) 

{\bf Acknowledgement.} The author is 
thankful to the referees and Neil O'Connell
for bringing to his attention  the papers by Barraez, Boucheron, and
Fernandez de la Vega \citeyear{BarBouFer00} 
and by Bollob{\'a}s, Grimmett, and Janson \citeyear{BolGriJan96}.

\appendix

\section{Summary of idempotent probability}
\label{sec:appendix}

This appendix relates some  facts of idempotent probability theory.
More detailed exposition is given in Puhalskii \citeyear{Puh01}.

Let $\Upsilon$ be a set. A  function $\mathbf{\Pi}$ from the power
set of $\Upsilon$ to $[0,1]$ is called an idempotent probability
 if $\mathbf{\Pi}(\Gamma)=\sup_{\upsilon\in\Gamma}
\mathbf{\Pi}(\{\upsilon\}),\,\Gamma\subset \Upsilon$ and
$\mathbf{\Pi}(\Upsilon)=1$. If in addition, $\Upsilon$ is a metric space and 
 the sets $\{\upsilon\in\Upsilon:\,\mathbf{\Pi}(\upsilon)\ge a\}$ are compact for all
$a\in(0,1]$, then $\mathbf{\Pi}$ is called a deviability.  
Obviously, $\mathbf{\Pi}$ is a deviability if and only if
$I(\upsilon)=-\log\mathbf{\Pi}(\{\upsilon\})$ is an action functional.
Below, we denote $\mathbf{\Pi}(\upsilon)=\mathbf{\Pi}(\{\upsilon\})$ and assume unless mentioned
otherwise that $\mathbf{\Pi}$ is an idempotent probability on $\Upsilon$.
A property $\mathcal{P}(\upsilon),\,\upsilon\in\Upsilon,$ 
pertaining to the elements of $\Upsilon$ is said to hold
$\mathbf{\Pi}$-a.e. if $\mathbf{\Pi}(\mathcal{P}(\upsilon)\text{ does not hold})=0$.
A $\tau$-algebra on $\Upsilon$ is defined as a subset of the power set of
$\Upsilon$ for which there exists a partitioning of $\Upsilon$ into
disjoint sets such that
every element of $\mathcal{A}$ is a union of the elements of
the partitioning. We call the elements of the partitioning  the
  atoms of $\mathcal{A}$ and denote as $[\upsilon]$ the atom containing
$\upsilon$. The power set of $\Upsilon$ is called the discrete $\tau$-algebra.
A $\tau$-algebra $\mathcal{A}$ is called complete (or $\mathbf{\Pi}$-complete,
or complete with
respect to $\mathbf{\Pi}$ if idempotent probability needs to be specified) if each
one-point set $\{\upsilon\}$ with $\mathbf{\Pi}(\upsilon)=0$ is an atom of
$\mathcal{A}$; the completion (or the $\mathbf{\Pi}$-completion,
or the completion with respect to $\mathbf{\Pi}$ 
if idempotent probability needs to be specified)
of a $\tau$-algebra $\mathcal{A}$ is
defined as the $\tau$-algebra obtained by taking as the
 atoms the points of  idempotent probability $0$ and
set-differences of the atoms of $\mathcal{A}$ and
sets of idempotent probability $0$; the completion of a $\tau$-algebra
is a complete $\tau$-algebra.
  If $\Upsilon'$ is another set equipped with idempotent probability
  $\mathbf{\Pi}'$ and $\tau$-algebra $\mathcal{A}'$, then the
product idempotent probability $\mathbf{\Pi}\times \mathbf{\Pi}'$ on $\Upsilon\times\Upsilon'$ is
defined by $(\mathbf{\Pi}\times\mathbf{\Pi}')(\upsilon,\upsilon')
=\mathbf{\Pi}(\upsilon)\mathbf{\Pi}'(\upsilon')$ for
$(\upsilon,\upsilon')\in\Upsilon\times\Upsilon'$, the product $\tau$-algebra
$\mathcal{A}\otimes \mathcal{A}'$ is defined as having the atoms 
$[\upsilon]\times[\upsilon']$, 
where $\upsilon\in\Upsilon$ and $\upsilon'\in\Upsilon'$.

A function $f$ from a set
$\Upsilon$ equipped with idempotent probability $\mathbf{\Pi}$ to a set $\Upsilon'$ 
is called an idempotent variable. 
If $\Upsilon$ and $\Upsilon'$ are equipped with
$\tau$-algebras $\mathcal{A}$ and $\mathcal{A}'$, respectively, 
the idempotent variable $f$ is said to be
$\mathcal{A}/\mathcal{A}'$-measurable, or simply measurable if the
$\tau$-algebras are understood, if $f^{-1}([\upsilon'])\in\mathcal{A}$
for any $\upsilon'\in\Upsilon'$.
We say that $f$ is $\mathcal{A}$-measurable if it is 
measurable for the discrete $\tau$-algebra on $\Upsilon'$.
The $\tau$-algebra of $\Upsilon$ generated by $f$ is defined by the atoms
$\{\upsilon\in\Upsilon:\,f(\upsilon)=\upsilon'\},
\,\upsilon'\in \Upsilon'$. 
The idempotent variable $f$
is thus  $\mathcal{A}$-measurable if 
$\{\upsilon\in\Upsilon:\,f(\upsilon)=\upsilon'\}\in\mathcal{A}$ for all
$\upsilon'\in \Upsilon'$.  
As in probability theory, we
routinely omit the argument $\upsilon$ 
in the notation for an idempotent variable.
 The idempotent distribution of an idempotent variable $f$ is defined
 as the set function $\mathbf{\Pi}\circ
f^{-1}(\Gamma)=\mathbf{\Pi}(f\in\Gamma),\,\Gamma\subset \Upsilon'$;
it is also called the image of $\mathbf{\Pi}$ under
$f$. If $\Upsilon$ is a metric space, 
$\mathbf{\Pi}$ is 
a deviability on $\Upsilon$, and $f$ is a continuous mapping from
$\Upsilon$ to a metric space $\Upsilon'$, then 
 $\mathbf{\Pi}\circ f^{-1}$ is a deviability on $\Upsilon'$.
In particular, if $\Upsilon'\subset\Upsilon$ with induced metric and 
$\mathbf{\Pi}(\Upsilon\setminus\Upsilon')=0$, then 
the restriction $\mathbf{\Pi}|_{\Upsilon'}$ of $\mathbf{\Pi}$ to
$\Upsilon'$  defined by $\mathbf{\Pi}|_{\Upsilon'}(\upsilon)
=\mathbf{\Pi}(\upsilon)$ for $\upsilon\in\Upsilon'$ 
is a deviability on $\Upsilon'$.
In general,  $f$ is said to be  Luzin if  $\mathbf{\Pi}\circ f^{-1}$ is a
deviability on $\Upsilon'$. 

Subsets $A$ and $A'$ of $\Upsilon$ 
are said to be independent 
if $\mathbf{\Pi}(A\cap A')=\mathbf{\Pi}(A)\mathbf{\Pi}(A')$; $\tau$-algebras
$\mathcal{A}$ and $\mathcal{A}'$ are said to be independent if events $A$ and
$A'$ are independent for any $A\in\mathcal{A}$ and
$A'\in\mathcal{A}'$;
$\Upsilon'$-valued idempotent 
variables $f$ and $f'$ are said to be independent if 
$\mathbf{\Pi}(f=\upsilon',\,f'=\upsilon'')=\mathbf{\Pi}(f=\upsilon')
\mathbf{\Pi}(f'=\upsilon'')$ for all $\upsilon',\upsilon''\in \Upsilon'$. An
idempotent variable $f$ and a $\tau$-algebra
$\mathcal{A}$ are said to be independent (or $f$ to be independent of
$\mathcal{A}$) if the $\tau$-algebra generated by $f$ and
$\mathcal{A}$ are independent. 
If $f$ is $\R_+$-valued, the idempotent expectation of $f$ is defined
by $\mathbf{S}f=\sup_{\upsilon\in\Upsilon}f(\upsilon)
\mathbf{\Pi}(\upsilon)$, it is also denoted as $\mathbf{S}_\mathbf{\Pi} f$
if the reference idempotent probability needs to be indicated. The
following analogue of the Markov inequality holds: $\mathbf{\Pi}(f\ge a)\le
\mathbf{S}f/a$, where $a>0$. If
$\R_+$-valued idempotent variables $f$ and $f'$ are independent, then
$\mathbf{S}(ff')=\mathbf{S}f\, \mathbf{S}f'$.
An $\R_+$-valued idempotent variable $f$ is
said to be maximable if $\lim_{b\to\infty}\mathbf{S}(f\ind(f>b))=0$. 
A collection $f_\alpha$ of $\R_+$-valued idempotent variables is called
uniformly maximable if $\lim_{b\to\infty}\sup_\alpha 
\mathbf{S}(f_\alpha\ind(f_\alpha>b))=0$. 
The conditional idempotent expectation of an $\R_+$-valued idempotent
variable $f$ given a $\tau$-algebra $\mathcal{A}$ is defined as
\begin{equation*}
       \mathbf{S}(f|\mathcal{A})(\upsilon)=
       \begin{cases}
\displaystyle         \sup_{\upsilon'\in[\upsilon]}
f(\upsilon')
\frac{\mathbf{\Pi}(\upsilon')}{\mathbf{\Pi}([\upsilon])}
&\text{ if }\mathbf{\Pi}([\upsilon])>0,\\
f'(\upsilon),&\text{  if }\mathbf{\Pi}([\upsilon])=0,
\end{cases}
\end{equation*}
where $f'(\upsilon)$ is an $\R_+$-valued function constant on the atoms of
$\mathcal{A}$. 
Conditional 
idempotent expectation is thus specified $\mathbf{\Pi}$-a.e. It has many of the
properties  of conditional expectation, 
in particular, $\mathbf{S}(f|\mathcal{A})$ is
$\mathcal{A}$-measurable, if $f$ is
$\mathcal{A}$-measurable then
$\mathbf{S}(f|\mathcal{A})=f$ $\mathbf{\Pi}$-a.e., 
and if $f$ and $\mathcal{A}$ are
independent then $\mathbf{S}(f|\mathcal{A})=\mathbf{S}f$ $\mathbf{\Pi}$-a.e.,
 Puhalskii \citeyear[Lemma 1.6.21]{Puh01}.
If for an $\R^d$-valued idempotent variable $f$ the conditional
idempotent expectation $\mathbf{S}(\exp(\lambda^T f)|\mathcal{A})$ is
$\mathbf{\Pi}$-a.e. constant on $\Upsilon$ for all $\lambda\in\R^d$
 and  is an essentially smooth
function of $\lambda$, then $f$ and $\mathcal{A}$ are independent,
Puhalskii \citeyear[Corollary 1.11.9]{Puh01}.

An $\R^d$-valued idempotent variable $f$ on $(\Upsilon,\mathbf{\Pi})$ 
is said to be Gaussian with
parameters $(m,\Sigma)$, where $m\in\R^d$ and $\Sigma$ is a positive
semi-definite $d\times d$ matrix, if
$\mathbf{S}\exp(\lambda^Tf)=\exp(\lambda^Tm+\lambda^T\Sigma\lambda/2)$ for all
$\lambda\in\R^d$. Equivalently,
$\mathbf{\Pi}(f=z)=\exp(-(z-m)^T\Sigma^\oplus(z-m)/2)$ if $z-m$ belongs to the
range of $\Sigma$ and $\mathbf{\Pi}(f=z)=0$ otherwise, where $\Sigma^\oplus$
denotes the pseudo-inverse of $\Sigma$,
 Puhalskii \citeyear[Lemma 1.11.12]{Puh01}.

A flow of $\tau$-algebras, or a $\tau$-flow, on $\Upsilon$ is defined
as a collection $\mathbf{A}=(\mathcal{A}_t,\,t\in\R_+)$ 
of $\tau$-algebras on $\Upsilon$ such that
$\mathcal{A}_s\subset\mathcal{A}_t$ for $s\le t$; the latter condition
is equivalent to the atoms of $\mathcal{A}_s$ being unions of the
atoms of $\mathcal{A}_t$. A $\tau$-flow is called complete if it
consists of complete $\tau$-algebras, the completion of a $\tau$-flow
is obtained by completing its $\tau$-algebras; the completion of a
$\tau$-flow is a complete $\tau$-flow. An idempotent variable 
$\sigma:\,\Upsilon\to\R_+$ is called an idempotent
$\mathbf{A}$-stopping time, or a stopping  time relative to
$\mathbf{A}$, if $\{\upsilon:\,\sigma(\upsilon)=t\}\in\mathcal{A}_t$ 
for $t\in\R_+$.
Given a $\tau$-flow $\mathbf{A}$ and an idempotent $\mathbf{A}$-stopping time
$\sigma$, we define  $\mathcal{A}_\sigma$ as the $\tau$-algebra with
atoms $[\upsilon]_{\mathcal{A}_{\sigma(\upsilon)}}$. 
 If $\Upsilon=\C(\R_+,\R^d)$, the canonical
$\tau$-flow is the $\tau$-flow $\mathbf{C}=(\mathcal{C}_t,\,t\in\R_+)$ 
with the $\mathcal{C}_t$ having the atoms 
$p_t^{-1}\bx,\,\bx\in\C(\R_+,\R^d),$ where $p_t:\,\C(\R_+,\R^d)
\to\C(\R_+,\R^d)$ 
is defined by $(p_t\bx)_s=\bx_{s\wedge t},\,s\in\R_+$.

A collection $(X_t,\,t\in\R_+)$ 
of $\R^d$-valued idempotent variables on $\Upsilon$ 
is called an idempotent process.
The functions $(X_t(\upsilon),\,t\in\R_+)$ 
for various $\upsilon\in\Upsilon$ are
called trajectories (or paths) of $X$.
An idempotent process $(X_t,\,t\in\R_+)$ is said to be 
$\mathbf{A}$-adapted if the $X_t$ are $\mathcal{A}_t$-measurable
for $t\in\R_+$. If $(X_t,\,t\in\R_+)$ is $\mathbf{A}$-adapted with
unbounded above continuous paths, 
then $\sigma=\inf\{t\in\R_+:\,X_t\ge a\}$, where $a\in\R$, is an
idempotent $\mathbf{A}$-stopping time,
  Puhalskii \citeyear[Lemma 2.2.18]{Puh01}.
 If $\Upsilon=\C(\R_+,\R^d)$, the canonical idempotent
process is defined by $X_t(\bx)=\bx_t$.
An $\mathbf{A}$-adapted $\R_+$-valued idempotent
process $M=(M_t,\,t\in\R_+)$ is said to be 
an $\mathbf{A}$-exponential maxingale,
or an exponential maxingale relative to $\mathbf{A}$, if the $M_t$ are
maximable and $\mathbf{S}(M_t|\mathcal{A}_s)=M_s$ 
$\mathbf{\Pi}$-a.e. for $s\le t$. If,
in addition, the collection $M_t,\,t\in\R_+,$ is uniformly maximable,
then $M$ is said to be a uniformly maximable exponential maxingale.
An $\mathbf{A}$-adapted $\R_+$-valued idempotent
process $M=(M_t,\,t\in\R_+)$ is called 
an $\mathbf{A}$-local exponential maxingale,
or a local exponential maxingale relative to $\mathbf{A}$, if there
exists a sequence $\tau_n$ of idempotent $\mathbf{A}$-stopping times such that
$\tau_n\uparrow\infty$ as $n\to\infty$ and the stopped idempotent
processes $(M_{t\wedge\tau_n},\,t\in\R_+)$ are uniformly maximable
$\mathbf{A}$-exponential maxingales. 
\begin{lemma}
  \label{le:timechange}
Let $M=(M_t,\,t\in\R_+)$ be an exponential maxingale relative to 
a $\tau$-flow $\mathbf{A}=(\mathcal{A}_t,\,t\in\R_+)$ and
$\sigma_t,\,t\in\R_+$, be a collection of bounded idempotent 
$\mathbf{A}$-stopping
times such that $\sigma_s\le \sigma_t$ for $s\le t$. Then the
idempotent process $(M_{\sigma_t},\,t\in\R_+)$ is an exponential
maxingale relative to the $\tau$-flow $(\mathcal{A}_{\sigma_t},\,t\in\R_+)$.
\end{lemma}
\begin{proof}
  By Corollary 2.3.10 in Puhalskii \citeyear{Puh01},
  $\mathbf{S}(M_{\sigma_t}|\mathcal{A}_{\sigma_s})=
M_{\sigma_s}$ $\mathbf{\Pi}$-a.e. for
  $s\le t$. Each $M_{\sigma_t}$ is maximable since by the boundedness
  of $\sigma_t$ there exists $T\ge \sigma_t$, so 
$M_{\sigma_t}=\mathbf{S}(M_T|\mathcal{A}_{\sigma_t})$, which is maximable by
  maximability of $M_T$, inclusion
  $\mathcal{A}_{\sigma_t}\subset\mathcal{A}_T$, and 
Lemma 1.6.21 in Puhalskii \citeyear{Puh01}.
\end{proof}
Given an $\R$-valued function
$G=(G_t(\lambda;\bx),\,t\in\R_+,\,\bx\in\C(\R_+,\R),\,\lambda\in\R)$,
where $G_t(\lambda;\bx)$ is $\mathcal{C}_t$-measurable in $\bx$,  we say
that a deviability $\mathbf{\Pi}$ on $\C(\R_+,\R)$
solves  the maxingale problem $(x,G)$, where $x\in\R$, if 
 $X_0=x$
$\mathbf{\Pi}$-a.e. and $(\exp(\lambda X_t-G_t(\lambda;X)),\,t\in\R_+)$ is a
$\mathbf{C}$-local exponential maxingale under $\mathbf{\Pi}$, where 
$X=(X_t,\,t\in\R_+)$ is the canonical 
idempotent process on $\C(\R_+,\R)$.
We have the following lemma.
\begin{lemma}
  \label{le:expmax}Let $\mathbf{\Pi}$ solve the maxingale problem $(x,G)$.
If the function $(G_t(\lambda;\bx),\,t\in\R_+,\,\bx\in\C(\R_+,\R))$ is
bounded in $(t,\bx)$ for all $\lambda\in\R$, then the process 
$(\exp(\lambda X_t-G_t(\lambda;X)),\,t\in\R_+)$ is a
$\mathbf{C}$-uniformly maximable exponential maxingale under $\mathbf{\Pi}$.
\end{lemma}
\begin{proof}Let $M_t(\lambda)=\exp(\lambda X_t-G_t(\lambda;X))$.
  By Lemma 2.3.13(3) in Puhalskii \citeyear{Puh01} 
it is enough to prove that the
  collection $(M_t(\lambda),\,t\in\R_+)$ is
  uniformly maximable. The definition of a local exponential maxingale
  and Lemma 1.6.22 in Puhalskii \citeyear{Puh01} imply that 
$\mathbf{S}_\mathbf{\Pi} M_t(2\lambda)\le1$. 
Therefore, denoting by $b$ an upper bound
  for $(\exp(-2G_t(\lambda;\bx)),\,t\in\R_+,\,\bx\in\C(\R_+,\R))$ and 
$(\exp(G_t(2\lambda;\bx)),\,t\in\R_+,\,\bx\in\C(\R_+,\R))$, we have
\begin{equation*}
  \mathbf{S}_\mathbf{\Pi} M_t(\lambda)^2=
\mathbf{S}_\mathbf{\Pi}\bl( M_t(2\lambda)
\exp(G_t(2\lambda;X))\exp(-2 G_t(\lambda;X))\br)\le b^2.
\end{equation*}
The uniform maximability now follows by Corollary 1.4.15 
in Puhalskii \citeyear{Puh01}.
\end{proof}
\begin{lemma}
  \label{le:prodmax}
Let  $(M_t,\,t\in\R_+)$ and $(M'_t,\,t\in\R_+)$ 
be exponential maxingales on $(\Upsilon,\mathbf{\Pi})$ and
$(\Upsilon',\mathbf{\Pi}')$, respectively, relative to the respective
$\tau$-flows $(\mathcal{A}_t,\,t\in\R_+)$ and 
$(\mathcal{A}'_t,\,t\in\R_+)$. Then   $(M_tM_t',\,t\in\R_+)$
is an exponential maxingale on
$(\Upsilon\times\Upsilon',\mathbf{\Pi}\times\mathbf{\Pi}')$ 
relative to the $\tau$-flow 
$(\mathcal{A}_t\otimes \mathcal{A}'_t,\,t\in\R_+)$.
\end{lemma}
\begin{proof}
  By Puhalskii \citeyear[Lemma 1.6.28]{Puh01},
  $\mathbf{S}_{\mathbf{\Pi}\times\mathbf{\Pi}'}
(M_tM_t'|\mathcal{A}_s\otimes \mathcal{A}'_s)=
\mathbf{S}_{\mathbf{\Pi}}(M_t|\mathcal{A}_s)
\mathbf{S}_{\mathbf{\Pi}'}(M_t'| \mathcal{A}'_s)$
  $\mathbf{\Pi}\times\mathbf{\Pi}'$-a.e. for $s\le t$. Maximability of 
$(M_tM_t',\,t\in\R_+)$ under $\mathbf{\Pi}\times\mathbf{\Pi}'$ is obvious.
\end{proof}
Poisson idempotent probability 
(or Poisson deviability) is a deviability on $\C(\R_+,\R)$   defined 
by
\begin{equation*}
  \mathbf{\Pi}^N(\bx)=
  \begin{cases}
\displaystyle    \exp\Bl(-\int_0^\infty\pi(\dot{\bx}_t)
\,dt\Br)&\begin{aligned}[t]&
\text{if }\bx\text{ is absolutely continuous and non-decreasing,}\\&
\text{and } \bx_0=0,  
\end{aligned}
\\
0&\text{ otherwise.}
  \end{cases}
\end{equation*}
A  Poisson idempotent
process on $(\Upsilon,\mathbf{\Pi})$ 
is defined as an idempotent process with idempotent
distribution  $\mathbf{\Pi}^N$. Thus, a Poisson
idempotent process has   absolutely continuous non-decreasing
trajectories $\mathbf{\Pi}$-a.e. 
  The
definition implies that the canonical idempotent process on $\C(\R_+,\R)$
is  Poisson under  $\mathbf{\Pi}^N$. If $N$ is a Poisson idempotent process on
$(\Upsilon, \mathbf{\Pi})$, then the idempotent process 
$M^N(\lambda)=(M^N_t(\lambda),\,t\in\R_+)$ defined by 
$M^N_t(\lambda)=\exp\bl(\lambda N_t-(e^\lambda-1)t\br)$ is an
exponential maxingale relative to the $\tau$-flow
$(\mathcal{A}^N_t,\,t\in\R_+)$, where the $\mathcal{A}^N_t$ are
the $\tau$-algebras 
generated by the $N_s,\,s\le t$ 
Puhalskii \citeyear[Theorem 2.4.16]{Puh01}. We say
that a continuous-path idempotent process $N$ is Poisson relative to a
$\tau$-flow $\mathbf{A}$ if $N_0=0$ and the idempotent process
$M^N(\lambda)$ is an $\mathbf{A}$-exponential maxingale for all
$\lambda\in\R$. If $N$  is idempotent Poisson relative to 
 $\mathbf{A}$, then it is idempotent Poisson,
 Puhalskii \citeyear[Corollary 2.4.19]{Puh01}.

  Wiener idempotent probability 
(or Wiener deviability) is a deviability on $\C(\R_+,\R)$   defined 
by
\begin{equation*}
  \mathbf{\Pi}^W(\bx)=
  \begin{cases}
\displaystyle    \exp\Bl(-\frac{1}{2}\int_0^\infty\dot{\bx}_t^2\,dt\Br)&
\text{ if }\bx\text{ is absolutely continuous and } \bx_0=0,\\
0&\text{ otherwise}.
  \end{cases}
\end{equation*}
A Wiener idempotent process  on $(\Upsilon,\mathbf{\Pi})$ is defined as 
an idempotent process with idempotent
distribution $\mathbf{\Pi}^W$. Thus, a Wiener idempotent  process has 
$\mathbf{\Pi}$-a.e. absolutely continuous paths.
 The
definition implies that the canonical idempotent process on $\C(\R_+,\R)$
is Wiener  under $\mathbf{\Pi}^W$.

Let $W=(W_t,\,t\in\R_+)$ 
be a Wiener idempotent
process on $(\Upsilon,\mathbf{\Pi})$. Then the idempotent process 
$\bl(\exp(\lambda W_t-\lambda^2t/2),\,t\in\R_+\br)$ is an exponential
maxingale relative to the flow
$\mathbf{A}^W=(\mathcal{A}^W_t,\,t\in\R_+)$, where the
$\mathcal{A}^W_t$ are the $\tau$-algebras generated by $W_s,\,s\le t$,
Puhalskii \citeyear[Theorem 2.4.2]{Puh01}. 
We say
that a continuous-path idempotent process $W$ is Wiener relative to a
$\tau$-flow $\mathbf{A}$ if $W_0=0$ and the idempotent process
$\bl(\exp(\lambda W_t-\lambda^2t/2),\,t\in\R_+\br)$
 is an $\mathbf{A}$-exponential maxingale for all
$\lambda\in\R$. If $W$  is idempotent Wiener relative to 
 $\mathbf{A}$, then it is idempotent Wiener,
 Puhalskii \citeyear[Corollary 2.4.6]{Puh01}.
In particular, $W_t-W_s$ for $t\ge s$ is
independent of $\mathcal{A}_s$ by the fact that 
$\mathbf{S}\bl(\exp(\lambda (W_t-W_s))|\mathcal{A}_s\br)
=\exp(\lambda^2(t-s)/2)$, which is a
smooth function of $\lambda$.

 Given a bounded $\R$-valued idempotent process
$\sigma_t,\,t\in\R_+$, we define the idempotent Ito integral 
$(\sigma\diamond W)_t$ by 
\begin{equation*}
  (\sigma\diamond W)_t(\upsilon)=
  \begin{cases}
\displaystyle    
\int_0^t\sigma_s(\upsilon)
\dot{W}_s(\upsilon)\,ds&\text{ if }\mathbf{\Pi}(\upsilon)>0,\\
Y(\upsilon)&\text{ otherwise },
  \end{cases}
\end{equation*}
where $Y(\upsilon)$ is an $\R$-valued idempotent variable and
$\dot{W}_s(\upsilon)$ denotes the Radon-Nikodym derivative in $s$ of the Wiener
idempotent trajectory. The integral is thus specified uniquely
$\mathbf{\Pi}$-a.e. 
The idempotent process 
$\bl((\sigma\diamond W)_t,\,t\in\R_+\br)$ has $\mathbf{\Pi}$-a.e.
continuous paths. If $(W_t,\,t\in\R_+)$ and $(\sigma_t,\,t\in\R_+)$ are
adapted to a complete $\tau$-flow $\mathbf{A}$,
then $((\sigma\diamond W)_t,\,t\in\R_+)$ is $\mathbf{A}$-adapted.
For clarity, we further use  $\int_0^t\sigma_s\dot{W}_s\,ds$ for
 $(\sigma\diamond W)_t$. In the next lemma, 
$\int_s^t \sigma_p\dot{W}_p\,dp=\int_0^t\sigma_p\ind(r\in[s,t])\dot{W}_p\,dp$.
\begin{lemma}
  \label{le:wieind}
Let $\sigma_s,\,s\in\R_+$ be an $\R$-valued bounded
Lebesgue-measurable function and $W=(W_t,\,t\in\R_+)$ 
be a Wiener idempotent
process on $(\Upsilon,\mathbf{\Pi})$ relative to a complete $\tau$-flow
$\mathbf{A}$. Then the
idempotent process $M=(M_t,\,t\in\R_+)$, where
$M_t=\exp(\lambda\int_0^t \sigma_s\dot{W}_s\,ds-
\lambda^2\int_0^t\sigma_s^2\,ds/2)$,
is an $\mathbf{A}$-exponential maxingale. In
particular, $\int_s^t \sigma_p\dot{W}_p\,dp$ is 
independent of $\mathcal{A}_s$ for $s\le t$.
\end{lemma}
\begin{proof}
The idempotent process $M$ is $\mathbf{A}$-adapted by
$M_t$ being constant on the atoms of $\mathcal{A}_t$ for $t\in\R_+$,
cf., Puhalskii \citeyear[Lemma 2.2.17]{Puh01}. If
the function $\sigma_s,\,s\in\R_+$, is piecewise constant, the
maxingale property follows by the properties of conditional idempotent
expectations in a standard manner. A limit argument shows that this
property carries over to continuous $\sigma_s,\,s\in\R_+$. 
The case of a Lebesgue measurable $\sigma_s,\,s\in\R_+$, follows
via Luzin's theorem. Maximability of the $M_t$ follows by
Lemma~\ref{le:expmax}. Finally, $\int_s^t \sigma_p\dot{W}_p\,dp$ is 
independent of $\mathcal{A}_s$ for $s\le t$ since by the
maxingale property $\mathbf{S}\bl(\exp(\lambda\int_s^t
\sigma_p\dot{W}_p\,dp)|\mathcal{A}_s\br)=
\exp\bl((\lambda^2/2) \int_s^t
\sigma_p^2\,dp\br)$, where  the latter is a smooth function of $\lambda$.
\end{proof}
Let $\sigma_t(x),\,x\in\R,\,t\in\R_+,$ and
$b_t(x),\,x\in\R,\,t\in\R_+,$ 
be real-valued  
functions, which are continuous in $x$ and Lebesgue-measurable in $t$. 
Let $W$ be a Wiener idempotent process on an
idempotent probability space $(\Upsilon,\mathbf{\Pi})$ relative to a complete
$\tau$-flow $\mathbf{A}$ and let $\overline{\mathcal{C}}^W_t$ for
$t\in\R_+$ denote the
completion of $\mathcal{C}_t$ with respect to the Wiener deviability
on $\C(\R_+,\R)$. We say that, given $x\in\R$,
 an idempotent process $X$ on
$(\Upsilon,\mathbf{\Pi})$ is  a strong solution to the Ito idempotent equation 
  \begin{equation}
    \label{eq:46}
  X_t=x+\int_0^tb_s(X_s)\,ds+\int_0^t\sigma_s(X_s)\,\dot{W}_s\,ds,\;t\in\R_+,
\end{equation}
where integrals are understood as Lebesgue integrals,
if equality \eqref{eq:46} holds $\mathbf{\Pi}$-a.e. and there exists a
  function  $J:\,\C(\R_+,\R)\to\C(\R_+,\R)$, which is 
$\overline{\mathcal{C}}^W_t/\mathcal{C}_t$-measurable  for every
$t\in\R_+$, such that  $X=J(W)$ $\mathbf{\Pi}$-a.e.  
 As a consequence, $X$ is
$\mathbf{A}$-adapted.   A strong solution is called
Luzin if the function $J$ is continuous in restriction to the sets 
$\{\bx\in\C(\R_+,\R):\,\mathbf{\Pi}^W(\bx)\ge a\}$ for $a\in(0,1]$.
We say that there exists a unique  strong solution 
(respectively, Luzin strong solution) 
if any strong solution (respectively, Luzin strong solution) can be
written as $X=J(W)$ $\mathbf{\Pi}$-a.e.  for the same function $J$.
Let us assume that 
$\sigma_t(x)$ and $b_t(x)$ are locally Lipshitz-continuous in $x$,
i.e., for every $a>0$ there
 exists an $\R_+$-valued Lebesgue-measurable in $t$
function $k^a_t,\,t\in\R_+,$ 
with $\int_0^tk^a_s\,ds<\infty$ for $t\in\R_+$ such
that $\abs{b_t(x)-b_t(y)}\le k^a_t\abs{x-y}$ and 
$\abs{\sigma_t(x)-\sigma_t(y)}^2\le k^a_t\abs{x-y}^2$ if $\abs{x}\le
a$ and $\abs{y}\le a$,
 and satisfy the linear-growth condition that
there exists an $\R_+$-valued 
Lebesgue-measurable function $l_t,\,t\in\R_+,$ 
with $\int_0^tl_s\,ds<\infty$ for
$t\in\R_+$ such that
$\abs{b_t(x)}\le l_t(1+\abs{x})$ and 
$\sigma_t(x)^2\le l_t(1+\abs{x}^2)$ for $x\in\R$. Then
\eqref{eq:46} has a unique strong solution, which is also a Luzin
strong solution,
 Puhalskii \citeyear[Theorems 2.6.21, 2.6.22 and 2.6.26]{Puh01}.

Let $\Upsilon$ be a metric space.
A net $\mathbf{\Pi}^\psi,\,\psi\in\Psi,$ where $\Psi$ is a directed set,
 of idempotent probabilities on $\Upsilon$
is said to converge weakly to idempotent probability $\mathbf{\Pi}$ on
$\Upsilon$ if 
$\lim_{\psi\in\Psi}\mathbf{S}_{\mathbf{\Pi}^\psi}f
=\mathbf{S}_\mathbf{\Pi} f$ for every
non-negative bounded and continuous function $f$ on $\Upsilon$;
equivalently,  Puhalskii \citeyear[Theorem 1.9.2]{Puh01},
$\limsup_{\psi\in\Psi}\mathbf{\Pi}^\psi (F)  \leq \mathbf{\Pi}(F)$
for all closed sets $F  \subset \Upsilon$ and 
$\liminf_{\psi\in\Psi} \mathbf{\Pi}^\psi (G)  \geq  \mathbf{\Pi}(G)$
 for all open sets $G  \subset \Upsilon$.  
A net of idempotent variables with values in the same metric
space is said to converge in idempotent distribution if their
idempotent distributions weakly converge. One has a continuous mapping
theorem for convergence in idempotent distribution: if a net
$X^\psi,\,\psi\in\Psi,$
of idempotent variables with values in $\Upsilon$
converges in idempotent distribution to an idempotent variable $X$
with values in $\Upsilon$ and $f$ is a continuous function from
$\Upsilon$ to a metric space $\Upsilon'$, then the net
$f(X^\psi),\,\psi\in\Psi,$ 
converges in idempotent distribution to $f(X)$.
 A  net $\mathbf{\Pi}^\psi,\,\psi\in\Psi$, of deviabilities on 
$\Upsilon$ is said to be
 tight  if
$\inf_{K\in\mathcal{K}}\limsup_{\psi\in\Psi}\mathbf{\Pi}^\psi(\Upsilon\setminus
K)=0$, where $\mathcal{K}$ denotes the collection of compact subsets
of $\Upsilon$. A tight net of deviabilities 
is weakly  compact, i.e., 
it contains  a subnet that converges weakly to a deviability, see
Puhalskii \citeyear[Theorem 1.9.27]{Puh01}
(if $\mathbf{\Pi}^\psi$ is a sequence, then it
contains a weakly convergent subsequence).
\bibliographystyle{annprob}
\bibliography{large,puh,stoch,rgraph,other,que}
\vspace{2.cm}
Anatolii Puhalskii, Mathematics Department, 
University of Colorado at Denver, Campus Box 170, P.O. Box 173364, 
Denver, CO 80217-3364, U.S.A.
 and Institute for Problems in Information
Transmission, Moscow, Russia. Phone: (1)-303-5564811. Email:
puhalski@math.cudenver.edu
\end{document}